\documentstyle[amssymb,10pt]{amsart}
\newtheorem{theorem}{Theorem}[section]
\newtheorem{lemma}[theorem]{Lemma}
\newtheorem{proposition}[theorem]{Proposition}
\newtheorem{corollary}[theorem]{Corollary} 
\theoremstyle{definition}  
\newtheorem{definition}[theorem]{Definition}

\newtheorem{example}[theorem]{Example}

\newtheorem{remark}[theorem]{Remark}

\newcommand{\id}{\operatorname{id}}

\newcommand{\End}{\operatorname{End}} 

\newcommand{\Hom}{\operatorname{Hom}} 
\newcommand{\Ad}{\operatorname{Ad}}

\newcommand{\Ind}{\text{Ind}}

\newcommand{\nc}{\newcommand}
\nc{\Symm}{{\on{Sym}}}

\newcommand{\on}{\operatorname}   
\newcommand{\eps}{\varepsilon}
 \nc{\cE}{{\cal E}}

\renewcommand{\a}{{\mathfrak a}}
\renewcommand{\b}{{\mathfrak b}}
\renewcommand{\c}{{\mathfrak c}}

\nc{\D}{{\mathfrak d}}
\nc{\SL}{{\mathfrak sl}}
\newcommand{\G}{{\mathfrak{g}}}\nc{\HH}{{\mathfrak h}}
\newcommand{\g}{{\mathfrak{g}}}
\newcommand{\x}{{\mathfrak{x}}}
\renewcommand{\u}{{\mathfrak{u}}}
\newcommand{\p}{{\mathfrak{p}}}
\newcommand{\vv}{{\mathfrak{v}}}

\newcommand{\z}{{\mathfrak{z}}}
\newcommand{\h}{{\mathfrak{h}}}
\renewcommand{\l}{{\mathfrak{l}}}
\renewcommand{\k}{{\mathfrak{k}}}
\newcommand{\n}{{\mathfrak{n}}}
\newcommand{\m}{{\mathfrak{m}}}
\renewcommand{\b}{{\mathfrak{b}}}
\renewcommand{\t}{{\mathfrak{t}}}

\nc{\wh}{\widehat}\nc{\wt}{\widetilde}

\newcommand{\la}{{\lambda}}

\newcommand{\ben}{\begin{enumerate}}
\newcommand{\een}{\end{enumerate}}
\newcommand{\ad}{{\text{ad}}}

\newcommand{\cO}{{\mathcal O}}
\newcommand{\cT}{{\mathcal T}}

\newcommand{\CC}{{\mathbb{C}}}

\newcommand{\RR}{{\mathbb{R}}}
\newcommand{\PP}{{\mathbb{P}}}

\newcommand{\al}{\alpha}
\newcommand{\ZZ}{{\mathbb{Z}}}

\newcommand{\cA}{{\mathcal A}}

\begin{document}

\title{Quantization of classical dynamical $r$-matrices with nonabelian base}

\begin{abstract} 
We construct some classes of dynamical $r$-matrices over a nonabelian base, 
and quantize some of them by constructing dynamical (pseudo)twists in the sense
of Xu.  This way, we obtain quantizations of $r$-matrices obtained in 
earlier work of the second author with Schiffmann and Varchenko. A part of 
our construction may be viewed as a generalization of the Donin-Mudrov 
nonabelian fusion construction. We apply these results to the construction 
of equivariant star-products on Poisson homogeneous spaces, which include 
some homogeneous spaces introduced by De Concini. 
\end{abstract}

\author{Benjamin Enriquez}
\address{IRMA (CNRS), rue Ren\'e Descartes, F-67084 Strasbourg, France}
\email{enriquez@@math.u-strasbg.fr}

\author{Pavel Etingof}
\address{Department of Mathematics, Massachusetts Institute of Technology,
Cambridge, MA 02139, USA}
\email{etingof@@math.mit.edu}

\maketitle

\section*{Introduction}

In this paper, we construct generalizations of some classes of classical 
dynamical
$r$-matrices with nonabelian base, and quantizations of some of them. 
This way, we obtain quantizations of dynamical $r$-matrices introduced 
in \cite{EV1,ES2}. We then apply these results to obtain explicit, 
equivariant star-products on some homogeneous spaces. In particular, 
we obtain quantizations of Poisson homogeneous spaces,
introduced by De Concini.  

The classes of $r$-matrices we consider are the following
(a), ..., (d) (all Lie algebras are assumed to be finite dimensional). 

\medskip 

(a) Let $\g = \l\oplus\u$ be a Lie algebra with a nondegenerate splitting
(see Section \ref{rat:class}). Then the natural map $\u\otimes\u \to\l$ 
can be "inverted" and yields a solution $r_\l^\g(\lambda)  
\in \big( \wedge^2(\g)
\otimes S^\cdot(\l)[1/D_0] \big)^\l$ of the classical dynamical 
Yang-Baxter equation (CDYBE)
$$
\on{CYB}(r_\l^\g(\lambda)) - \on{Alt}(\on{d} r_\l^\g(\lambda)) = 0. 
$$
$r_\l^\g(\lambda)$ is a rational function in $\lambda$, homogeneous of degree $-1$, 
and is a generalization of the rational classical dynamical $r$-matrices
of \cite{EV1}. 

$r_\l^\g(\lambda)$ also plays a role in "composing" $r$-matrices: 
\begin{proposition} \label{prop:comp:r:matrices} \label{prop:comp}
(see \cite{EV1}, Theorem 3.14 and \cite{FGP}, Proposition 1.)
Let $\l\subset \g$ be an inclusion of Lie algebras. 
Assume that $\l$ has a nondegenerate splitting $\l = \k\oplus\m$.  
Given $Z\in\wedge^3(\g)^\g$, let us say that a $(\l,\g,Z)$-$r$-matrix 
is an $\l$-invariant function $\rho : \l^* \to \wedge^2(\g)$, 
solution of $\on{CYB}(\rho(\lambda)) - \on{Alt}(\on{d}\rho(\lambda)) = Z$. 
Set 
$$
\sigma(\lambda) := r_\k^\l(\lambda) + \rho_{|\k^*}(\lambda).   
$$
Then $\rho\mapsto\sigma$ is a map $\{(\g,\l,Z)$-$r$-matrices$\}
\to \{(\g,\k,Z)$-$r$-matrices$\}$. 
\end{proposition}

Here by a function $\l^*\to \wedge^2(\g)$, we understand an element of
$\wedge^2(\g) \otimes \wh S^\cdot(\l)[1/\Delta]$, where $\Delta$ 
is a suitable nonzero element of $\wh S^\cdot(\l)$. 
We set $\on{CYB}(\rho): = [\rho^{1,2},\rho^{1,3}]
+ [\rho^{1,2},\rho^{2,3}] + [\rho^{1,3},\rho^{2,3}]$, 
and $\on{d}(x\otimes y \otimes z_1\cdots z_l) := \sum_{i=1}^l 
x\otimes y \otimes z_i \otimes z_1\cdots \check z_i \cdots z_l$. 
If $\xi\in\g^{\otimes 3}$ is antisymmetric in two tensor factors, we set
$\on{Alt}(\xi\otimes f) = (\xi + \xi^{2,3,1} + \xi^{3,1,2})\otimes f$. 

Proposition \ref{prop:comp:r:matrices} enables us to construct 
new $r$-matrices from known ones. 

\medskip 

(b) Let $(\g = \l \oplus \u,t\in S^2(\g)^\g)$ be a quadratic Lie 
algebra with a nondegenerate splitting (we do not assume $t$ to be 
nondegenerate).  
We may apply Proposition \ref{prop:comp:r:matrices} to $\rho := $
the Alekseev-Meinrenken $r$-matrix of $\g$ (\cite{AM1}), 
$(\l,\k,\m) := (\g,\l,\u)$, and obtain: 

\begin{corollary} \label{cor:rat:trigo} \label{cor:1:2}
(see also \cite{FGP}.)
Let $c\in\CC$, and $\lambda\in\l^*$, set 
$$
\rho_c(\lambda) = r_\l^\g(\lambda) 
+ c \big( f(c \on{ad}(\lambda^\vee))\otimes \id\big)(t). 
$$
Then we have $\on{CYB}(\rho_c) - \on{Alt}(\on{d}\rho_c) = -\pi^2 c^2 Z$, 
where $Z = [t^{1,2},t^{2,3}]$.
\end{corollary}
Here we set $\lambda^\vee = (\lambda\otimes \id)(t)$ and $f(x) = -1/x + 
\pi\on{cotan}(\pi x)$.

\medskip 
(c) Let $(\g,t\in S^2(\g)^\g)$ be a quadratic Lie algebra, equipped with 
$\sigma \in \on{Aut}(\g,t)$. We assume that $\sigma - \id$ is invertible on 
$\g / \g^\sigma$. We set $\l := \g^\sigma$, $\u := \on{Im}(\sigma - \id)$, 
so $\g = \l\oplus\u$, $t = t_\l + t_\u$, 
$t_\x \in S^2(\x)$ for $\x = \l,\u$. The following result can be found in 
\cite{AM2} (see also \cite{S} and \cite{ES2}, Theorem A1). 
\begin{proposition} \label{lemma:sigma} \label{lemma:0:3}
Set 
$$
\rho_{\sigma,c}(\lambda) := (c f(c\on{ad}(\lambda^\vee)) \otimes \id)(t_\l)
+ i\pi c ({{e^{2\pi ic\on{ad}(\lambda^\vee)} \circ \sigma + \id}\over{
e^{2\pi ic\on{ad}(\lambda^\vee)} \circ \sigma - \id}}\otimes \id)(t_\u) . 
$$ 
Here we set 
$\lambda^\vee = (\lambda\otimes \id)(t_\l)$ for $\lambda\in\l^*$. 
Then $\rho_{\sigma,c}$ is a solution of $\on{CYB}(\rho_{\sigma,c}) 
- \on{Alt}(\on{d}\rho_{\sigma,c}) = -\pi^2 c^2 Z$. 
\end{proposition}

Note that if $\chi : \l\to\CC$ is a character, then 
$\chi^\vee$ is central in $\l$, and if $\g = \l\oplus \u$
is nondegenerate and $t_\u$ is nondegenerate, then 
$\rho_{\on{exp}(\on{ad}(\chi^\vee))}(\lambda)$
coincides with $\rho_c(\lambda - \chi)$, with $\rho_c$ as in Corollary
\ref{cor:1:2}. 

If now $\l$ has a nondegenerate splitting $\l = \k\oplus \m$, 
then Proposition \ref{prop:comp:r:matrices} implies that 
$r_\k^\l + (\rho_{\sigma,c})_{|\k^*}$ is a 
$(\k,\g,-4\pi^2c^2 Z)$-$r$-matrix.  

\medskip 
(d) Let $\g = \l\oplus \u$ be a Lie algebra with a splitting. 
Assume that $t\in S^2(\g)^\g$ decomposes as $t_\l + t_\u$, with 
$t_\x\in S^2(\x)$ for $\x = \l,\u$. Let us say that $C\in \End(\u)$
is a Cayley endomorphism if it satisfies the following axioms:
$C$ is a $\l$-module endomorphism, and for any $x,y\in\u$, we have 
$$
[C(x),C(y)]_\u = C([C(x),y]_\u) + C([x,C(y)]_\u) - [x,y]_\u.  
$$
Here for $x\in \g$, we denote by $x_\u$ its projection on
$\u$ parallel to $\l$. 
If $\sigma$ is as in (3), then 
$C = C(\sigma) := (\sigma + \id) / (\sigma - \id)$ 
is a Cayley endomorphism. (Such Cayley endomorphisms are exactly those 
which do not contain $\pm 1$ in their spectrum.) More generally, a limit of 
$C(\exp(\on{ad}(\chi)))$, where some eigenvalues of $\chi$ tend to 
$\pm\infty$, is a Cayley endomorphism. 

\begin{proposition} \label{prop:C} \label{prop:0:4}
Assume that $(C\otimes \id + \id \otimes C)(t_\u) = 0$ (when $C = C(\sigma)$, 
this condition means that $\sigma$ preserves $t_\u$). Set 
$$
\rho_{C,c}(\lambda) := 
(c f(c\on{ad}(\lambda^\vee)) \otimes \id)(t_\l)
+ i\pi c
\Big(  
{{C + i \on{tan}(\pi c \on{ad}(\lambda^\vee))}\over
{1 + i \on{tan}(\pi c \on{ad}(\lambda^\vee))C}}
\otimes \id \Big) (t_\u). 
$$
Then $\rho_{C,c}$ is a $(\l,\g,-\pi^2c^2 Z)$-$r$-matrix. 
\end{proposition}

\medskip 

We define quantizations of solutions of the (modified) CDYBE as 
solutions of suitable (pseudo)twist equations (Sections \ref{sect:quant}, 
\ref{sect:pseudo}, and also \cite{EE}, equation 9). We construct 
quantizations for the above $r$-matrices in the following cases. 

\medskip 

(a') Rational $r$-matrices. We construct quantizations $J_\l^\g$ 
of the rational
$r$-matrices introduced above in the particular case when $\g$ is 
{\it polarized}, i.e.,  $\u$ decomposes as a sum of $\l$-submodules
$\u_+\oplus \u_-$, such that $\u_\pm$ are Lie subalgebras of $\g$
(Corollary \ref{cor:quant:rat}). 
We do so by working out a nonabelian generalization of
the fusion construction of \cite{EV2} (whose ideas originate 
from \cite{Fad,AF}). Recently J. Donin and A. Mudrov \cite{DM} (see also
\cite{AL}) extended 
this construction to the case when $\h$ is replaced
with a Levi subalgebra $\l\subset \g$; their work relies on 
Jantzen's computation of the 
Shapovalov form for induced modules. To generalize their result, 
we work directly in (microlocalizations of) universal enveloping 
algebras. 

\medskip 

(b') We quantize the rational-trigonometric $r$-matrices of Corollary
\ref{cor:rat:trigo} in the following situation: $\g$ is 
polarized, and $t\in S^2(\g)^\g$ decomposes as $t_\l + s + s^{2,1}$, 
where $t_\l\in S^2(\l)$ and $s\in \u_+ \otimes \u_-$ (Theorem \ref{thm:barK}). 
Our argument is based on nonabelian versions of the ABRR identities 
(see \cite{ABRR,EV1,ES2}), which are satisfied by $J_\l^\g$ when 
$\g$ is polarized and quadratic, and the use of Drinfeld associators
(\cite{Dr2}). When $\l = \g$, our construction coincides with the quantization
of the Alekseev-Meinrenken $r$-matrix (\cite{EE}), which is based on
renormalizing an associator.

\medskip 

(c') We quantize the $r$-matrix $\rho_{\sigma,c}$ defined in 
Proposition \ref{lemma:sigma} (Theorem \ref{thm:pseudo}). 
We also quantize the $r$-matrix 
$(\rho_{\sigma,c})_{|\k^*} + r_\k^\l$ under the assumption that 
$\l = \k\oplus \m_+ \oplus \m_-$ is quadratic polarized (Proposition 
\ref{prop:eta}). For this, 
we introduce a compatible differential system, generalizing the
Knizhnik-Zamolodchikov (KZ) system, and we adapt Drinfeld's proof that 
the KZ associator satisfies the pentagon equations (\cite{Dr2}). 

In Section \ref{sect:5:6}, we explain why the quantizations
obtained in (b'), (c') may be interpreted in terms of 
infinite-dimensional ABRR equations for extended (twisted) 
loop algebras.

\medskip 

In Section \ref{sect:homog}, we apply these results to the construction of
equivariant star-products on Poisson homogeneous spaces. In particular, 
we quantize Poisson homogeneous spaces introduced by De Concini. 
\medskip

\begin{remark}
An earlier version of this paper is available
at www-math.mit.edu/~~${}_{\widetilde{}}$~~~~~etingof/ee.tex; this version is less general
but uses a somewhat more intuitive representation theoretic language.
\end{remark}

\noindent
{\bf Acknowledgments.} The authors thank J. Donin and 
A. Mudrov for references. P.E. is grateful 
to C. De Concini for a very useful discussion
on Poisson homogeneous spaces related to automorphisms, 
which was crucial for writing Section \ref{sect:homog} of this paper. 
P.E. is indebted to IRMA (Strasbourg) for hospitality. 
The work of P.E. was partially supported by the NSF grant DMS-9988796.

\newpage

\section{Rational classical dynamical $r$-matrices}
\label{rat:class} \label{sect:rat}

In this section, we introduce the notion of a (nondegenerate) Lie algebra with
a splitting $\g = \l\oplus \u$. We associate to each such nondegenerate Lie
algebra a rational $r$-matrix $r_\l^\g$. We show that in the case of a double
inclusion $\k\subset \l\subset\g$ of Lie algebras, $r_\k^\l$ plays a role in a
restriction theorem for $r$-matrices. We introduce the notions of polarized Lie
algebras and of quantizations of the rational $r$-matrices $r_\l^\g$. As
$r_\l^\g$ is singular at $\lambda = 0\in \l^*$, the latter notion involves
a microlocalization of $U(\l)$ (in the sense of \cite{Spr}). 

\medskip
\noindent {\it Notation.} If $A$ is a Hopf algebra, we use Sweedler's notation:
$\Delta(x) = \sum x^{(1)} \otimes x^{(2)}$. 
If $x\in A$, we write $x^{(2)} := 1 \otimes x \otimes 1\cdots\in A^{\otimes n}$, 
and if $x = \sum_i x'_i \otimes x''_i\in A^{\otimes 2}$, 
we set $x^{3,2} := \sum_i 1 \otimes x''_i \otimes x'_i \otimes 1\cdots \in
A^{\otimes n}$, $x^{3,21} := \sum_i (x''_i)^{(2)} \otimes (x''_i)^{(1)}
\otimes x'_i \otimes 1\cdots \in A^{\otimes n}$, etc.

\subsection{A family of classical dynamical $r$-matrices} \label{splitted}

Let $\g$ be a finite dimensional Lie algebra. Assume that we have a 
decomposition $\g = \l \oplus \u$, where $\u$ is an $\l$-invariant 
complement of $\l$ in $\g$; that is, $[\l,\u] \subset \u$. 
Such a triple $(\g,\l,\u)$ is called a "Lie algebra with a splitting". 

We have a linear map $\l^* \to \wedge^2(\u)^*$, taking 
$\lambda\in \l^*$ to $\omega(\la) : x\wedge y \mapsto \lambda([x,y])$. 
The triple $(\g,\l,\u)$ is called {\it nondegenerate} 
if for a generic $\lambda\in \l^*$, $\omega(\lambda)$ is nondegenerate. 
The algebraic translation of this condition is the following: identify
$\wedge^2(\u)^*$ with a subspace of $\End(\u)$ using any 
linear isomorphism $\u \simeq \u^*$, then the map 
$\lambda \mapsto \det \omega(\la)$ does not vanish identically. 
This map is a degree $d := \dim(\u)$ 
polynomial on $\l^*$, i.e., an element of $S^d(\l)$. If $(\g,\l,\u)$
is nondegenerate, then $d$ is even. 

If $E$ is an even dimensional vector space, denote by 
$\wedge^2(E)_{\on{nondeg}}$ the space of nondegenerate tensors of $\wedge^2(E)$. 
Then we have a bijection $\wedge^2(E^*)_{\on{nondeg}} \to \wedge^2(E)_{\on{
nondeg}}$, $\omega\mapsto \omega^{-1}$, taking a tensor $\omega$ to its
image under the inverse of the linear isomorphism $E^*\to E$
induced by $\omega$. 

\begin{proposition}  \label{prop:r:mat} \label{prop:2:1} 
(see \cite{FGP}, Proposition 1 and \cite{Xu}, Theorem 2.3.)
Let $(\g,\l,\u)$ be a nondegenerate Lie algebra with a splitting. 
Then we have a rational map 
$$
r_\l^\g : \l^* \supset U \to \wedge^2(\u), 
$$ 
defined by $r_\l^\g(\la) := -\omega(\la)^{-1}$. 
It is homogeneous of total degree $-1$ in $\lambda$. Here $U\subset \l^*$
is the $\l$-invariant open subset $\{ \la | \det \omega(\la) \neq 0\}
\subset \l^*$. 

Then $r_\l^\g$ is $\l$-invariant, and is a solution of the CDYBE, 
i.e., $\on{CYB}(r_\l^\g) + \on{Alt}(\on{d}r_\l^\g) = 0$. 
\end{proposition}

{\em Proof.} Set $D_0 := \det \omega(\la) \in S^d(\l)$. Then 
$r_\l^\g = \sum_i u_i \otimes v_i \otimes \ell_i$ belongs to 
$\wedge^2(\u) \otimes S^\cdot(\l)[1/D_0]$, and is uniquely determined by the 
equivalent conditions 
$$
\forall x,y\in \u, \; 
\sum_i h(x,u_i) h(v_i,y) \ell_i = - h(x,y)
\quad \on{(equality\ in\ }S^\cdot(\l)[1/D_0])
$$
or 
\begin{equation} \label{repr} \label{4:years}
\forall x\in \g, \; \sum_i h(x,u_i) \ell_i \otimes v_i = - 1 \otimes x_\u
\quad \on{(equality\ in\ }S^\cdot(\l)[1/D_0] \otimes \u). 
\end{equation}
Here we denote by $x_\u,x_\l$ the components of $x\in\g$ in $\u,\l$, 
and by $h : \g\times \g \to \l$ the map $(x,y)\mapsto [x,y]_\l$. 

We define a bilinear map $\langle -,- \rangle : 
(\g^{\otimes 3} \otimes S^\cdot(\l)[1/D_0]) \times 
\G^{\otimes 3} \to S^\cdot(\l)[1/D_0]$ by 
$\langle a\otimes b \otimes c \otimes \ell, x \otimes y \otimes z \rangle 
:= h(x,a)h(y,b)h(z,c)\ell$. This pairing is left-nondegenerate, so 
we will prove that the pairing of 
$\on{CYB}(r_\l^\g) - \on{Alt}(\on{d}r_\l^\g)$
with $x\otimes y\otimes z\in \g^{\otimes 3}$ is zero. 

We have $\langle [(r_\l^\g)^{1,2},(r_\l^\g)^{1,3}], x\otimes y \otimes z
\rangle = h([y_\u,z_\u],x)$, therefore $\langle 
\on{CYB(r_\l^\g)}, x\otimes y \otimes z\rangle = 
[[y_\u,z_\u],x]_\l + \on{c.p.} = [h(y,z),x_\l] + \on{c.p.}$; 
the last equality follows from the Jacobi identity and the 
$\l$-invariance of $\u$. 

On the other hand,  differentiating (\ref{repr}), and pairing the 
resulting identity with $z\otimes y$, we get 
\begin{equation} \label{part}
\sum_i h(x,u_i) h(v_i,y)h(\eps,z) \ell_i^\eps
+ \sum_i [h(x,u_i),z_\l] \ell_i h(v_i,y) = 0. 
\end{equation}   
Here we set $\on{d}(\ell_i) = \sum \eps \otimes \ell_i^\eps$. 

We have $\langle \on{d}r_\l^\g, x\otimes y \otimes z\rangle
= - \sum_{i,\eps} h(x,u_i)h(y,v_i)h(z,\eps)\ell_i^\eps$, 
so by (\ref{part}) this is equal to 
$$
\sum_i \ell_i h(v_i,y) [h(x,u_i),z_\l].
$$ 
Now (\ref{repr}) 
implies that this is equal to $[h(x,y_\u),z_\l]$. Finally, 
$\langle \on{Alt}(\on{d}r_\l^\g), x\otimes y \otimes z \rangle = 
[h(x,y_\u),z_\l] + \on{c.p.} = [h(x,y),z_\l] + \on{c.p.}$, 
where the last equality follows from $\l$-invariance of $\u$ and the 
Jacobi identity. Finally, we get $\langle \on{CYB}(r_\l^\g) 
- \on{Alt}(\on{d}r_\l^\g), x\otimes y \otimes z \rangle = 0$, as wanted.  
\hfill \qed \medskip

\begin{remark} \label{rem12}
The nondegeneracy condition  means that for a generic $\lambda\in\l^*$, 
the tangent space $T_\la({\cal O}_\lambda)$ of the coadjoint orbit of $\lambda$ 
contains $\u^*$; this means that a generic element of $\G^*$ is conjugate to 
an element of $\l^*$. 
\end{remark}

\begin{remark}
$D_0$ satisfies $\on{ad}(a)(D_0) = \chi_0(a) D_0$, 
where $\chi_0 : \l \to \CC$ is the character of $\l$ defined by
$\chi_0(a) = \on{tr} (\on{ad}(a)_{\u})$, so $D_0$ is $\l$-equivariant. 
\end{remark}

\subsection{Composition of $r$-matrices}

Let us prove Proposition \ref{prop:comp:r:matrices}.  

Let us first prove that the restriction $\rho_{|\k^*}$
is well-defined. The singular locus $\{\lambda\in\l^* | \Delta(\lambda) = 0\}$
is $\l$-invariant, so it cannot contain $\k^*$; therefore $\Delta_{|\k^*}$
is nonzero, and $\rho_{|\k^*} \in \wedge^2(\g) \otimes \wh
S^\cdot(\k)[1/\Delta_{|\k^*}]$ is well-defined. Both $r_\k^\l$
and $\rho_{|\k^*}$ are $\k$-invariant, hence so is $\sigma$.

Let us write  
$r_\k^\l(\lambda)$ as $\sum_i u_i(\lambda) \otimes e^i$, 
where $u_i(\lambda)\in \u\otimes S^\cdot(\l)[1/D_0]$, and show that  
\begin{equation} \label{step:30}
\on{ad}^*(u_i(\lambda))(\lambda) = - \eps_i.
\end{equation}
This equality means that
for any $x\in\u$, we have $\sum_i \lambda([x,u_i(\lambda)]) e^i = 
- x$, i.e., if $u_i(\lambda) = \sum_{j} e^j \otimes f_{i,j}(\lambda)$, that 
$\sum_{i,j} \lambda([x,e^j]) f_{i,j}(\lambda) \otimes e^i = 
- 1\otimes x$. Taking into account the identification of the function 
$\lambda \mapsto \lambda(x)$ with $x\in S^1(\l)$, (\ref{step:30}) 
now follows from (\ref{4:years}). 

We now show that if $f\in \big( \wedge^2(\g) \otimes 
\wh S^\cdot(\l)[1/\Delta] \big)^\l$, then 
\begin{equation} \label{identity:d}
(\on{d}f)_{|\k^*} - \on{d}(f_{|\k^*}) = 
- [r_{\k,\l}(\lambda)^{1,3} + r_{\k,\l}(\lambda)^{2,3}, 
(f_{|\k^*})^{1,2}]. 
\end{equation}
The l.h.s., evaluated at $\lambda\in\k^*$, is equal to 
$\sum_i {{\on{d}}\over{\on{d}t}}_{|t = 0} f(\lambda + t \eps^i)^{1,2}
(e^i)^3$, where $(\eps_i),(e^i)$ are dual bases of $\u^*$ and $\u$. 
According to (\ref{step:30}), this l.h.s. is equal to 
$- \sum_i {{\on{d}}\over{\on{d}t}}_{|t=0} f(\Ad(e^{t u_i(\lambda)})
(\lambda))^{1,2} (e^i)^3$, which by  invariance of $f$
is the r.h.s. of (\ref{identity:d}). 

Then we get 
\begin{align*}
& \on{CYB}(\sigma) - \on{Alt}(\on{d}\sigma) = 
\big( \on{CYB}(r_{\k,\l}) - \on{Alt}(\on{d}r_{\k,\l}) \big) 
\\ & 
+ \big( \on{CYB}(\rho) - \on{Alt}(\on{d}\rho_{|\k^*}) \big) 
+ \big( \on{CYB}(r_{\k,\l},\rho) - \on{Alt}( (\on{d}\rho)_{|\k^*} 
- \on{d}(\rho_{|\k^*}) ) \big).  
\end{align*}  
(Here $\on{CYB}(a,b)$ is the bilinear form derived from the quadratic form
$\on{CYB}$.) 
In this equality, the first term is zero by Proposition \ref{prop:2:1}, 
the second term is equal to  $Z$, and the last term is zero by 
(\ref{identity:d}). 
\hfill \qed \medskip 

\begin{remark} \label{rem:2:3}
In the case where $\g = \l\oplus\u$ is a nondegenerate Lie 
algebra with a splitting, $Z = 0$ and $\rho = r_{\l,\g}$, then 
$\sigma = r_{\k,\g}$. In the polarized case, a quantum analogue of this 
statement is Proposition \ref{prop:comp:formula}.  
\end{remark}

\subsection{Polarized Lie algebras}

We say that the Lie algebra with a splitting 
$(\g,\l,\u)$ is {\it polarized} if 
we are given a decomposition $\u = \u_+ \oplus \u_-$ of $\u$ as a sum of 
two $\l$-submodules, such that $\u_+$ and $\u_-$ are Lie subalgebras of 
$\g$. We denote by $\p_\pm$ the "parabolic"
Lie subalgebras $\p_\pm  = \l \oplus \u_\pm \subset \g$.     

Assume that $(\g,\l,\u)$ is nondegenerate and polarized; then 
$\on{dim}(\u_+) = \on{dim}(\u_-)$. In that case, $\lambda \mapsto 
r_\l^\g(\lambda)$ takes
its values in $\big((\u_+\otimes \u_-) \oplus(\u_-\otimes\u_+) \big) 
\cap \wedge^2(\u)$. Therefore $r_\l^\g = r' - (r')^{2,1}$, where 
$r'\in \u_+ \otimes \u_- \otimes S^\cdot(\l)[1/D_0]$. 
We will call $r'$ the "half $r$-matrix" of $(\g,\l,\u_+,\u_-)$. 

\subsection{Quantization} \label{sect:quant} 

Let $(\g,\l,\u)$ be a nondegenerate Lie algebra with a splitting.  
Let $D\subset U(\l)$ be a degree $\leq d$ element with 
symbol $D_0$. Define $\wh U$ as the microlocalization of $U(\l)$, 
obtained by inverting $D$ (\cite{Spr}). $U(\l)$ embeds
into $\wh U$, and $\wh U$ is independent on the choice of $D$ up 
to isomorphism. $\wh U$ is a complete filtered algebra, with 
associated graded $S^\cdot(\l)[1/D_0]$. 

Here is a description of $\wh U$. An element of $\wh U$
is represented by a series $\sum_{i\in \ZZ} a_i D^{-i}$,
where $a_i\in U(\l)$ vanish for $-i$ large enough, and the sequence 
$\on{deg}(a_i) - i d$ tends to $-\infty$ as $i\to \infty$. Two such series are
equivalent if they differ by a sum $\sum_i x_i$, where 
$x_i$ has the form $\sum_{k = p(i)}^{q(i)} \alpha_k(i)D^{-k}$,  
$\sum_{k=p(i)}^{q(i)} \alpha_k(i) D^{q(i)-k} = 0$ (equality in $U(\l)$)
and $\on{max}_{k=p(i)}^{q(i)}(\on{deg}\alpha_k(i) - k d) \to -\infty$
as $i\to\infty$.  
The degree of $f$ is the minimum of all $\on{max}_i(\on{deg}(a_i) - id)$
running over all $\sum_i a_i D^{-i}$ representing $f$. 
The product of two elements  by 
$$
\big( \sum_{i\in\ZZ} a_i D^{-i} \big) 
\big( \sum_{j\in\ZZ} b_j D^{-j} \big) 
= \sum_{k\in\ZZ} \big( \sum_{\alpha\geq 0}
\pmatrix -i \\ \alpha \endpmatrix a_i \ad^{\alpha}(D)(b_j)
\big) D^{-k}. 
$$

We denote by $\wh U_{\leq k}$ the degree $\leq k$ part of 
$\wh U$. Then if $V$ is a vector space, we set $V\wh\otimes \wh U 
= \on{lim}_{\leftarrow} (V\otimes \wh U) / (V\otimes \wh U_{\leq k})$. 
The coproduct map $\Delta$ of $U(\l)$ extends to a map 
$\wh U \to U(\l) \wh\otimes \wh U$, where the image of 
$D^{-1}$ is $\sum_{i\geq 0} (-1)^i 
(1\otimes D^{-1}) a \cdots a (1\otimes D^{-1})$, and 
$a = \Delta(D) - 1\otimes D \in U(\l) \otimes U(\l)_{\leq d-1}$.  

\begin{definition}
{\it A quantization of $r_\l^\g$ is a $\l$-invariant 
element $J\in U(\g)^{\otimes 2}
\wh\otimes \wh U_{\leq 0}$, satisfying the dynamical twist equation 
$$
J^{12,3,4} J^{1,2,34} = J^{1,23,4} J^{2,3,4},
$$
and the following conditions:  
$J-1 \in U(\g)^{\otimes 2} \wh\otimes \wh U_{\leq -1}$, 
the reduction $j$ of $J-1$ modulo 
$U(\g)^{\otimes 2} \wh\otimes \wh U_{\leq -2}$
(an element of $U(\g)^{\otimes 2} \otimes S^\cdot(\l)[1/D_0]_{-1}$) 
satisfies $\on{Alt}(j) = r_\l^\g$. 
} \end{definition}

Then $R:= (J^{2,1,3})^{-1} J^{1,2,3}$ is a solution of the 
dynamical quantum Yang-Baxter equation 
$R^{1,2,4} R^{1,3,24} R^{2,3,4} = R^{2,3,14} R^{1,3,4} R^{1,2,34}$. 
 
The PBW star-product on $\l^*$ may be described as follows: 
$\hbar \l[[\hbar]] \subset \l[[\hbar]]$ is a Lie subalgebra, 
then the $\hbar$-adic completion of $U(\hbar \l[[\hbar]]) 
\subset U(\l)[[\hbar]]$ is a flat deformation of $\wh S^\cdot(\l)$, 
which we denote by $\wh S^\cdot(\l)_\hbar$ 
(it is the quantized formal series algebra 
associated to the trivial deformation of $U(\l)$). Then 
$U(\l)((\hbar))$ identifies with $\wh S^\cdot(\l)_\hbar[\hbar^{-1}]$; 
moreover, this identification takes 
$U(\l)_{\leq k}[[\hbar]]$ into $\hbar^{-k} \wh S^\cdot(\l)_\hbar$. 

This discussion can be localized. We denote by 
$\wh S^\cdot(\l)[1/D_0]_\hbar$ the $\hbar$-adic completion of the 
subalgebra of $\wh U((\hbar))$ generated by $\hbar \l[[\hbar]]$
and $(\hbar^d D)^{-1}$. This is a flat deformation of 
$\wh S^\cdot(\l)[1/D_0]$. Moreover, $\wh U((\hbar))$ 
identifies with $\wh S^\cdot(\l)[1/D_0]_\hbar[\hbar^{-1}]$, 
and $\wh U_{\leq k}[[\hbar]]$ goes into 
$\hbar^{-k}\wh S^\cdot(\l)[1/D_0]_\hbar$. 

It follows that $J$ gives rise to an element 
$J(\la)$ of $U(\g)^{\otimes 2} \wh\otimes 
\wh S^\cdot(\l)[1/D_0]_\hbar$, with the 
expansion $J(\lambda) = 1 + \hbar j(\lambda) + O(\hbar^2)$, 
and $\on{Alt}(j(\lambda)) = r_\l^\g$. 
    
In Section  \ref{quant:pol}, we will quantize the classical dynamical
$r$-matrices arising from nondegenerate polarized Lie algebras. 

\begin{remark} \label{rem:char}
The continuous characters $\chi : \wh U \to \CC((\hbar))$ are all
of the following form: $\lambda : \l \to \CC((\hbar))$ is a character of 
$\l$, of the form $\lambda = \sum_{i\geq v} \hbar^i \lambda_i$, 
with $v<0$ and $D_0(\lambda_v)\neq 0$. Then $\lambda$ is a character $U(\l)
\to \CC((\hbar))$, it extends to a character $\wh U\to\CC((\hbar))$, 
which restricts to $\wh U_{\leq 0}\to \CC[[\hbar]]$. 
\end{remark}

 \begin{remark} {\it Microlocalization.} Springer's microlocalization
 associates to a pair $(A,f)$, where $A$ is a $\ZZ$-filtered algebra
 with $\on{gr}(A)$ integral commutative and $f\in A$ is nonzero, a complete 
 separated $\ZZ$-filtered algebra $A_f$, such that $\on{gr}(A_f) = 
 \on{gr}(A)[1/\bar f]$
 (here $\bar f$ is the symbol of $f$, i.e., its nonzero homogeneous component
 with maximal degree). $(A_f)_{\leq 0}$ is a subalgebra of $A_f$ and contains
 $(A_f)_{\leq -1}$ as an ideal. 
 
 $A_f$ has the following universal property: if $B$ is a $\ZZ$-filtered, 
 complete separated algebra (i.e., $\cap_i B_i = \{0\}$ and $B =
 \lim_{\leftarrow i}(B/B_i)$), and $\mu : A \to B$ is a morphism of 
 filtered algebras, such that $\mu(f)$ is invertible, then $\mu$
 extends to a morphism of topological filtered algebras $A_f \to B$. 

 Actually, $A_f$ depends only on $\bar f$, and when $A$ is graded, $A_f$
 is the completion of its associated graded. E.g., if 
 $A = \CC[x_1,\ldots,x_n]$ and $f\in A - \{0\}$ is homogeneous, these
 algebras can be described as follows. Let $C(f) = C\subset
 \CC^n$ be the cone defined by the equation $f = 0$. Then 
 $\on{gr}(A_f)$ is the ring on functions on $\CC^n - C$.  
 The projective space $\PP^n$ decomposes
 as $\CC^n \cup H$, where $H$ is the hyperplane at infinity, and the closure
 $\overline C$ of
 $C$ in $\PP^n$ decomposes as $C \cup C_\infty$, where $C_\infty = C\cap H$.
 Then $\on{gr}((A_f)_{\leq 0})$ is the ring of functions on $\PP^n - 
 \overline C$,
 the quotient $(A_f)_{\leq 0} / (A_f)_{\leq -1}$ is the ring of
 $\CC^\times$-invariant functions on $\CC^n - C$. Finally, $(A_f)_{\leq 0}$
 (resp., $A_f$) is the ring of functions on the formal (resp., formal 
 punctured) neighborhood on $H - C_\infty$ in $\PP^n - \overline C$. 
 
 In general, if $A$ is a $\ZZ_+$-filtered commutative algebra and 
 $X = \on{Spec}(A)$, then $X$ has a 
 compactification $\overline X = X \cup X_\infty$. Here $\overline X 
 = \on{Proj}(R(A))$, where $R(A)$ is the Rees algebra of $A$, and $X_\infty =
 \on{Proj}(\on{gr}(A))$. If $g\in A - \{0\}$, then 
  $A_g$ (resp., $(A_g)_{\leq 0}$) is
 the ring of functions on the formal (resp., formal punctured) neighborhood of 
 $X_\infty - C_\infty(g)$, where $C_\infty(g) = \overline{V(g)} \cap X_\infty$, 
 and $V(g)\subset X$ is the zero-set of $g$. $C_\infty(g)$ depends only on
 $\bar g$, which explains why the same is true about $A_g$. 
 \end{remark}

\subsection{Examples}

\subsubsection{Lie algebras with a splitting}

(1) An inclusion $\l\subset \g$ of simple Lie algebras with the same Cartan 
algebra $\h$ is called a Borel-de Siebenthal pair (\cite{BS}). Then $\l$
has an invariant complement $\u$. If $\lambda\in\h^*$, the bilinear form 
$\g^2\to\CC$, $(x,y)\mapsto \langle\lambda,[x,y]\rangle$
is nondegenerate for $\lambda$ generic, and is the sum of two bilinear forms
$\l^{2} \to \CC$ and $\u^2 \to\CC$, which are therefore nondegenerate. 
In particular, $\u^2\to\CC$, $(x,y)\mapsto \langle\lambda,[x,y]\rangle$
is nondegenerate. So 
$(\g,\l,\u)$ is a nondegenerate Lie algebra with a splitting. 

(2) If $\g$ is a finite dimensional Lie algebra and $r\in\wedge^2(\g)$
is a nondegenerate triangular $r$-matrix, then the dual of $r$ is a 
$2$-cocycle on $\g$. Let $\wh\g = \g \oplus\CC c$ be the corresponding 
central extension. Then $\wh\g$ is a nondegenerate Lie algebra with 
a splitting
with $\l = \CC c$, $\u = \g$. The corresponding $r$-matrix is 
$\lambda \mapsto r/\lambda$. 

(3) We generalize (2) to the case when $\l$ is no longer $1$-dimensional. 
Let $\wh\g$ be a Lie algebra, let $\z \subset \g$ be a central subalgebra, 
set $\g:= \wh\g/\z$, and let $\pi : \wh\g\to\g$ be the canonical projection. 
Let $\u\subset \wh\g$ be a complement of $\z$. Then $(\wh\g,\z,\u)$
is a Lie algebra with a splitting; let us assume it is nondegenerate. 
Set $r:= (\pi\otimes \pi\otimes \id)(r_\l^\g)$. Then $r$ satisfies 
$\on{CYB}(r) = 0$. In particular, for any $\lambda\in \l^*$ such that 
$D_0(\lambda)\neq 0$, $r_\lambda := (\id\otimes\id\otimes\lambda)(r)$
is a triangular $r$-matrix (we identify $\lambda$ with a character of 
$S^\cdot(\l)[1/D_0]$). 
If $J$ is a quantization of $r_\l^\g$, and $\chi : 
\wh U_{\leq 0} \to \CC[[\hbar]]$ is a character as in Remark
\ref{rem:char}, then $F_\chi:= (\pi\otimes\pi\otimes\chi)(J)$ 
is a solution of the twist equation, quantizing $r_\lambda$.

\subsubsection{Polarized Lie algebras} \label{2.5.2}

(1) If $\g$ is a semisimple Lie algebra and $\l\subset \g$ 
is a Levi subalgebra, then $(\g,\l)$ gives rise to a nondegenerate 
polarized Lie algebra, which was studied in \cite{DM}. Then 
$r_\l^\g : \l^* \to \wedge^2(\u)$ is defined by 
$$
r_\l^\g(\lambda) = - \sum_{\alpha\in \Delta_+(\g) - \Delta_+(\l)}
(\on{ad}(\lambda^\vee))^{-1}(e_\alpha) \wedge f_\alpha,
$$ for 
$\lambda\in\l^*$ such that $\ad(\lambda^\vee)_{|\u} \in \End(\u)$
is invertible. $r_\l^\g$ is also uniquely determined by the requirements 
that it is an $\l$-equivariant rational function, such that 
$$
\forall \lambda\in\h^*, \quad 
r_\l^\g(\lambda) = - \sum_{\alpha \in \Delta_+(\g) - \Delta_+(\l)}
{{e_\alpha \wedge f_\alpha}\over{(\lambda,\alpha)}}. 
$$
Here $\Delta_+(\g),\Delta_+(\l)$ are the sets of positive roots of 
$\g,\l$, and $x\wedge y = x\otimes y - y \otimes x$. 

(2) Let $\g$ be a finite dimensional Lie algebra, which can be
decomposed (as a vector space) as $\g = \g_+ \oplus \g_-$, where 
$\g_\pm \subset \g$ are Lie subalgebras. Let $\wt r\in\g_+\otimes \g_-$
be a nondegenerate tensor, such that $r:= \wt r - (\wt r)^{2,1}$ is a 
triangular $r$-matrix (i.e., it satisfies the CYBE). Then 
we may construct $\wh\g$ as above. If we set $\l = \CC c$, 
$\u_\pm = \g_\pm$, we get a nondegenerate 
polarized Lie algebra. 

(3) Let $\g = \l \oplus\u_+\oplus\u_-$ be a nondegenerate polarized
Lie algebra, and $A = \CC[t] / (t^n)$, then $\g\otimes A = (\l \otimes A)
\oplus (\u_+\otimes A) \oplus (\u_- \otimes A)$ is a nondegenerate polarized 
Lie algebra.   

\subsubsection{Infinite dimensional examples}

The definitions of 
polarized Lie algebras, their $r$-matrices and quantizations 
generalize to the case of graded Lie algebras with 
finite dimensional graded parts. 

(1) The Virasoro algebra $\on{Vir}$ decomposes as $\l \oplus \u_+
\oplus \u_-$, where $\l = \CC c \oplus \CC L_0$ and 
$\u_\pm = \oplus_{i>0} \CC L_{\pm i}$. 
Then $(\on{Vir},\l,\u_+,\u_-)$
is an infinite dimensional polarized Lie algebra.  

(2) If $\g$ is a Kac-Moody Lie algebra and $\l\subset \g$
is a Levi subalgebra, then $(\g,\l)$ gives rise to an infinite 
dimensional polarized Lie algebra. 

\begin{remark} The proof of Proposition \ref{prop:comp} shows that 
if $\g = \l\oplus \u$ is a Lie algebra with a nondegenerate splitting, 
$U\subset \g^*$ is an invariant open subset and $r : U \to \wedge^2(\g)$
is $\g$-invariant, such that $r_{|\g^*} + r_\l^\g$ is a 
$(\l,\g,Z)$-$r$-matrix, then $r$ is a $(\g,\g,Z)$-$r$-matrix. 
This leads to the following $r$-matrix (a quantization of 
which is unknown). 

Let $\g$ be a semisimple Lie algebra and $t\in S^2(\g)^\g$ be nondegenerate. 
For $\xi\in \g^*$, set $\xi^\vee = (\xi\otimes \id)(t)$. If $\h'\subset \g$
is a Cartan subalgebra, let $t_{\h'}$ be the part of $t$ corresponding to 
$\h'$. Set $\g^*_{\on{ss}} = \{\xi\in\g^* | \xi^\vee$ is semisimple$\}$. 
If $x\in\g$ is semisimple, let $\h_x = \{h\in\g | [h,x]=0\}$ be the Cartan
subalgebra associated to $x$. 
Then the map $\g^*_{\on{ss}} \to \wedge^2(\g)$, $\xi \mapsto 
(\on{ad}(\xi)^{-1} \otimes \id)(t - t_{\h(\xi)})$ is a 
$(\g,\g,0)$-$r$-matrix. 
\end{remark}

\section{Dynamical twists in the polarized case} \label{quant:pol}
\label{sect:K} 
  
In this section, we constuct a dynamical twist $J_\l^\g$ quantizing 
$r_\l^\g$. For this, we first construct an element $K$; it is defined by 
algebraic requirements, related with the Shapovalov form. We then 
construct $J = J_\l^\g$ and show that it obeys the dynamical twist equation. We then 
show that $J$ satisfies nonabelian versions of the ABRR equations.

\subsection{Construction of $K$}
Let $\g = \l\oplus \u_+ \oplus \u_-$ be a nondegenerate polarized 
Lie algebra. Denote by $H : U(\g) \to U(\l)$ the Harish-Chandra map, 
defined as the unique linear map such that 
$H(x_+ x_0 x_-) = \eps(x_+) \eps(x_-) x_0$ if 
$x_0 \in U(\l)$ and $x_\pm\in U(\u_\pm)$. Here $\eps : 
U(\g)\to\CC$ is the counit map. 

Let $d' = \on{dim}(\u_\pm)$ and let $D'_0\in S^{d'}(\l)$
be the polynomial taking $\lambda\in\l^*$ to $\on{det}
(\lambda\circ\omega) \circ i$, where $\omega : \u_+ \otimes \u_- \to \l$
is the Lie bracket followed by projection, $\lambda\circ \omega$ is viewed
as a linear map $\u_+^*\to \u_-$ and $i$ is a fixed linear isomorphism 
$\u_- \to \u_+^*$. 

The relation with the objects introduced 
in the previous section is 
$d' = d/2$ and $(D'_0)^2 = D_0$. In particular, 
$\wh U$ identifies with the microlocalization of $U(\l)$
with respect to a lift $D'\in U(\l)_{\leq d'}$ of $D'_0$.

\begin{theorem} \label{thm:J0}
There exists a unique element $K\in \big(U(\u_+) \otimes U(\u_-)\big)
\wh\otimes \wh U$, such that if we set $K = \sum_i e_i^+ \otimes e_i^-
\otimes \ell_i$, then we have 
\begin{equation}
\forall x\in U(\p_-), \forall y\in U(\p_+), \;
\sum_i H(xe_i^+) \ell_i H(e_i^- y) = H(xy). 
\end{equation}
Equivalently, we have for any  $x_\pm\in U(\u_\pm)$, 
$x_0\in U(\l)$, 
$$
\sum_i H(xe_i^+) \ell_i \otimes e_i^- = \eps(x_+)x_0 \otimes x_-, 
\; 
\sum_i e_i^+ \otimes  \ell_i H(e_i^- y) = x_+ \otimes x_0 \eps(x_-),   
$$
where $x = x_+x_0 x_-$. 
$K$ has also the following properties. $K$ is invariant under the 
adjoint action of $\l$. $K$ is a sum $\sum_{n\geq 0} K_n$, where 
$K_n \in \big(U(\u_+)_{\leq n} \otimes U(\u_-)_{\leq n}\big) 
\wh\otimes \wh U_{\leq -n}$, 
and the image of $K_n$ in $S^n(\u_+) \otimes S^n(\u_-) \otimes
S^\cdot(\l)[1/D'_0]_{-n}$ under the tensor product of the projection maps 
$U(\u_\pm)_{\leq n} \to U(\u_\pm)_{\leq n} / U(\u_\pm)_{\leq n-1}$
and $\wh U_{\leq -n} \to \wh U_{\leq -n} / \wh U_{\leq -n-1}$, coincides with 
${1\over {n!}} (r'')^n$, 
where $r'' := -r' \in \u_+\otimes \u_- \otimes S^\cdot(\l)[1/D'_0]_{-1}$
is the opposite of the 
"half $r$-matrix" of $(\g,\l,\u_+,\u_-)$ (the index $k$ means the 
homogeneous part of degree $k$; the index $\leq k$ means the part
of degree $\leq k$; the algebra structure of 
$S^\cdot(\u_+) \otimes S^\cdot(\u_-) \otimes S^\cdot(\l)[1/D'_0]$
is understood). 
\end{theorem}

{\em Proof.} $r''$ is a sum $\sum_i a_i \otimes b_i \otimes 
P_i(D'_0)^{-1}$, with $P_i \in S^{d'-1}(\l)$. Let $\bar P_i\in 
U(\l)_{\leq d'-1}$ be a lift of $P_i$, and let $\bar r:= \sum_i a_i \otimes
b_i \otimes \bar P_i (D')^{-1}$. Then $\bar r\in \u_+ \otimes \u_- 
\otimes \wh U_{\leq -1}$ is a lift of $r''$. Let us set $\bar K_0 :=
\exp(\bar\rho)$. 

\begin{lemma} 
Set $\bar K = \sum_i \bar e^+_i \otimes \bar e^-_i \otimes \bar\ell_i$.  
If $x\in U(\u_-)$, set $T(x) = \sum_i H(x\bar e_i^+) \bar\ell_i 
\otimes \bar e_i^-$. Then $T$ is a linear map 
$U(\u_-) \to \wh U_{\leq 0} \wh\otimes U(\u_-)$, such that if $x$ has 
degree $\leq n$, then 
\begin{equation} \label{**}
T(x) - 1\otimes x\in \wh U_{\leq 0}
\wh\otimes U(\u_-)_{\leq n-1} + \wh U_{\leq -1} \wh\otimes U(\u_-) . 
\end{equation}   
\end{lemma}

{\em Proof.} The map $H$ is such that if $x_\pm\in U(\u_\pm)$
have degrees $\leq n_\pm$, then $H(x_- x_+)\in U(\l)$ has degree
$\leq \on{min}(n_+,n_-)$. The  
bilinear map $U(\u_-) \otimes U(\u_+)
\to U(\l)$, $x_- \otimes x_+ \mapsto H(x_-x_+)$ therefore 
induces a collection of bilinear maps 
$S^n(\u_-) \otimes S^n(\u_+) \to S^n(\l)$, which turn out to be the 
symmetric powers of $h : \u_+ \otimes \u_- \to \l$. 

Write $\bar r^n = \sum_i a_i^{(n)} \otimes b_i^{(n)} \otimes \ell_i^{(n)}$, 
where $a_i^{(n)}$, $b_i^{(n)}$ have degree $\leq n$ and $h_i^{(n)}$ has degree
$\leq -n$. Then if $x$ has degree $\leq k$, $H(xa_i^{(n)}) \ell_i^{(n)}$
has degree $\leq \on{min}(k,n) - n \leq 0$. Moreover, this degree tends to 
$-\infty$ as $n\to\infty$, so $T$ is well-defined and maps to 
$\wh U_{\leq 0} \wh\otimes U(\u_-)$. 

Let us prove (\ref{**}) when $x\in\u_-$. We have $\sum_i h(x\otimes a_i) P_i 
\otimes b_i = D'_0\otimes x$ (equality in $S^{d'}(\l) \otimes \u_-$), 
so $\sum_i H(xa_i) \bar P_i \otimes b_i\in D' \otimes x 
+ U(\l)_{\leq d'-1} \otimes \u_-$. Then $\sum_i
H(xa_i) \bar P_i (D')^{-1} \otimes b_i \in 1 \otimes x + 
\wh U_{\leq -1} \otimes
\u_-$.   On the other hand, if $n>1$, then $H(xa_i^{(n)}) \ell_i^{(n)}
\otimes b_i^{(n)}$ has degree $\leq 1-n \leq -1$. So $T(x) - 1\otimes x\in 
\wh U_{\leq -1} \wh\otimes U(\u_-)$. 

Let now $x\in U(\u_-)$ be of degree $k$. If $n<k$, then 
$H(xa_i^{(n)}) \ell_i^{(n)}$ has degree $\leq 0$, so 
$\sum_i H(xa_i^{(n)}) \ell_i^{(n)} \otimes b_i^{(n)}
\in \wh U_{\leq 0} \otimes U(\u_-)_{\leq n}$, so the 
contribution of $\sum_{n<k} \bar r^n/n!$ lies in
$\wh U_{\leq 0} \otimes U(\u_-)_{\leq k-1}$. If $n = k$, then 
$\sum_i H(xa_i^{(k)})
\ell_i^{(k)} \otimes b_i^{(k)}\in \wh U_{\leq 0} \otimes U(\u_-)_{\leq k}$, 
and its class modulo $\wh U_{\leq -1} \otimes U(\u_-)_{\leq k} +  
\wh U_{\leq 0} \otimes U(\u_-)_{\leq k-1}$ is $1\otimes x$, by the filtration
properties of $x\otimes y \mapsto H(xy)$. If $n>k$, then $H(xa_i^{(n)})
\ell_i^{(n)}$ has degree $\leq k-n \leq -1$. This shows that $T(x) - 
1\otimes x$ has the required degree properties. 
\hfill \qed \medskip 

\begin{lemma} $T$ extends uniquely to a continuous endomorphism 
$\wt T$ of $\wh U_{\leq 0} \wh\otimes U(\u_-)$, such that 
$\wt T(\ell \otimes x) = (\ell \otimes 1) T(x)$ for any 
$\ell\in \wh U_{\leq 0}$ and $x\in U(\u_-)$. $\wt T$ is invertible, and 
$T' := (\wt T^{-1})_{| 1 \otimes U(\u_-)}$ has the same degree properties 
as $T$: if $x\in U(\u_-)$ has degree $\leq n$, then 
$$
T'(x) - 1 \otimes x \in \wh U_{\leq 0} \wh\otimes U(\u_-)_{\leq n-1}
+ \wh U_{\leq -1} \wh\otimes U(\u_-). 
$$   
\end{lemma}

{\em Proof.} Clear. \hfill \qed \medskip 

{\em End of proof of Theorem \ref{thm:J0}.}   
We now set $K := \sum_i (\bar e_i^+ \otimes 1 \otimes \bar\ell_i)
(1\otimes T'(\bar e_i^-)^{2,1})$. Set $K = \sum_i e_i^+ \otimes e_i^- 
\otimes \ell_i$. Then if $x\in U(\u_-)$, we have 
\begin{align*}
& \sum_i H(xe_i^+) \ell_i \otimes e_i^- 
= \sum_i (H(x\bar e_i^+)\bar \ell_i \otimes 1) T'(\bar e_i^-)
\\ & 
= \sum (T(x)^{(1)} \otimes 1) T'(T(x)^{(2)}) = \wt T^{-1}(T(x)) = x, 
\end{align*}
where we have set $T(x) = \sum T(x)^{(1)} \otimes T(x)^{(2)}$. 

The properties of $T'$ then imply the following. If $n\geq 0$, 
then the reduction of 
$K$ modulo $\big( U(\u_+) \otimes U(\u_-)\big) 
\wh\otimes \wh U_{\leq -n-1}$ lies in $\big( U(\u_+)_{\leq n} \otimes 
U(\u_-)_{\leq n} \big) 
\otimes (\wh U_{\leq 0} / \wh U_{\leq -n-1})$, and its reduction modulo 
$(U(\u_+) \otimes U(\u_-))_{<2n} \otimes (\wh U_{\leq 0} / 
\wh U_{\leq -n-1})$ lies in $S^n(\u_+) \otimes S^n(\u_-) \otimes 
S^\cdot(\l)[1/D'_0]_{-n}$; it identifies with ${1\over {n!}} (r'')^n$. 
This implies the claim on the decomposition on $K$. 

We now prove the uniqueness of $K$. If $K'$ has the same properties as 
$K$, then $K'':= K' - K = \sum_i a'_i \otimes b'_i \otimes \ell'_i$
is such that for any $x\in U(\u_-)$, $\sum_i H(xa'_i) \ell'_i \otimes 
b'_i = 0$. Let $(e_I^-)$ be a basis of $U(\u_-)$, and set 
$K'' = \sum_{I,i} a_{i,I} \otimes e^-_I \otimes \ell_{i,I}$, then 
$\sum_i H(xa_i) \ell_{i,I} = 0$ for any $I$. We now prove: 

\begin{lemma}
If $\xi = \sum_i a_i \otimes \ell_i \in U(\u_+)\wh\otimes\wh U$
is such that $\sum_i H(xa_i) \ell_i = 0$ for any $x\in U(\u_-)$, then 
$\xi = 0$.
\end{lemma}   

{\em Proof of Lemma.} Set $\xi = \sum_\al \xi_\alpha$, where
$\on{deg}(\xi_\alpha) = -\alpha$, Let $\alpha_0$ be the largest integer such
that $\xi_{\alpha_0} \neq 0$. We have $\xi_{\alpha_0} = \sum_{s=0}^N \eta_s$,
where $\eta_s\in U(\u_-)_{\leq s} \wh\otimes \wh U_{\alpha_0 - s}$. Then if
$x\in U(\u_-)$ has degree $\leq n$, then $m \circ (H\otimes \on{id})
\big( (x\otimes 1)\eta_s\big)\in \wh U$ has degree $\leq \on{min}(n,s)
+ \alpha_0 - s$. Pairing $\eta_s$ with $U(\u_-)_{\leq N}$, 
$U(\u_-)_{\leq N-1}$, etc.,  we get $\eta_N\in U(\u_+)_{\leq N-1} 
\wh\otimes \wh U_{\alpha_0 - N}$, etc. Finally $\xi_{\alpha_0} = 0$, 
and $\xi = 0$. 
\hfill \qed \medskip 

Therefore $K'' = 0$, so $K$ is unique. Then its $\l$-invariance follows from 
the $\l$-invariance of $H$. 
\hfill \qed \medskip

\subsection{The dynamical twist equation}

If $K = \sum_i a_i \otimes b_i \otimes \ell_i$, set 
$J = J_\l^\g := \sum_i a_i \otimes S(b_i)S(\ell_i^{(2)}) 
\otimes S(\ell_i^{(1)})$. 
Then $J\in \big( U(\u_+) \otimes U(\p_-) \big) 
\wh\otimes \wh U_{\leq 0}$. 

\begin{proposition} \label{prop:K:twist:eqn}
$J$ satisfies the dynamical twist equation 
\begin{equation} \label{twist:eq}
J^{12,3,4} J^{1,2,34} = J^{1,23,4} J^{2,3,4}. 
\end{equation}
\end{proposition}

This proposition has a representation-theoretic interpretation 
in terms of intertwiners, analogous to that of the abelian case 
(see \cite{EV2} or \cite{ES1}, Proposition 2.3). 

\medskip 
{\em Proof.} Let us set $K = \sum_i a_i \otimes b_i \otimes 
\ell_i$. Then (\ref{twist:eq}) can be written as follows
\begin{align*}
& \sum_{i,j} a_i^{(1)} a_j \otimes a_i^{(2)} S(b_j) S(\ell_j^{(3)}) 
\otimes S(b_i) S(\ell_i^{(2)}) S(\ell_j^{(2)}) \otimes 
S(\ell_i^{(1)}) S(\ell_j^{(1)})
\\ & = \sum_{i,j}
 a_i \otimes S(b_i^{(2)}) S(\ell_i^{(3)}) a_j 
 \otimes S(b_i^{(1)}) S(\ell_i^{(2)}) S(b_j) S(\ell_j^{(2)})
 \otimes S(\ell_i^{(1)}) S(\ell_j^{(1)}) . 
\end{align*} 
Since $K$ is $\l$-invariant, the right-hand side is rewritten as 
$$
\sum_{i,j} a_i \otimes S(b_i^{(2)}) a_j S(\ell_i^{(3)})
\otimes S(b_i^{(1)}) S(b_j) S(\ell_j^{(2)}) S(\ell_i^{(2)})
\otimes S(\ell_j^{(1)}) S(\ell_i^{(1)}). 
$$
Now both sides belong to the image of the map 
$$
\big( U(\u_+) \otimes U(\g)\otimes U(\u_-) \big) \wh\otimes \wh U 
\to 
\big( U(\u_+) \otimes U(\g)\otimes U(\p_-) \big) \wh\otimes \wh U,  
$$
$$
x\otimes y \otimes z \otimes t\mapsto x \otimes y \otimes S(z) 
S(t^{(2)}) \otimes S(t^{(1)}).  
$$
So we have to prove the equality  
\begin{equation} \label{twist'}
\sum_{i,j} a_i^{(1)} a_j \otimes a_i^{(2)} S(b_j) S(\ell_j^{(2)})
\otimes b_i \otimes \ell_j^{(1)} \ell_i
= \sum_{i,j} a_i \otimes S(b_i^{(2)}) a_j S(\ell_i^{(2)})
\otimes b_j b_i^{(1)} \otimes \ell_i^{(1)} \ell_j 
\end{equation}
in $\big( U(\u_+) \otimes U(\g)\otimes U(\u_-)\big) \wh\otimes\wh U$.

The linear map 
$$
\big( U(\u_+) \otimes U(\g)\otimes U(\u_-)\big) \wh\otimes\wh U
\to 
\on{Hom}_\CC(U(\u_-) \otimes U(\u_+), U(\g) \wh\otimes \wh U), 
$$
$$
A \otimes B \otimes C \otimes D \mapsto \Big( x\otimes y \mapsto 
B S\big(H(xA)^{(2)} \big) \otimes H(xA)^{(1)} D H(By)\Big)  
$$
is injective. This map takes the l.h.s. of (\ref{twist'})
to 
$$
\alpha : x\otimes y \mapsto 
\sum_i y^{(2)} S\big( (xy^{(1)})_{-,i}\big) 
S\big( \big( (xy^{(1)})_{0,i} \big)^{(2)} \big)
\otimes \big( (xy^{(1)})_{0,i} \big)^{(1)} \eps((xy^{(1)})_{+,i}), 
$$ 
and the r.h.s. of (\ref{twist'}) to 
$$
\beta : x\otimes y \mapsto \sum_i S(x^{(2)}) (x^{(1)}y)_{+,i}
\otimes (x^{(1)}y)_{0,i} \eps((x^{(1)}y)_{-,i}). 
$$
Here we denote by $\sum_i x_{+,i} \otimes x_{0,i} \otimes x_{-,i}$
the image of $x\in U(\g)$ in $U(\u_+)\otimes U(\l) \otimes U(\u_-)$
by the inverse of the product map. 

To prove that $\alpha = \beta$, we will prove that the maps 
$(x,y) \mapsto (S(y^{(2)}) \otimes 1)\alpha(x\otimes y^{(1)})$
and $(x,y) \mapsto (S(y^{(2)}) \otimes 1) \beta(x\otimes y^{(1)})$
coincide. 
The first map takes $(x,y)$ to 
$$
\sum_i S((xy)_{-,i}) S\big( \big( (xy)_{0,i} \big)^{(2)} \big) 
\otimes \big( (xy)_{0,i} \big)^{(1)} 
\eps((xy)_{+,i}) 
$$ 
and the second map takes $(x,y)$ to 
$$
\sum_i S((xy)^{(2)}) ((xy)^{(1)})_{+,i} \otimes ((xy)^{(1)})_{0,i}
\eps(((xy)^{(1)})_{-,i}). 
$$   
To prove the equality of both maps, it suffices to prove that the maps 
$U(\g) \to U(\g) \otimes U(\l)$, 
$$
a\mapsto \sum_i S(a_{-,i}) S((a_{0,i})^{(2)}) \otimes (a_{0,i})^{(1)}
\eps(a_{+,i})
$$
and
$$
a\mapsto \sum_i S(a^{(2)}) (a^{(1)})_{+,i} \otimes (a^{(1)})_{0,i}
\eps((a^{(1)})_{-,i})
$$
coincide. If $a_0\in U(\l)$ and $a_\pm\in U(\u_\pm)$, then the 
first map takes $a_+ a_0 a_-$ to 
$S(a_-)S((a_0)^{(2)}) \otimes (a_0)^{(1)} \eps(a_+)$, 
and the second map takes $a_+ a_0 a_-$
to $S((a_-)^{(2)}) S((a_0)^{(2)}) S((a_+)^{(2)}) (a_+)^{(1)}
\otimes (a_0)^{(1)} \eps((a_-)^{(1)})$, so both maps coincide. 
\hfill \qed \medskip 

Together with the valuation results of Theorem \ref{thm:J0}, 
and taking into account the change of sign induced by $S$, 
Proposition \ref{prop:K:twist:eqn} implies: 

\begin{corollary} \label{cor:quant:rat}
$J$ is a quantization of $r_\l^\g$, in the sense of 
Section \ref{sect:quant}. 
\end{corollary}

\begin{example} If $\g$ is the Heisenberg algebra, spanned by $x_+,x_-,c$, 
with $[x_+,x_-] = c$, $\u_\pm = \CC x_\pm$, $\l = \CC c$, then 
$K = \on{exp}(-x_+ \otimes x_- \otimes c^{-1})$, so that 
$J = \on{exp}(-x_+^{(1)} x_-^{(2)} (c^{(2)} + c^{(3)})^{-1})$, 
i.e., $J(\lambda) = \on{exp}(-x_+ \otimes x_-(\lambda + c)^{-1})$. 
\end{example}

\subsection{$K$, singular vectors, and fusion of intertwiners}

If $u\in \u_-$, then $x\mapsto [u,x]$ maps $U(\u_+)$ to $U(\u_+)\l\subset
U(\g)$. 

\begin{proposition}
If $u\in\u_-$, then $[u^{(1)},K] - K u^{(2)} \in \on{Im}(\varphi)$, 
where $\varphi : \big( U(\u_+) \otimes U(\u_-) \otimes \l\big) \wh\otimes
\wh U_{\leq -1} \to \big( U(\p_+) \otimes U(\u_-)\big) \wh\otimes \wh U_{\leq 0}$
is the map 
taking $x_+ \otimes x_- \otimes \ell \otimes \wh x$ to 
$x_+\ell \otimes x_- \otimes \wh x - x_+ \otimes x_- \otimes \ell \wh x$.  
\end{proposition}

{\em Proof.} Set $K = \sum_i a_i \otimes b_i \otimes \ell_i$. 
Then if $x\in U(\u_-)$, $y\in U(\u_+)$, we have 
$$
\sum_i H(x [u,a_i]) \ell_i H(b_i y) = \sum_i H((xu)a_i) \ell_i H(b_iy), 
$$
because $H(\xi u) = 0$ for any $\xi\in U(\g)$. Now 
$$
\sum_i H((xu)a_i) \ell_i H(b_iy) = H(xuy) = \sum_i H(xa_i) \ell_i
H(b_i u y). 
$$

So if $L = [u^{(1)},K] - K u^{(2)}$ is decomposed as $\sum_i 
\alpha_i \otimes \beta_i \otimes \lambda_i$,  we get 
$\sum_i H(x\alpha_i) \lambda_i H(\beta_i y) = 0$, therefore 
for any $x\in U(\u_-)$, $\sum_i H(x\alpha_i)\lambda_i \otimes \beta_i = 0$. 
Now if $\sum_i \alpha'_i \otimes \lambda'_i \in \big( U(\u_+) \oplus 
U(\u_+) \l\big) \wh\otimes\wh U_{\leq 0}$ is such that $\sum_i H(x\alpha'_i) \lambda'_i =
0$, there exists $\sum_i \alpha''_i \otimes \ell_i \otimes \lambda''_i \in
\big( U(\u_+) \otimes \l \big) \wh\otimes \wh U_{\leq -1}$, such that 
$\sum_i \alpha'_i \otimes \lambda'_i = \sum_i (\al''_i \ell_i) \otimes
\lambda''_i - \alpha''_i \otimes (\ell_i\lambda''_i)$. 
\hfill \qed \medskip 

If now $Y$ is a topological $\wh U$-module and $V$ is a 
$\g$-module, the morphism 
$\wh U \to U(\l) \wh\otimes \wh U$ extending the coproduct of $U(\l)$ 
allows to view $Y\wh\otimes V$ as a $\wh U$-module. 
$\wh U$-modules can be constructed as follows: let $\lambda\in\l^*$
be a character such that $D'_0(\lambda)\neq 0$, 
and let $(V,\rho_V)$ be a $\l$-module. 
Then $V((\hbar))$ is a $\wh U$-module, where $x\in\l$ acts as 
$\rho_V(x) + \hbar^{-1} \lambda(x)\id_V$. 

Denote by $\wh{U(\g)}$ the microlocalization of $U(\g)$ 
associated with $D'$. Then $\wh{U(\g)}$ is isomorphic to 
$\big( U(\u_+) \otimes U(\u_-) \big) \wh\otimes \wh U$. Let 
$\wh U(\p_-) \subset \wh{U(\g)}$ be the subalgebra $U(\u_-)
\wh\otimes \wh U$. Then any $\wh U$-module $Y$ may be viewed as a
$\wh{U(\p_-)}$-module. We associate to it the $\wh{U(\g)}$-module 
$\wh Y := \Ind_{\wh{U(\p_-)}}^{\wh{U(\g)}}(Y)$. 

The coproduct of $U(\g)$ also extends to a morphism 
$\wh{U(\g)} \to \wh{U(\g)} \wh\otimes U(\g)$, so if $Z$ is a 
$\wh{U(\g)}$-module and $V$ is a $\g$-module, then $Z\wh\otimes V$ is a 
$\wh{U(\g)}$-module. 

\begin{proposition}  \label{prop:K}
If $Y,Y'$ are $\wh U$-modules and $V$ is a $\g$-module, and if 
$\xi\in \on{Hom}_{\wh U}(Y,Y'\otimes V)$, then there is a
$\wh{U(\g)}$-module morphism $\Phi^\xi : \wh Y \to \wh Y' \wh\otimes V$, 
such that $\Phi^\xi_{|Y} = \xi + $ higher degree terms. 
Set $K = \sum_i a_i \otimes b_i \otimes \ell_i$, 
and set $J := \sum_i a_i \otimes S(b_i)S(\ell_i^{(2)}) \otimes 
S(\ell_i^{(1)})$. If we write 
$J = \sum_i \alpha_i \otimes\beta_i \otimes \lambda_i$, then 
$\Phi^\xi_{|Y} = \sum_i (\alpha_i \otimes\beta_i) \circ \xi \circ 
S(\lambda_i)$. 
\end{proposition} 

{\em Proof.} The properties of $K$ imply that for any $u\in\u_-$, 
$(u^{(1)} + u^{(2)}) J = J u^{(1)} + L$, where $L$ has the form 
$\sum_{i,\alpha} (k'_{i,\alpha} \otimes k''_{i,\alpha} \otimes 1)
(e_\alpha^{(1)} + e_\alpha^{(2)} - e_\alpha^{(3)}) 
(1\otimes 1\otimes k'''_{i,\alpha})$, 
where $(e_\alpha)_\alpha$ is a basis of $\l$. 
\hfill \qed \medskip 

When $Y = \CC$, this proposition shows how to construct singular vectors in
tensor products. 

We now show that $J$ also controls the fusion of intertwiners. 
Let $Y,Y',Y''$ be $\wh U$-modules and let $V',V''$ be $\g$-modules. 
Let $\xi\in \on{Hom}_{\wh U}(Y,Y'\wh\otimes V')$ and $\xi'\in 
\on{Hom}_{\wh U}(Y',Y''\wh\otimes V'')$. Then $(\Phi^{\xi'} \otimes \on{id})
\circ \Phi^{\xi}\in \Hom_{\wh{U(\g)}}(\wh Y,\wh Y'\wh\otimes (V''\otimes V'))$, 
and 
\begin{equation} \label{exp:value}
\langle (\Phi^{\xi'} \otimes \on{id})
\circ \Phi^{\xi}\rangle = \sum_i (\id\otimes a_i \otimes b_i)
\circ (\xi'\otimes \id) \circ \xi \circ S(\ell_i), 
\end{equation}
where $K = \sum_i a_i \otimes b_i \otimes \ell_i$, and $\langle ... \rangle$
denote the component $Y \to Y'\wh\otimes (V''\otimes V')$ of an intertwiner. 
If $A(\xi,\xi')$ is the r.h.s. of (\ref{exp:value}), we even have 
$$
(\Phi^{\xi'} \otimes \id) \circ \Phi^\xi = \Phi^{A(\xi,\xi')}. 
$$
All this follows from the fact that $J$ satisfies the dynamical twist
equation. 

\subsection{Microlocalized Harish-Chandra map}

To state the composition formula, we need microlocalized versions of the 
Harish-Chandra map and the PBW isomorphism, which we now prove. 

Let $\a = \b \oplus \c_+ \oplus \c_-$ be a polarized Lie algebra, 
let $D\in S^d(\a)$ be a nonzero element, such that $D_{|\b^*}\in S^d(\b)$
is nonzero. Let $\wh U_\a$, $\wh U_\b$ be the microlocalizations of 
$U(\a)$, $U(\b)$ w.r.t. lifts of $D,D_{|\b^*}$. 

Define a product on $\big( U(\c_+) \otimes U(\c_-) \big) \wh\otimes 
\wh U_\b$ as follows: 
\begin{align*}
& \mu = (132) \circ 
\big(  m_{U(\c_+)} \otimes m_{\wh U_\b}^{(3)} \otimes m_{U(\c_-)} \big) 
\circ 
(\id \otimes e_+ \otimes \id \otimes e_-^{-1} \otimes \id)
\circ (\id\otimes \id \otimes \pi \otimes \id \otimes \id)
\\ & \circ \big( (132) \otimes (132) \big).  
\end{align*}
Here $m_{A}$ is the product map of an algebra $A$, $m_A^{(3)} : A^{\otimes 3}
\to A$ is $(m_A\otimes \id) \circ m_A$, 
$e_\pm : \wh U_\b \wh\otimes U(\c_\pm) \to U(\c_\pm) 
\wh\otimes \wh U_\b$ are the exchange maps defined as the unique continuous
extensions of $\wh U_\b \otimes U(\c_\pm) \ni f \otimes x \mapsto 
\sum_i x_i \otimes f_i \in U(\c_\pm) \wh\otimes \wh U_\b$, such that 
$f x = \sum_i x_i f_i$ (identity in the microlocalization of $U(\b\oplus
\c_\pm)$ w.r.t. a lift of $D_{|\b^*}$), $\pi : U(\c_-) \otimes U(\c_+)
\to U(\c_+) \otimes U(\b) \otimes U(\c_-)$ is the composition of 
$U(\c_-) \otimes U(\c_+) \to U(\a)$ with the inverse of 
$U(\c_+) \otimes U(\b)\otimes U(\c_-) \to U(\a)$ (both maps are 
inclusions followed by the product of $U(\a)$).

\begin{lemma} \label{PBW:loc}
$\mu$ is an associative, continuous product on 
$\big( U(\c_+) \otimes U(\c_-) \big) \wh\otimes 
\wh U_\b$. The subspace $U(\c_+) \otimes U(\c_-) \otimes U(\b)$
is a subalgebra of $\Big( \big( U(\c_+) \otimes U(\c_-) 
\big) \wh\otimes \wh U_\b, \mu\Big)$, and is isomorphic to $(U(\a),
m_{U(\a)})$ under $\alpha : x_+\otimes x_- \otimes f
\mapsto x_+ f x_-$. 

There is a unique morphism of topological algebras
$\wh{\on{PBW}} : \wh U_\a\to \Big( \big( U(\c_+) \otimes U(\c_-) 
\big) \wh\otimes \wh U_\b, \mu\Big)$, extending the inverse of 
the isomorphism $\alpha$.  
\end{lemma}

{\em Proof.} The associativity of the transport of $m_{U(\a)}$ 
on $U(\c_+) \otimes U(\c_-) \otimes U(\a)$ may be viewed as
a consequence of the commutativity of diagrams involving 
$U(\c_\pm)$ and $U(\b)$. These diagrams still commute when 
$U(\b)$ is replaced by $\wh U_\b$, which implies the associativity 
of $\mu$. 

Let us choose lifts $\wt D$, $\wt D_0$ of $D$ and $D_{|\b^*}$ 
in $U(\a),U(\b)$ is such a way that $\wt D\in U(\a)_{\leq d}$, 
$\wt D_0\in
U(\b)_{\leq d}$ and $H(\wt D) = \wt D_0$, and let us 
construct an inverse of $\alpha(\wt D)$. Set $\xi_0 := \wt D - 
\wt D_0$, 
and define inductively $\xi_n$, $n\geq 0$ by 
$$
\xi_n = - (\id\otimes m^{(2)}_{\wh U_\b} \otimes \id) \circ (e_+ \otimes \id
\otimes e_-^{-1})(\wt D_0^{-1} \otimes \xi_{n-1} \otimes \wt D_0^{-1}). 
$$
The partial degree of $\xi_0$ in $\wh U_\b$ is $\leq d-1$, by construction, 
and the partial degree of $\xi_n$ in $\wh U_\b$ in $\leq d-1 - 2nd$
(because $e_\pm$ has partial degree $0$ for the filtration 
by the $\wh U_\b$-degree; actually its associated graded 
for this filtration is the identity). Therefore the sum 
$1\otimes 1 \otimes \wt D_0^{-1} + \sum_{n\geq 1} \xi_n^{1,3,2}$
converges in $\big( U(\c_+) \otimes U(\c_-) 
\big) \wh\otimes \wh U_\b$, and one shows that it is inverse to 
$\alpha(\wt D)$. The construction of $\wh{\on{PBW}}$ then follows from the 
universal property of Springer's microlocalization. 
\hfill \qed \medskip 

\begin{remark}
Set $\wh H := (\eps\otimes \eps\otimes \id) \circ \wh{\on{PBW}}$, 
then $\wh H : \wh U_\a \to \wh U_\b$ is a continuous map, extending the
Harish-Chandra map $H$. Moreover, $\wh{\on{PBW}}$ can be recovered from 
$\wh H$ using the formula 
$$
\wh{\on{PBW}} = (\pi_+\otimes \wh H \otimes \pi_-)
\circ (\Delta_l \otimes \id) \circ \Delta_r.
$$
Here $\Delta_l : \wh U_\a\to U(\a)\wh\otimes
\wh U_\a$, $\Delta_r : \wh U_\a \to \wh U_\a \wh\otimes U(\a)$ are 
the left- and right-comodule structures of $\wh U_\a$ under $U(\a)$, 
and $\pi_\pm : U(\a) \to U(\c_\pm)$ are the maps 
$U(\a)\to U(\a) \otimes_{U(\b\oplus \c_-)} \CC \to U(\c_+)$, 
$U(\a) \to \CC \otimes_{U(\b\oplus\c_+)} U(\a) \to U(\c_-)$, 
induced by the natural projections and the inverses of the maps 
$x_+\mapsto x_+\otimes 1$, $x_-\mapsto 1\otimes x_-$.  
In particular, $\wh{\on{PBW}}$ is a left $U(\c_+)$-module and 
right $U(\c_-)$-module morphism. 
\end{remark}

\begin{remark}  
If $d\in \ZZ$, let $X_d$ be the subspace of 
$\big( U(\c_+) \otimes U(\c_-) \big) \wh\otimes 
\wh U_\b$, topologically generated by the 
$U(\c_+)_{\alpha} \otimes U(\c_-)_{\beta} \otimes (\wh U_\b)_{\leq 
d-\alpha - \beta}$, 
where $\alpha,\beta\geq 0$. Then $X_d\subset X_{d+1}$, $\mu(X_d \otimes 
X_{d'}) \subset X_{d+d'}$, and $X_d$ is contained in the degree $\leq d$
part of $\big( U(\c_+) \otimes U(\c_-) \big) \wh\otimes 
\wh U_\b$, therefore if we set $X :=
\cup_{d\in\ZZ} X_d$, $X$ is a topological subalgebra of 
$\Big( \big( U(\c_+) \otimes U(\c_-) \big) \wh\otimes 
\wh U_\b,\mu\Big)$. One can check that $\wh{\on{PBW}}$ factors through a 
morphism $\wh U_\a\to X$. 
\end{remark}

\subsection{The composition formula} \label{sect:comp} \label{sect:34}

Assume that $\g = \l \oplus \u_+ \oplus \u_-$ and 
$\l = \k \oplus \m_+ \oplus \m_-$ are nondegenerate polarized Lie 
algebras. 

Set $\vv_\pm = \u_\pm \oplus \m_\pm$. Then $\g = \k \oplus 
\vv_+ \oplus \vv_-$ is a polarized Lie algebra.

\begin{lemma}
$\g = \k \oplus \vv_+ \oplus \vv_-$ is nondegenerate. 
\end{lemma}

{\em Proof.}
Let $d'_\m = \on{dim}(\m_\pm)$, $d'_\u = \on{dim}(\u_\pm)$,
$d'_\vv = d'_\m + d'_\u$, and 
$\underline D^\g_\k\in S^{d'_\vv}(\k)$, 
$\underline D^\g_\l\in S^{d'_\u}(\l)$, 
$\underline D^\k_\l\in S^{d'_\m}(\k)$
be the determinants associated to each polarized Lie algebra. 
Then $[\m_\pm,\u_\mp] \subset \u_\mp$, therefore 
$$
\underline D^\g_\k = (\underline D^\g_\l)_{|\k^*} \underline D^\l_\k.    
$$
Here the map $x\mapsto x_{|\k^*}$ is the algebra morphism 
$S^\cdot(\l) \to S^\cdot(\k)$, taking $x_+x_0 x_-$ to 
$\eps(x_+)\eps(x_-)x_0$, where $x_\pm\in S^\cdot(\m_\pm)$
and $x_0\in S^\cdot(\k)$ 
(it is the associated graded of the Harish-Chandra map, 
and corresponds to the inclusion $\k^* \subset \l^*$
attached to the decomposition of $\l$). Since $\underline D^\g_\k$ and 
$\underline D^\l_\k$ 
are nonzero, so is $\underline D^\g_\l$. 
\hfill \qed \medskip

Let us denote by: 

$\bullet$  $\wh U_\l$, the microlocalization of $U(\l)$
w.r.t. a lift of $\underline D_\l^\g$

$\bullet$ $\wh U_\k$ (resp., $\wh U'_\k$, $\wh U''_\k$), the microlocalization of 
$U(\k)$ w.r.t. a lift of $\underline D_\k^\g$ (resp., $\underline D_\k^\l$, 
$(\underline D_\l^\g)_{|\k^*}$).

\begin{lemma}
1) We have natural inclusions $\wh U'_\k \subset \wh U_\k$, $\wh U''_\k\subset
\wh U_\k$ of complete filtered rings. 

2) $\wh{\on{PBW}}$ is a continuous map 
$\wh U_\l \to \big( U(\vv_+) \otimes U(\vv_-) \big) 
\wh\otimes \wh U''_\k$ of degree $\leq 0$.   
\end{lemma}

{\em Proof.} 1) is clear. 2) follows from Lemma \ref{PBW:loc}. 
\hfill \qed \medskip

Let us denote by $J_\k^\l$, $J_\k^\g$ and $J_\l^\g$ the dynamical twists
associated to the polarized Lie algebras $\l = \k \oplus \m_+\oplus\m_-$, 
$\g = \k \oplus \vv_+ \oplus \vv_-$, $\g = \l \oplus \u_+ \oplus \u_-$. 

We denote by $\eta$ the linear map $\big( U(\u_+) \otimes U(\p_-)\big) 
\wh\otimes \wh U_\l \to \big( U(\vv_+) \otimes U(\p_-)\big) \wh\otimes 
\wh U_\k''$, taking $\alpha\otimes \beta \otimes \lambda$ 
to $\sum_i \sum  \alpha S(\la^{(2)}_{+,i}) \otimes \beta S(\lambda^{(1)}_{+,i})
\otimes \eps(\lambda_{-,i})\lambda_{0,i}$. Here $\alpha\in U(\u_+)$, 
$\beta\in U(\p_-)$, and 
$\wh{\on{PBW}}(\lambda) = \sum \lambda_{-,i} \otimes \lambda_{+,i}
\otimes \lambda_{0,i}$, where $\lambda_{\pm,i}\in U(\m_\pm)$ and 
$\lambda^0\in \wh U_\k$. Here the map $\wh{\on{PBW}}$ is relative to 
the polarization $\l = \k\oplus \m_- \oplus \m_+$, so 
if $\lambda\in U(\l)$, we have $\lambda = \sum_i \lambda_{-,i}
\lambda_{0,i} \lambda_{+,i}$.   

\begin{proposition} \label{prop:comp:formula}
We have $J_\k^\g = \eta(J_\l^\g) J_\k^\l$. This is an equality in 
$\big( U(\vv_+) \otimes U(\p_-)\big) \wh\otimes \wh U_\k$, where 
$\eta(J_\l^\g)$ (resp., $J_\k^\l$) is viewed as an element of this 
algebra using the injection $\wh U''_\k \subset \wh U_\k$ (resp., 
$\wh U'_\k\subset \wh U_\k$). 
\end{proposition}

\begin{remark}
This formula allows one to recover $J_\l^\g$ uniquely from $J_\k^\l,J_\k^\g$.
Indeed, let $\eta' : \big( U(\vv_+) \otimes U(\p_-)\big) \wh\otimes \wh U_\k 
\to \big( U(\u_+) \otimes U(\p_-)\big) \wh\otimes \wh U_\k$ be the map 
taking $u_+\lambda^+ \otimes p_-\otimes k$
to $u_+ \otimes p_- S(\lambda^+)\otimes k$, where $\u_+\in U(\u_+)$, 
$p_-\in U(\p_-)$, $\lambda^+\in U(\m_+)$, 
$k\in U(\k)$, then $\eta' \circ \eta = \id\otimes\id\otimes H_\k^\l$. Now
$J_\l^\g$ can be uniquely recovered from its image by $\id \otimes \id \otimes
H_\k^\l$ using its $\l$-invariance, because the map $S^\cdot(\l)^\l \to
S^\cdot(\k)$, $f\mapsto f_{|\k^*}$ is injective (see Remark \ref{rem12}). 
\end{remark}

\begin{remark}
One can prove that the classical limit of $\eta(J_\l^\g)$ is 
$(r_\l^\g)_{| \k^*}$, so the classical limit of 
Proposition \ref{prop:comp:formula} is 
$r_\k^\g = (r_\l^\g)_{|\k*} + r_\k^\l$ (see Remark \ref{rem:2:3}).  
\end{remark}

{\em Proof.}
We set 
$$
K_\k^\l = \sum_i \alpha_i \otimes \beta_i \otimes \kappa_i 
\in \big( U(\m_+) \otimes U(\m_-)\big) \wh\otimes \wh U'_\k, 
$$
$$
K_\l^\g = \sum_j a_j \otimes b_j \otimes 
c_j^+ c_j^0 c_j^- 
\in \big( U(\u_+) \otimes U(\u_-)\big) \wh\otimes \wh U_\l, 
$$
where $c_j^\pm \in U(\m_\pm)$, $c_j^0\in \wh U_\k''\subset \wh U_\k$. 

Then $J_\k^\l = \sum_i \alpha_i \otimes S(\beta_i) S(\kappa_i^{(2)})
\otimes S(\kappa_i^{(1)})$, and 
$$
J_\l^\g = \sum_j a_j \otimes S(b_j) S(c_j^{-(2)}) S(c_j^{0(2)})
S(c_j^{+(2)}) \otimes S(c_j^{-(1)}) S(c_j^{0(1)}) S(c_j^{+(1)}), 
$$
therefore 
$$
\eta(J_\l^\g) = \sum_j a_j c_j^+ \otimes S(b_j) S(c_j^-) S(c_j^{0(2)})
\otimes S(c_j^{0(1)}), 
$$
and we want to prove that 
$$
J_\k^\g = \sum_{i,j} a_j c_j^+ \alpha_i \otimes S(b_j) S(c_j^-) S(c_j^{0(2)})
S(\beta_i) S(\kappa_i^{(2)}) \otimes S(c_j^{0(1)}) S(\kappa_i^{(1)}), 
$$
i.e., that 
\begin{equation} \label{wanted}
K_\k^\g = \sum_{i,j} a_j c_j^+ \alpha_i 
\otimes S(c_j^{0(2)}) \beta_i c_j^{0(3)} c_j^- b_j 
\otimes \kappa_i c_j^{0(1)}. 
\end{equation}
To prove (\ref{wanted}), we will prove that: 

(a) the r.h.s. of (\ref{wanted}) belongs to $\big( U(\vv_+) \otimes U(\vv_-)
\big) \wh\otimes \wh U_\k$, 

(b) for any $x\in U(\vv_-)$, $y\in U(\vv_+)$, we have 
\begin{equation} \label{id:k:g}
\sum_{i,j} H_\k^\g(x a_j c_j^+ \alpha_i) \kappa_i c_j^{0(1)}
H_\k^\g\big( S(c_j^{0(2)}) \beta_i c_j^{0(3)} c_j^- b_j y\big) 
= H_\k^\g(xy). 
\end{equation}
Here $H_\k^\g$ is the Harish-Chandra map $U(\g) \to U(\k)$. 

Let us now prove (a). We have $a_j\in U(\u_+)$, $c_j^+\in U(\m_+)$, 
$\alpha_i \in U(\m_+)$, so the first factor of the r.h.s. of 
(\ref{wanted}) belongs to $U(\vv_+)$. Since $[\k,\m_-] \subset \m_-$, 
$S(c_j^{0(2)})\beta_i c_j^{0(3)}\in U(\m_-)$; we also have
$c_j^-\in U(\m_-)$ and $b_j\in U(\u_-)$, therefore the second factor of the 
r.h.s. of (\ref{wanted}) belongs to $U(\vv_-)$. 
Finally, since $\kappa_i\in \wh U'_\k$ and $c_j^{0(1)} \in \wh U_\k$, 
the third factor of the r.h.s. of (\ref{wanted}) belongs to 
$\wh U_\k$. This proves (a). 

Let us now prove (b), i.e., identity (\ref{id:k:g}). Since $H_\k^\g$
is a left $U(\k)$-module morphism, $c_j^{0(1)}$ can be inserted in the argument
of $H_\k^\g$, so (\ref{id:k:g}) is equivalent to the identity
\begin{equation} \label{id:k:g:1}
\sum_{i,j}
H_\k^\g(x a_j c_j^+ \alpha_i) \kappa_i H_\k^\g(\beta_i c_j^0 c_j^-
b_j y) = H_\k^\g(xy). 
\end{equation}   
We now prove: 
\begin{lemma}
If $z\in U(\g)$ and $t\in U(\l)$, then 
$$
H_\k^\g(zt) = H_\k^\l\big( H_\l^\g(z) t\big)
\; \on{and} \ 
H_\k^\g(tz) = H_\k^\l \big( t H_\l^\g(z)\big).   
$$
\end{lemma}

{\em Proof of Lemma.} We way assume that $z = z_+ z_0 z_-$, with 
$z_\pm \in U(\u_\pm)$, $z_0\in U(\l)$. Then $H_\k^\g(zt) = 
\eps(z_+) H_\k^\g(z_0 z_- t) = \eps(z_+) \eps(z_-) H_\k^\g(z_0 t)
= H_\k^\l(H_\l^\g(z)t)$. 
The second identity is proved in the same way. 
\hfill \qed \medskip 

It follows that the l.h.s. of (\ref{id:k:g:1}) is equal to 
$$
\sum_{i,j} H_\k^\l \big( H_\l^\g(xa_j c_j^+)\alpha_i\big) \kappa_i 
H_\k^\l\big( \beta_i H_\l^\g(c_j^0 c_j^- b_j y)\big) , 
$$
which is equal to 
\begin{equation} \label{id:k:g:2}
\sum_j H_\k^\l \big( H_\l^\g(xa_j c_j^+) H_\l^\g(c_j^0 c_j^- b_j y)\big) . 
\end{equation}
Now $c_j^+$ and $c_j^0 c_j^-$ belong to $U(\l)$, and $H_\l^\g$ is a 
$U(\l)$-bimodule map, so 
$$
\sum_j H_\l^\g(xa_j c_j^+) H_\l^\g(c_j^0 c_j^- b_j y) 
= \sum_j H_\l^\g(xa_j) c_j^+ c_j^0 c_j^- H_\l^\g(b_j y) =
H_\l^\g(xy). 
$$
Therefore l.h.s. of (\ref{id:k:g:2}) $ = H_\k^\l(H_\l^\g(xy)) = 
H_\k^\g(xy)$. This proves (b). 
\hfill \qed \medskip

\subsection{The ABRR equation} \label{sect:ABRR}

We assume now that $\g = \l\oplus \u_+ \oplus \u_-$ is a polarized Lie 
algebra, equipped with $t\in S^2(\g)^\g$, such that $t$
decomposes as $t = t_\l + s + s^{2,1}$, where $t_\l\in S^2(\l)$ and 
$s\in \u_+ \otimes \u_-$. Then $t_\l$ is $\l$-invariant. We then say that 
$(\g,t)$ is a quadratic polarized Lie algebra. 
We set $\bar s := s^{2,1}$. 

Let $\mu$ be the Lie bracket, and set $\gamma := - {1\over 2}\mu(s)$. Then

\begin{lemma} \label{pties:gamma} \label{3.14}
1) $[\gamma,\l] = 0$. 

2) $[\gamma,\u_\pm] \subset \u_\pm$. 
\end{lemma}

{\em Proof.} Let us prove 1). If $x\in \l$, 
we have $[s,x^{(1)} + x^{(2)}] \in \u_+ \otimes\u_-$, 
$[t_\l,x^{(1)} + x^{(2)}] \in \l\otimes \l$, and $[\bar s,x^{(1)} + x^{(2)}]
\subset \u_- \otimes \u_+$. Since the sum of these terms is 
zero, each of them is zero. Applying $\mu$ to 
$[s,x^{(1)} + x^{(2)}] = 0$, we get 1). 

Let us prove 2). If $x\in \u_+$, we have $[s,x^{(1)}] \in 
\u_+ \otimes \u_-$, $[s,x^{(2)}]\in \u_+ \otimes \g$, 
$[t_\l,x^{(1)}]\in \u_+ \otimes \l$, $[t_\l,x^{(2)}]\in\l\otimes \u_+$, 
$[\bar s,x^{(1)} + x^{(2)}]\in \g\otimes\u_+$. Therefore 
$([s,x^{(2)}])_{\u_+ \otimes \u_-} = - [s,x^{(1)}]$, 
so $[s,x^{(1)} + x^{(2)}]\in \u_+ \otimes \p_+$. Applying $\mu$
to this relation, we get $[\gamma,x] \in \u_+$.
One proves $[\gamma,\u_-] \subset \u_-$ in the same way.      
\hfill \qed\medskip

Assume now that $\g$ is nondegenerate (as a polarized Lie algebra). 

\begin{lemma} 
Let us set $K = \sum_i a_i \otimes b_i \otimes \ell_i$, 
$s = \sum_{\sigma} u_\sigma^+ \otimes u_\sigma^-$, 
$t_\l = \sum_\lambda I_\lambda \otimes I_\lambda$. Then 
$$
\sum_{i,\sigma} u_\sigma^+ a_i \otimes \ell_i \otimes b_i u_\sigma^- 
= \sum_i a_i \otimes \ell_i \otimes [\gamma,b_i]
+  \sum_{i,\lambda} a_i \otimes \ell_i I_\lambda \otimes [I_\lambda,b_i]
- {1\over 2}
\sum_{i,\lambda} a_i \otimes \ell_i \otimes  [I_\lambda,[I_\lambda,b_i]]. 
$$
\end{lemma}

{\em Proof.} Let $\delta$ be the difference of both sides, then $\delta^{1,3,2}$
belongs to $(U(\u_+) \otimes U(\u_-))\wh\otimes \wh U_\l$. Set $\delta = 
\sum_i \delta'_i \otimes \delta''_i \otimes \delta'''_i$, it will suffice to 
prove that 
for any $x,y\in U(\g)$, we have $\sum_i H(x\delta'_i) \delta'''_i 
H(\delta''_i y)= 0$. 

Let $x,y\in U(\g)$, then  
\begin{equation} \label{start}
\sum_{i,\sigma} H(x u_\sigma^+ a_i) \ell_i H(b_i u_\sigma^- y)
= H(x m(s) y), 
\end{equation}
where $m$ is the product map of $U(\g)$. 

Set $C_\g = {1\over 2} m(t_\g)$, $C_\l = {1\over 2} m(t_\l)$. Then 
$C_\g = C_\l + m(s) + \gamma$, so (\ref{start}) is equal to 
$$
H(x(C_\g - C_\l - \gamma)y).
$$ 
Since $C_\g$ is central, this is  
$H(xy C_\g) - H(x(C_\l + \gamma)y)$. Using again $C_\g = C_\l + 
m(s) + \gamma$ and the fact that $H(zm(s)) = 0$ for any $z\in U(\g)$, 
we rewrite (\ref{start}) as $H(xy (C_\l + \gamma)) - H(x(C_\l + \gamma)y)$, 
and since $H$ is a right $U(\l)$-module map, this is  
$$
H(xy) C_\l - H(xC_\l y) + H(x[y,\gamma]).
$$ 
Now we have $
H(xy)C_\l  = 
\sum_{i,\lambda} H(xa_i) \ell_i C_\l  H(b_i y)$,  
since $C_\l$ commutes with $U(\l)$. 
Moreover,  
$$
H(xC_\l y) = {1\over 2}\sum_{i,\lambda} H(xI_\lambda a_i)
\ell_i H(b_i I_\lambda y) 
= {1\over 2}\sum_{i,\lambda}
H(x([I_\lambda,a_i] + a_i I_\lambda)) \ell_i 
H(([b_i,I_\lambda] + I_\lambda b_i)y). 
$$
On the other hand, for any $\xi\in U(\g)$, we have 
$H([\gamma,\xi]) = [\gamma,H(\xi)]$. 
Indeed, if $\xi = \xi^+ \xi^0 \xi^-$, then according to 
Lemma \ref{pties:gamma}, the triangular decomposition of 
$[\gamma,\xi]$ is $[\gamma,\xi^+] \xi^0 \xi^- + 
\xi^+ [\gamma,\xi^0] \xi^- + \xi^+ \xi^0 [\gamma,\xi^-]$, 
and since $\eps([\gamma,\xi^\pm]) = 0$, we get 
$H([\gamma,\xi]) = \eps(\xi^+) [\gamma,\xi^0] \eps(\xi^-) = 
[\gamma,H(\xi)]$. 

Therefore $H(x[y,\gamma]) = \sum_i H(xa_i) \ell_i H(b_i[y,\gamma])
= \sum_i H(xa_i) \ell_i H([\gamma,b_i]y)$. Therefore 
$$
(\ref{start}) = \sum_i H(xa_i) \ell_i C_\l H(b_i y) 
- {1\over 2}\sum_{i,\lambda}
H(x([I_\lambda,a_i] + a_i I_\lambda)) \ell_i 
H(([b_i,I_\lambda] + I_\lambda b_i)y) 
+ \sum_i H(xa_i) \ell_i H([\gamma,b_i]y). 
$$
Therefore 
\begin{align*}
& \sum_{i,\sigma} u_\sigma^+ a_i \otimes \ell_i \otimes 
b_i u_\sigma^- = \sum_i a_i \otimes \ell_i 
\otimes [\gamma,b_i]
\\ & + {1\over 2}\sum_{i,\lambda}
a_i \otimes [\ell_i,I_\lambda] I_\lambda \otimes b_i
+ 
a_i \otimes I_\lambda \ell_i \otimes [I_\lambda,b_i]
+ [a_i,I_\lambda] \otimes \ell_i I_\lambda \otimes b_i
- [a_i,I_\lambda] \otimes \ell_i \otimes [I_\lambda,b_i]. 
\end{align*}
Then we use the $\l$-invariance of $K$ to transform the two last terms. 
\hfill \qed \medskip 

Recall that $J = \sum_i a_i \otimes S(b_i)S(\ell_i^{(2)}) \otimes 
S(\ell_i^{(1)})$. 

\begin{corollary} \label{pln} \label{ABRR} (The nonabelian 
ABRR equation.) We have 
$$
s^{1,2} J = [-\gamma^{(2)} + {1\over 2} m(t_\l)^{(2)} +  t_\l^{2,3},J]. 
$$
\end{corollary}

{\em Proof.} Uses the facts that $t_\l$ commutes with 
$\Delta(U(\l)) \subset U(\l)^{\otimes 2}$ and that $m(t_\l)$
is central in $U(\l)$. 
\hfill \qed \medskip

\begin{remark}
A quadratic polarized Lie algebra $\g$ such that $\l = 0$ and $t$ is
nondegenerate, is the same as a Manin triple, i.e., as a Lie bialgebra 
structure on $\u_+$ (or $\u_-$). Such a polarized Lie algebra is 
degenerate (unless $\g = 0$) and does not lead to a classical dynamical
$r$-matrix. 
\end{remark}

\begin{remark}
Corollary \ref{pln} may be written in a "normally ordered way"
\begin{equation} \label{no:modif}
\big( t_\l^{2,3} - s^{1,2} - m(\bar s)^{(2)} \big) J  = 
J \big( t_\l^{2,3} - m(\bar s)^{(2)} \big).  
\end{equation}
Here "normally ordered" means that all expressions involving $s = \sum_\sigma
u_\sigma^+ \otimes u_\sigma^-$ are such that $u_\sigma^-$ appears before 
$u_\sigma^+$ if both of them are in the same factor. 
\end{remark}

\begin{remark} {\it Expression of the $r$-matrix.}
Assume that $t$ is nondegenerate. 
Let $t^\vee : \g^*\to \g$ be the map $\lambda\mapsto 
(\lambda \otimes \id)(t)$. Then $t^\vee$ is an isomorphism and
restricts to an isomorphism $\l^* \to \l$. So if $\ell$ is
a generic element of $\l$, the bilinear form $\u_+ \times \u_-\to\CC$, 
$(x,y)\mapsto 
\big( (t^\vee)^{-1}(\ell) \big) \big( [x,y] \big) = \ell([x,y])$ 
is nondegenerate. By invariance of the scalar product, 
it follows that for such an $\ell$, 
the operators $\ad(\ell)\in \End(\u_\pm)$ are invertible. 
If we identify $\wedge^2(\g)$ with a subspace of $\End(\g)$
using the scalar product, the $r$-matrix of Proposition \ref{prop:r:mat}
is $\lambda \mapsto {1\over {\ad(t^\vee(\lambda))}} P$, where 
$P$ is the projection on $\u_+\oplus \u_-$ along $\l$ and 
$\on{ad}(t^\vee(\lambda))$ is viewed as an automorphism of 
$\u_+\oplus\u_-$. The same applies in the case of a Lie algebra with a 
splitting and a nondegenerate $t\in S^2(\g)^\g$. 
\hfill \qed \medskip  
\end{remark}

\subsection{Multicomponent ABRR equations}

Here $\g$ is still a quadratic polarized Lie algebra, 
nondegenerate as a polarized Lie algebra. 

\begin{proposition} \label{semi:abrr}
We have 
\begin{equation} \label{ABRR:2} \label{13}
J^{12,3,4} (t_\l^{2,3} + t_\l^{2,4} - s^{1,2} - m(\bar s)^{(2)})
= (t_\l^{2,3} + t_\l^{2,4} - s^{1,2} + s^{2,3} - m(\bar s)^{(2)}) J^{12,3,4}  
\end{equation}
and
\begin{equation} \label{ABRR:1} \label{14}
J^{1,23,4} (t_\l^{3,4} - s^{2,3} - m(\bar s)^{(3)})
= (t_\l^{3,4} - s^{1,3} - s^{2,3} - m(\bar s)^{(3)}) J^{1,23,4} . 
\end{equation}
\end{proposition}

{\em Proof.} 
Let us prove the first identity. 
Recall that $K = \sum_i a_i \otimes b_i \otimes \ell_i$. 
The identity follows from  
\begin{align} \label{ABRR:3}
& \sum_{i,\lambda} a_i^{(1)} \otimes a_i^{(2)} I_\lambda 
\otimes I_\lambda \ell_i \otimes b_i 
+ \sum_{i,\sigma} a_i^{(1)} u_\sigma^+ \otimes a_i^{(2)} u_\sigma^-
\otimes \ell_i \otimes b_i 
+ 
a_i^{(1)} \otimes a_i^{(2)} u_\sigma^- u_\sigma^+
\otimes \ell_i \otimes b_i 
\\ & =  \nonumber 
\sum_{i,\lambda} a_i^{(1)} \otimes I_\lambda 
a_i^{(2)} \otimes \ell_i I_\lambda 
\otimes b_i 
\\ & \nonumber 
+ \sum_{i,\sigma} u_\sigma^+ a_i^{(1)} \otimes u_\sigma^- a_i^{(2)} 
\otimes \ell_i \otimes b_i 
+ a_i^{(1)} \otimes u_\sigma^+ a_i^{(2)} 
\otimes \ell_i \otimes b_i u_\sigma^-
+ a_i^{(1)} \otimes u_\sigma^- u_\sigma^+ a_i^{(2)} 
\otimes \ell_i \otimes b_i 
\\ & \nonumber 
- \sum_{i,\lambda}  a_i^{(1)} \otimes I_\lambda a_i^{(2)} 
\otimes \ell_i \otimes [I_\lambda,b_i]. 
\end{align}
Write the difference of both sides of (\ref{ABRR:3}) as $\sum_i A_i \otimes B_i 
\otimes C_i \otimes D_i$. It will suffice to show that for any $x\in U(\u_+)$, 
we have $\sum_i A_i \otimes B_i \otimes C_i H(D_i x) = 0$. 

This means that for any $x\in U(\u_-)$,  
\begin{align} \label{this:id}
& \sum_\lambda x^{(1)} \otimes x^{(2)} I_\lambda \otimes I_\lambda
+ \sum_\sigma x^{(1)} u_\sigma^+ \otimes x^{(2)} u_\sigma^- \otimes 1
+ x^{(1)} \otimes x^{(2)} u_\sigma^- u_\sigma^+ \otimes 1
\\ & \nonumber 
= \sum_\sigma 
u_\sigma^+ x^{(1)} \otimes u_\sigma^- x^{(2)} \otimes 1 
+ x^{(1)} \otimes u_\sigma^- u_\sigma^+ x^{(2)} \otimes 1
+ \sum_\lambda 
x^{(1)} \otimes I_\lambda x^{(2)} \otimes I_\lambda
+ [I_\lambda,x]^{(1)} \otimes I_\lambda [I_\lambda,x]^{(2)}
\otimes 1
\\ & \nonumber
+ \sum_{\lambda,\sigma} 
A_{\sigma,\lambda}(x)^{(1)} \otimes u_\sigma^+ A_{\sigma,\lambda}(x)^{(2)}
\otimes I_\lambda 
+ \sum_\sigma 
A_{\sigma}(x)^{(1)} \otimes u_\sigma^+ A_{\sigma}(x)^{(2)}
\otimes 1,  
\end{align}
where $A_\sigma, A_{\lambda,\sigma}$ are the linear endomorphisms of
$U(\u_+)$ defined by the condition that 
$u_\sigma^- x - (A_\sigma(x) + \sum_\lambda A_{\sigma,\lambda}(x)I_\lambda)$
belongs to $U(\g)\u_-$. 

This identity decomposes into two parts. The first part is 
\begin{equation} \label{part:1}
x^{(1)} \otimes x^{(2)}I_\lambda 
= x^{(1)} \otimes I_\lambda x^{(2)} + 
\sum_\sigma A_{\sigma,\lambda}(x)^{(1)} 
\otimes u_\sigma^+ A_{\sigma,\lambda}(x)^{(2)}, 
\end{equation}
which we prove as follows: since $x^{(1)} \otimes A_{\sigma,\lambda}(x^{(2)})
= A_{\sigma,\lambda}(x)^{12}$, it suffices to prove that 
$[x,I_\lambda] = \sum_{\sigma} u_\sigma^+ A_{\sigma,\lambda}(x)$.
This collection of identities in $U(\u_+)$ is equivalent to the 
identity in $U(\u_+)\l$, $\sum_\lambda [x,I_\lambda] I_\lambda = 
\sum_{\lambda,\sigma} u_\sigma^+ A_{\sigma,\lambda}(x) I_\lambda$. 
The last identity is proved as follows: we have 
$[m(t),x] = 0$, where $m(t) = \sum_\lambda (I_\lambda)^2 + 
2 \sum_{\sigma} u_\sigma^+ u_\sigma^- + 2\gamma$. This gives
\begin{equation} \label{source}
\sum_{\lambda} 2  [I_\lambda,x] I_\lambda - [[I_\lambda,x],I_\lambda]
+ 2 [\gamma,x] + 2 \sum_{\sigma} u_\sigma^+\big( A_\sigma(x)
+ \sum_{\lambda} A_{\sigma,\lambda}(x) I_\lambda \big) \in U(\g)\u_-
\end{equation} 
which implies $\sum_{\lambda} [I_\lambda,x] I_\lambda
+ \sum_{\sigma,\lambda} u_\sigma^+ A_{\sigma,\lambda}(x)I_\lambda = 0$,
as wanted.  Notice for later use that (\ref{source}) also implies 
\begin{equation} \label{19}
\sum_{\lambda} - [[I_\lambda,x],I_\lambda] + 2 [\gamma,x]
+ 2 \sum_\sigma u_\sigma^+ A_\sigma(x) = 0. 
\end{equation}

The second part of (\ref{this:id}) is 
\begin{equation} \label{20}
[x^{(1)} \otimes x^{(2)},\sum_\sigma u_\sigma^+ \otimes u_\sigma^-
+ 1 \otimes u_\sigma^- u_\sigma^+] 
= \sum_\lambda [I_\lambda,x]^{(1)} \otimes I_\lambda [I_\lambda,x]^{(2)}
+ \sum_\sigma A_\sigma(x)^{(1)} \otimes u_\sigma^+ A_\sigma(x)^{(2)}.
\end{equation}
Before we prove (\ref{20}), we prove 
\begin{equation} \label{21} 
[x,\sum_\sigma u_\sigma^- u_\sigma^+] = \sum_\lambda I_\lambda 
[I_\lambda,x] + \sum_\sigma u_\sigma^+ A_\sigma(x), 
\; \forall x\in U(\u_+).  
\end{equation}
According to (\ref{19}), this is written as 
$$
[x,\sum_\sigma u_\sigma^- u_\sigma^+]
= {1\over 2} \sum_\lambda \big( I_\lambda [I_\lambda,x]
+ [I_\lambda,x]I_\lambda\big) - [\gamma,x],  
$$ 
which follows from the fact that $m(t)$ is central. This implies
(\ref{21}). 

Let us now prove (\ref{20}). The difference between (\ref{20})
and (\ref{19}) applied to the second factor of $x^{(1)} \otimes x^{(2)}$
is
\begin{equation} \label{22}
[x^{(1)} \otimes x^{(2)} , \sum_\sigma u_\sigma^+ \otimes u_\sigma^-]
= \sum_\lambda [I_\lambda,x^{(1)}] \otimes I_\lambda x^{(2)}
+ \sum_\sigma (1\otimes u_\sigma^+)(A_\sigma(x)^{(1)} \otimes A_\sigma(x)^{(2)} 
- x^{(1)} \otimes A_\sigma(x^{(2)})). 
\end{equation}     
So we should prove (\ref{22}). 
 
Since $A_\sigma(x)^{12} = (A_\sigma \otimes \id + \id\otimes A_\sigma)
(x^{12})$, (\ref{22}) is rewritten as 
\begin{equation} \label{23}
[x^{(1)} \otimes x^{(2)} , \sum_\sigma u_\sigma^+ \otimes u_\sigma^-]
= \sum_\lambda [I_\lambda,x^{(1)}] \otimes I_\lambda x^{(2)}
+ \sum_\sigma A(x^{(1)}) \otimes u_\sigma^+ x^{(2)} . 
\end{equation}     
If $x\in \u_+$, then $[x^1 + x^2,\sum_\sigma u_\sigma^+ \otimes u_\sigma^-]
\in \u_+ \otimes \p_+$, so if $x\in U(\u_+)$, then the l.h.s. of (\ref{23})
belongs to $U(\u_+) \otimes \big( U(\u_+) \oplus U(\u_+)\l\big)$. 

Now 
\begin{align*} 
& [x^{(1)} \otimes x^{(2)},\sum_\sigma u_\sigma^+ \otimes u_\sigma^-]  
= - [x^{(1)} \otimes x^{(2)},\sum_\lambda I_\lambda \otimes I_\lambda]
- [x^{(1)} \otimes x^{(2)}, \sum_\sigma u_\sigma^- \otimes u_\sigma^+]
\\ & = 
\sum_\lambda [I_\lambda,x^{(1)}] \otimes I_\lambda x^{(2)}
+ \sum_\lambda x^{(1)} I_\lambda \otimes [I_\lambda,x^{(2)}]
+ \sum_\sigma u_\sigma^- x^{(1)} \otimes u_\sigma^+ x^{(2)}
\; \on{modulo}\; U(\g)\u_- \otimes U(\u_+)
\\ & = 
\sum_\lambda [I_\lambda,x^{(1)}] \otimes I_\lambda x^{(2)}
+ \sum_\lambda x^{(1)} I_\lambda \otimes [I_\lambda,x^{(2)}]
\\ & 
+ \sum_{\sigma} \big( \sum_\lambda A_{\sigma,\lambda}(x^{(1)}) I_\lambda
+ A_\sigma(x^{(1)})\big)  \otimes u_\sigma^+ x^{(2)}
\; \on{modulo}\; U(\g)\u_- \otimes U(\u_+).
\end{align*}
Projecting this identity on $U(\u_+) \otimes U(\g)$ parallel to 
$U(\u_+)\l \otimes U(\g)$, we get (\ref{23}). 

Let us now prove the second identity. Using $[\sum_{\sigma} u_\sigma^+ 
\otimes u_\sigma^-,x^1 + x^2] = 0$ for $x\in\l$ and the 
$\l$-invariance of $K$, we transform this identity 
into the analogue of (\ref{this:id}) with $\u_+,\u_-$ exchanged, 
which also holds. 
\hfill \qed \medskip

Corollary \ref{ABRR} has a multicomponent version. Namely, let 
$$
J^{[n]} := 
J^{1,\{2,...,n\},n+1} \cdots J^{n-2,\{n-1,n\},n+1}J^{n-1,n,n+1}
$$
(this element corresponds to fusing $n$ intertwiners). 
\footnote{Here, for example, 
$J^{1,\{2,...,n\},n+1}$ means that we put the first component 
of $J$ in component 1, the second in components $2...n$ 
(after taking the coproduct $n-2$ times), and the third in component $n+1$.} 

\begin{theorem}\label{ABRRM:bis} (The multicomponent ABRR equation)
For $i=1,\ldots,n$, the element $J^{[n]}$ satisfies 
the equations
\begin{equation} \label{multi:ABRR}
[\sum_{j=i+1}^{n+1} t_\l^{i,j} - \gamma^{(i)}
+\frac{1}{2}m(t_\l)^{(i)},J^{[n]}]=
(\sum_{j=1}^{i-1} s^{j,i} - \sum_{j = i+1}^n s^{i,j}) J^{[n]}.
\end{equation}
\end{theorem}

{\em Proof.} We will treat the case $n=3$. Then $J^{[3]} = 
J^{1,23,4}J^{2,3,4} = J^{12,3,4}J^{1,2,34}$. Then 
\begin{align*}
& 
(t_\l^{3,4} - s^{1,3} - s^{2,3} - m(\bar s)^{(3)}) J^{[3]}
= 
(t_\l^{3,4} - s^{1,3} - s^{2,3} - m(\bar s)^{(3)}) J^{1,23,4} J^{2,3,4}
\\ & 
= 
J^{1,23,4}
(t_\l^{3,4} - s^{2,3} - m(\bar s)^{(3)}) J^{2,3,4}
= 
J^{1,23,4} J^{2,3,4}
(t_\l^{3,4} - m(\bar s)^{(3)}) 
\\ & = 
J^{[3]}
(t_\l^{3,4} - m(\bar s)^{(3)}) ,  
\end{align*}
where we have used (\ref{14}) and (\ref{no:modif}). This proves 
(\ref{multi:ABRR}) when $i=3$. 

Let us treat the case $i=2$. 
\begin{align*}
& 
(t_\l^{2,3} + t_\l^{2,4} - s^{1,2} + s^{2,3} - m(\bar s)^{(2)}) J^{[3]}
= 
(t_\l^{2,3} + t_\l^{2,4} - s^{1,2} + s^{2,3} - m(\bar s)^{(2)}) 
J^{12,3,4} J^{1,2,34}
\\ & 
= 
J^{12,3,4}
(t_\l^{2,3} + t_\l^{2,4} - s^{1,2} - m(\bar s)^{(2)}) J^{2,3,4}
= 
J^{12,3,4} J^{1,2,34}
(t_\l^{2,3} + t_\l^{2,4} - m(\bar s)^{(2)}) 
\\ & = 
J^{[3]}
(t_\l^{2,3} + t_\l^{2,4} - m(\bar s)^{(2)}),  
\end{align*}
where we have used (\ref{13}) and $(\ref{no:modif})^{1,2,34}$. This proves 
(\ref{multi:ABRR}) when $i=2$.

In general, (\ref{multi:ABRR}) for $i=1$ is a consequence of 
(\ref{multi:ABRR}) for $i = 2,\ldots,n$, and of the $\l$- and 
$\gamma$-invariances of $J$, and of $[\gamma,\l] = 0$. 
We have already proven the $\l$-invariance of 
$J$, and its $\gamma$-invariance follows from that of $K$, 
which in its turn follows from the identity $H([\gamma,x]) = [\gamma,H(x)]$ 
for $x\in U(\g)$.   
\hfill \qed \medskip 

\begin{proposition} \label{ABRR:compat} (Compatibility of multicomponent ABRR) 
Write the multicomponent ABRR equations as 
$a_i^{[n]} J^{[n]} = J^{[n]} b_i^{[n]}$, for $i = 2,\ldots,n$. 
This is a compatible system, i.e., $[a_i^{[n]},a_j^{[n]}] = 
[b_i^{[n]},b_j^{[n]}] = 0$ for any pair $i,j \in \{2,\ldots,n\}$.  
\end{proposition}

{\em Proof.} The vanishing of these brackets follows from the identities
\begin{equation} \label{mod:CYB}
[s^{1,2},s^{1,3}] + [s^{1,2},s^{2,3}] + [s^{1,3},s^{2,3}] 
= [t_\l^{2,3},s^{1,3}]
\end{equation}
and 
\begin{equation} \label{to:prove}
[s^{1,2},m(\bar s)^{(1)} + m(\bar s)^{(2)} - t_\l^{1,2}] = 0. 
\end{equation}
To prove (\ref{mod:CYB}), one may assume that $t$ is nondegenerate. 
Both sides of (\ref{mod:CYB}) belong to $\u_+ \otimes \g\otimes \u_-$. 
Then (\ref{mod:CYB}) follows from its pairing with $x_- \otimes \id 
\otimes x_+$, with $x_\pm\in \u_\pm$. 

Let us prove (\ref{to:prove}). Let us project $[\gamma^{(1)} 
+ \gamma^{(2)},t] = 0$ on $\u_+ \otimes\u_-$. Lemma \ref{3.14}
says that this projection is $[\gamma^{(1)} + \gamma^{(2)},s] = 0$. 
(\ref{to:prove}) then follows from this identity, together with 
$m(\bar s) = {1\over 2} (m(t) - m(t_\l)) + \gamma$, the fact that 
$m(t)$ is central in $U(\g)$, and $[s,m(t_\l)^{12}] = 0$, 
which follows from the $\l$-invariance of $s$.   
\hfill \qed \medskip

\section{Dynamical pseudotwists associated to a quadratic 
polarized Lie algebra} \label{sect:pseudo}

As we noted in the Introduction, Proposition \ref{prop:comp:r:matrices}
together with \cite{AM1} implies: 

\begin{lemma} \label{lemma:rho:c}
Let $(\g = \l\oplus\u,t)$ be a quadratic Lie algebra 
with a nondegenerate splitting. 
Let $c\in\CC$, then $\rho_c \in\wedge^2(\g) \otimes 
\wh S^\cdot(\l)[1/D_0]$ defined by 
$$
\rho_c(\lambda) := r_\l^\g(\lambda) 
+ c \big( f(c\on{ad}\lambda^\vee) \otimes \id\big) (t)
$$
for $\lambda\in \l^*$, is a solution of $\on{CYB}(\rho_c) 
+ \on{Alt}(\on{d}\rho_c) = -\pi^2 c^2 Z$. Here we set 
$\lambda^\vee = (\lambda\otimes\id)(t)$ 
and $f(x) = - 1/x + \pi{\on{cotan}}(\pi x)$. 
\end{lemma}

In this section, we assume that $(\g = \l\oplus \u_+ \oplus \u_-,
t\in S^2(\g)^\g)$ is a quadratic polarized Lie algebra, such that 
$\g$ is nondegenerate as a polarized Lie algebra (see Section
\ref{sect:ABRR}). Recall that this means that $t$ decomposes as $t_\l + s +
s^{2,1}$, with $t_\l\in S^2(\l)$ and $s\in\u_+ \otimes \u_-$. 
We will construct a dynamical pseudotwist quantizing
$\rho_c$ in this situation. 

We fix a formal parameter $\hbar$ and a complex 
parameter $c$. We set $\kappa = \hbar c$.  
If $\Phi(A,B)$ is a Lie associator (\cite{Dr2}), we set 
$\Phi_\kappa^{-1}(A,B) = \Phi(\kappa A,\kappa B)^{-1}$. 
An example of an associator is 
the KZ associator, i.e., the renormalized holonomy from $0$ to $1$ 
of the differential equation ${{dF}\over{dz}} = ({A\over z} + 
{B\over {z-1}})F$. 

\begin{theorem} \label{thm:barK}
Set 
$$
\bar J = \Phi_\kappa^{-1}\big( t^{1,2},
t_\l^{2,3} - s^{1,2} - m(\bar s)^{(2)}\big) 
J_\l^\g. 
$$
Then $\bar J\in U(\g)^{\otimes 2} \wh\otimes \wh U_\l[[\hbar]]$
is a solution of the dynamical pseudotwist equation
\begin{equation} \label{pseudotwist}
\bar J^{12,3,4} \bar J^{1,2,34} = \Phi_\kappa^{-1}(t^{1,2},t^{2,3})
\bar J^{1,23,4} \bar J^{2,3,4}. 
\end{equation}
\end{theorem}

{\em Proof.} Drinfeld's algebra $\cT_4$ is defined by generators 
$\tau_{i,j}$, $1\leq i\neq j\leq 4$, and relations 
$\tau_{i,j} = \tau_{j,i}$, $[\tau_{i,j},\tau_{k,l}] = 0$
if $\{i,j,k,l\} = \{1,2,3,4\}$, and 
$[\tau_{i,j} + \tau_{i,k},\tau_{j,k}] = 0$ if card$\{i,j,k\} = 3$.
Then we have the pentagon relation
$$
\Phi_\kappa^{-1}(\tau^{1,2},\tau^{2,3})
\Phi_\kappa^{-1}
(\tau^{1,2} + \tau^{1,3}, \tau^{2,4} + \tau^{3,4})
\Phi_\kappa^{-1}(\tau^{2,3},\tau^{1,2})
= 
\Phi_\kappa^{-1}(\tau^{1,3} + \tau^{2,3}, \tau^{3,4})
\Phi_\kappa^{-1}(\tau^{1,2},\tau^{2,3} + \tau^{2,4}). 
$$

Then we have an algebra morphism $\cT_4\to U(\g)^{\otimes 3}
\otimes U(\l)$, with $\tau^{i,j} \mapsto t^{i,j}$ for $i,j\neq 4$, 
$\tau^{2,4} \mapsto t_\l^{2,4} - s^{1,2} - s^{3,2} - m(\bar s)^{(2)}$, 
$\tau^{3,4} \mapsto t_\l^{3,4} - s^{1,3} - s^{2,3} - m(\bar s)^{(3)}$, 
$\sum_{1\leq i <j\leq 4} \tau_{i,j} \mapsto 0$. 

Taking the image of the pentagon relation by this morphism, 
we get an identity in $U(\g)^{\otimes 2} \wh\otimes \wh U_\l[[\hbar]]$. 
Multiply it from the right by the identity 
$(J_\l^\g)^{1,23,4} (J_\l^\g)^{2,3,4} = 
(J_\l^\g)^{12,3,4} (J_\l^\g)^{1,2,34}$. Then using the identities 
of Proposition \ref{semi:abrr} to put $(J_\l^\g)^{1,23,4}$ before 
the image of $\Phi_\kappa^{-1}(\tau^{2,3},\tau^{3,4})$ in the l.h.s., and 
$(J_\l^\g)^{12,3,4}$ before 
the image of 
$\Phi_\kappa^{-1}(\tau^{1,2},\tau^{2,3} + \tau^{2,4})$ in the r.h.s., 
we get the result. 
\hfill \qed \medskip

Let us study the classical limit of $\bar J$. In 
Section \ref{sect:quant}, we introduced quasi-commutative algebras
$\wh S^\cdot(\l)_\hbar \subset \wh S^\cdot(\l)[1/D_0]_\hbar$, 
and inclusions $\wh S^\cdot(\l)_\hbar \subset U(\l)[[\hbar]]$, 
$(\wh U_\l)_{\leq 0} \subset \wh S^\cdot(\l)[1/D_0]_\hbar$. 

Then $\bar J$ belongs to $\cA := U(\G)^{\otimes 2} \bar\otimes 
\wh S^\cdot(\l)[1/D_0]_\hbar$ (here $\bar\otimes$ is 
the "formal series" tensor product). 

\begin{proposition} (Classical limit.)
$\bar J-1$ belongs to $\hbar \cA$, and the reduction of 
$\on{Alt}_{1,2}(\bar J -1)/\hbar$ modulo $\hbar$ belongs to 
$\wedge^2(\g) \otimes \wh S^\cdot(\l)[1/D_0]$. It coincides
with the expansion at origin of the meromorphic function 
$\rho_c : \l^* \to \wedge^2(\g)$, defined in Lemma \ref{lemma:rho:c}. 
\end{proposition}

{\em Proof.} It will be enough to compute the classical limit of 
$X := \Phi_\kappa^{-1}(t^{1,2},t_\l^{2,3} - s^{1,2} - m(\bar s)^{(2)})$.
We have $\Phi^{-1} = \exp(\phi)$, where $\phi$ is a Lie formal series in 
$A,B$. $X$ is the specialization of $\Phi$ under 
$A\mapsto A_0 := \hbar c t^{1,2}$, $B\mapsto B_0 
:= \hbar c(t_\l^{2,3} - s^{1,2} - m(\bar s)^{(2)})$. 
We have $A_0\in \hbar \cA$ and $B_0\in \cA$, moreover if
$B'_0 := \hbar c t_\l^{2,3}$, then $B_0 = B'_0$ modulo $\hbar \cA$. 

Since $\phi$ is a Lie series, it is a sum of homogeneous components 
with partial degree $\geq 1$
in $A$ and $B$, so $\phi(A_0,B_0) \in \hbar\cA$. In particular, 
$X \in 1 + \phi(A_0,B_0) + \hbar^2 \cA$. Moreover, $\phi(A_0,B_0)
= \phi(A_0,B'_0) = \phi_1(A_0,B'_0)$ modulo $\hbar^2 \cA$, where 
$\phi_1$ is the part of $\phi$, linear in $A$. 

Set $\phi_1 = \sum_{k\geq 1} c_k \on{ad}(B)^k(A)$. 
Since $\wh S^\cdot(\l)[D_0^{-1}]_\hbar$ is quasi-commutative, 
$\hbar^{-1}\on{ad}(B'_0)^k (A_0) \in \cA$ is equal 
modulo $\hbar \cA$, to 
$c^{k+1} 
\sum_{\alpha,\lambda_1,\ldots,\lambda_k} 
e_\alpha \otimes \on{ad}(I_{\lambda_1})
\cdots \on{ad}(I_{\lambda_k})(e_\alpha)
\otimes (\hbar I_{\lambda_1}) \cdots (\hbar I_{\lambda_k})$ 
(here $t = \sum_\alpha
e_\alpha \otimes e_\alpha$). 

Its reduction modulo $\hbar\cA$ is the linear map $\l^*\to \g^{\otimes 2}$, 
$\lambda \mapsto c^{k+1} \ad(1\otimes t^\vee(\lambda))^k(t^{1,2})$. 
This map takes values in $S^2(\g)$ if $k$ is even and in $\wedge^2(\g)$
if $k$ is odd. Only the "even $k$" part remains after antisymmetrization, 
and the result follows from the fact that the $c_{2k}$
coincide with the Taylor coefficients of $f$ (see e.g. \cite{EE}).  
\hfill \qed \medskip

Let $\eta : U(\g)^{\otimes 3} \to U(\g)^{\otimes 2} \otimes U(\l)$
be the linear map associated to $\g = \l\oplus\u_+ \oplus \u_-$
(see Section \ref{sect:comp}). 

\begin{proposition} \label{prop:4:3}
If $P$ is any noncommutative polynomial in two variables, we have 
$$
\eta(P(t^{1,2},t^{2,3})) = P(t^{1,2},t_\l^{2,3} - s^{1,2} - m(\bar s)^{(2)}). 
$$
In particular, $\eta(\Phi_\kappa^{-1}(t^{1,2},t^{2,3})) = 
\Phi_\kappa^{-1}(t^{1,2},t_\l^{2,3} - s^{1,2} - m(\bar s)^{(2)})$.  
\end{proposition}

{\em Proof.} We have for any $x\in U(\g)^{\otimes 3}$, 
\begin{equation} \label{stat:1}
\eta(t^{1,2} x) = t^{1,2} \eta(x).
\end{equation}

Let us now show that if $x\in U(\g)^{\otimes 3}$ is $\g$-invariant, 
then 
\begin{equation} \label{stat:2}
\eta(t^{2,3} x) = (t_\l^{2,3} - s^{1,2} - m(\bar s)^{(2)}) \eta(x). 
\end{equation}
Let us write $x = \sum_i A_i \otimes B_i \otimes \lambda_i^-
\lambda_i^0\lambda_i^+$. Since $\eps(u_\sigma^- \lambda_i^-) = 0$, 
we have $\eta(s^{2,3} x) = 0$. 

We have 
$t_\l^{2,3} x = \sum_{i,\lambda} A_i \otimes 
I_\lambda B_i \otimes [I_\lambda,\lambda_i^-] \lambda_i^0 \lambda_i^+
+ \sum_{i,\lambda} A_i \otimes 
I_\lambda B_i \otimes \lambda_i^- (I_\lambda\lambda_i^0) \lambda_i^+$, 
so since $\eps([I_\lambda,\lambda_i^-]) = 0$, we get 
$$
\eta(t_\l^{2,3} x) = \sum_{i,\lambda}
A_i S(\lambda_i^{+(2)}) \otimes I_\lambda B_i S(\lambda_i^{+(1)})
\otimes \eps(\lambda^-) I_\lambda \lambda_i^0 = t_\l^{2,3} \eta(x).  
$$ 

Finally, since $x$ is invariant, we have 
$$
\bar s^{2,3} x = \sum_{i,\sigma} A_i \otimes u_\sigma^- B_i \otimes 
\lambda_i^- \lambda_i^0 (\lambda_i^+ u_\sigma^+)
+ [A_i,u_\sigma^+] \otimes u_\sigma^- B_i \otimes \lambda_i^- \lambda_i^0
\lambda_i^+
+ A_i \otimes u_\sigma^- [B_i,u_\sigma^+] \otimes \lambda_i^- \lambda_i^0
\lambda_i^+, 
$$
therefore 
\begin{align*}
&  \eta(\bar s^{2,3}x) = \sum_{i,\sigma}
- A_i S(\lambda_i^{+(2)}) \otimes u_\sigma^- B_i u_\sigma^+ S(\lambda_i^{+(1)})
\otimes \eps(\lambda_i^-) \lambda_i^0
\\ & 
- A_i u_\sigma^+ S(\lambda_i^{+(2)}) \otimes u_\sigma^- B_i S(\lambda_i^{+(1)})
\otimes \eps(\lambda_i^-) \lambda_i^0
+ [A_i, u_\sigma^+] S(\lambda_i^{+(2)}) \otimes u_\sigma^- B_i S(\lambda_i^{+(1)})
\otimes \eps(\lambda_i^-) \lambda_i^0
\\ & + A_i S(\lambda_i^{+(2)}) \otimes u_\sigma^- [B_i, u_\sigma^+] S(\lambda_i^{+(1)})
\otimes \eps(\lambda_i^-) \lambda_i^0
\\ & 
= - \sum_\sigma (u_\sigma^+ \otimes u_\sigma^- \otimes 1
+ 1 \otimes u_\sigma^- u_\sigma^+ \otimes 1) \eta(x). 
\end{align*}
Adding up these results, we get (\ref{stat:2}). The proposition now follows 
from (\ref{stat:1}), (\ref{stat:2}) and $\eta(1) = 1$. 
\hfill \qed \medskip 

\begin{remark}
The relation $\eta\big(\Phi^{\on{KZ}}_\kappa(t^{1,2},t^{2,3})^{-1}\big) = 
\Phi^{\on{KZ}}_\kappa(t^{1,2},t_\l^{2,3} - s^{1,2} - m(\bar s)^{(2)})^{-1}$, 
where $\Phi^{\on{KZ}}$ is the KZ associator, can be derived from the 
results of Section \ref{connection:ABRR} 
(in the untwisted case) together with the  
composition formula.  
\hfill \qed \medskip 
\end{remark}

\begin{remark} \label{rem:connection}
If $J\in U(\g)^{\otimes 2}[[\hbar]]$ satisfies 
$J^{12,3} J^{1,2} = \Phi_{\kappa}^{-1}(t^{1,2},t^{2,3})
J^{1,23}J^{2,3}$, then $U_\hbar(\g) := 
(U(\g)[[\hbar]],\on{Ad}(J^{-1}) \circ \Delta_0)$ is a Hopf 
algebra. 

Assume that $\bar J \in U(\g)^{\otimes 2}\wh\otimes \wh U_\l[[\hbar]]$
satisfies the pseudotwist equation (\ref{pseudotwist}), and set 
$J_\hbar := (J^{1,2})^{-1} \bar J$, then $J_\hbar$ satisfies the twist 
equation in $U_\hbar(\g)$
\begin{equation} \label{twist:quantum}
J_\hbar^{12,3,4} J_\hbar^{1,2,34} = 
J_\hbar^{1,23,4} J_\hbar^{2,3,4}.  
\end{equation}
\end{remark}

\begin{remark} 
When $\g$ is a semisimple Lie algebra and $\l\subset\g$ is a 
Cartan subalgebra, $c = 1$, $Z = {1\over 4}[t^{1,2},t^{2,3}]$, 
we have 
$$
\rho_c(\lambda) = - {1\over 2} \sum_{\alpha\in\Delta_+}
(e_\alpha \wedge f_\alpha) \on{coth}{{(\lambda,\alpha)}\over 2}. 
$$
Here $x\wedge y = x\otimes y - y \otimes x$. 
Then a quantization of $\rho_c(\lambda)$ is 
$\bar J := \Phi_{\hbar}^{-1}(t^{1,2},t_\l^{2,3} - r^{1,2} - m(\bar r)^{(2)})
J_\h^\g$ (here $r$ is the standard $r$-matrix of $\g$). 
Therefore, if $J$ is as in Remark \ref{rem:connection}, 
then $J_\hbar := (J^{1,2})^{-1}\bar J$ satisfies (\ref{twist:quantum}). 
On the other hand, we know from \cite{EV2} another solution 
$J'_\hbar$ of the same equation, obtained by a quantum analogue of
the construction of $J_\l^\g$. It is natural to conjecture that 
$J_\hbar$ and $J'_\hbar$ are gauge-related in $U_\hbar(\g)$
(or, which is the same, that $\bar J$ and $J^{1,2} J'_\hbar$ 
are gauge-related in $U(\g)[[\hbar]]$). 
\end{remark}

\begin{remark} {\it (Expression of the $r$-matrix.)} 
If $t\in S^2(\g)^\g$ is nondegenerate and we use it for identifying
$\wedge^2(\g)$ with a subspace of $\on{End}(\g)$, then $\rho_c$
identifies with 
$$
\lambda \mapsto {1\over {\on{ad} t^\vee(\lambda)}} P 
+ c f(c \on{ad}( t^\vee(\lambda) )).   
$$
\end{remark}

\section{Dynamical pseudotwists associated to a quadratic Lie algebra 
with an automorphism} \label{sect:twists}

In this section, we quantize the $r$-matrix $\rho_{\sigma,c}$ of Proposition 
\ref{lemma:0:3}, as well as $(\rho_{\sigma,c})_{|\k^*} + r_\k^\l$, if 
$(\l = \k\oplus \m_+ \oplus \m_-,t_\l)$ is a polarized quadratic Lie algebra, 
nondegenerate as a polarized Lie algebra. 

\subsection{Quadratic Lie algebras with an automorphism} \label{sect:51}

Let $\g$ be a finite dimensional complex Lie algebra, 
equipped with $t\in S^2(\g)^\g$ and $\sigma\in \on{Aut}(\g)$, 
such that $(\sigma\otimes\sigma)(t) = t$. We assume that 
$\sigma - \on{id}$ is invertible on $\g / \g^\sigma$. 

Set $\l := \g^\sigma$, $\u := \on{Im}(\sigma - \id)$. Then 
$\g = \l\oplus \u$ is a Lie algebra with a splitting (see Section 
\ref{splitted}). 
Moreover, we have $t\in S^2(\l) \oplus S^2(\u)$. 
We denote by $t_\l$ the component of $t$ in $S^2(\l)$. We have 
$t_\l\in S^2(\l)^\l$. 

\begin{example}\label{ex1} $\g,\l$ are as in Example \ref{2.5.2} (1), and  
$\sigma= \on{exp}(\on{ad}(\chi))$,
where $\chi$ is a generic central element
of $\l$ (see Section 2).
\end{example}

\begin{example}\label{ex2} 
$\g$ is a simple Lie algebra, $\sigma$ is an involution of $\g$. 
Then $G/L$ is a symmetric space for $G$. 
\end{example}

\begin{example} \label{ex3}
$\g$ is a simple Lie algebra, 
and $\l\subset \g$ is a semisimple Lie subalgebra of the same rank 
(``a Borel-De Siebenthal pair'' \cite{BS}).
In this case there is an automorphism (of degree $2,3$ or $5$)
such that $\l=\g^\sigma$. 
\end{example}

\begin{example} \label{ex4}
$\g$ is a simple Lie algebra of simply laced type, 
$\sigma$ is induced by a Dynkin diagram automorphism
(with no fixed edges). Then 
$\l=\g^\sigma$ is the Lie algebra corresponding to the quotient diagram. 
\end{example}

\begin{example} \label{ex5}
More generally, in the setting of Example \ref{ex4}, 
one may consider the automorphism $\sigma_\beta=
\sigma \circ \on{exp} (\on{ad}(\beta))$, where $\beta$ is a generic  
element of $\h^\sigma$. In this case, $\l=\h^\sigma$. 
\end{example}

\medskip 

Let ${\bf D}:= \{a+ib\in\CC | a\in [0,1[$ and $b\geq 0,$ or 
$a\in ]0,1]$ and $b<0\}$. There is a unique operator $\log(\sigma)$
in $\End(\g)$, such that $e^{\log(\sigma)} = \sigma$, and whose eigenvalues 
belong to $2\pi i {\bf D}$. 

\subsection{The main result: definition and properties of $\Psi_\kappa$}
\label{sect:52}

Let us define $\cO_{]0,1[}$ as the ring of analytic functions on 
$]0,1[$. We define in the same way $\cO_{\RR_+^\times - \{1\}}$. 
Set 
$$
X(z) := {{(z^{\log(\sigma)/2\pi i} \otimes \id)(t-t_\l)}\over{z-1}}
+ {z\over{z-1}} t_\l. 
$$
Then $X(z) \in \g^{\otimes 2} \otimes \cO_{\RR_+^\times - \{1\}}$. 
More precisely, the first term of $X(z)$ is a linear combination 
of products of powers of $\log(z)$ (of degree $<\on{dim}(\g)$) with 
$z^\alpha/(z-1)$, where $\alpha$ is an eigenvalue of $\log(\sigma)/2i\pi$.  

Let $\kappa$ be a formal parameter and let $\Psi_\kappa$ be the 
renormalized holonomy from $0$ to $1$ of the equation 
\begin{equation} \label{basic:eqn}
z {{\on{d}G}\over{\on{d}z}} = \kappa \big( X(z)^{1,2} 
+ t_\l^{2,3} + {1\over 2} m(t_\l)^{(2)}\big) G(z),  
\end{equation}
where $G\in U(\g)^{\otimes 2} \otimes U(\l)\otimes \cO_{]0,1[}[[\kappa]]$. 

More precisely, if $\sigma$ has no strictly positive eigenvalues on
$\g/\g^\sigma$, 
there are unique solutions $G_0,G_1$ of (\ref{basic:eqn}), 
such that 
$G_0(z) = z^{ \kappa( t_\l^{2,3} + {1\over 2} m(t_\l)^{(2)} )}
(1 + o(1))$ 
as $z\to 0$, 
$G_1(z) = (1-z)^{\kappa t^{1,2}}  
(1 + o(1))$ as $z\to 1$, 
and $\Psi_\kappa \in U(\g)^{\otimes 2}
\otimes U(\l)[[\kappa]]$ is defined by $\Psi_\kappa = G_1(z)^{-1} G_0(z)$
for any $z$. If $\sigma$ has strictly positive eigenvalues, 
$G_0,G_1,\Psi_\kappa$ are defined similarly, replacing $[0,1]$ by a
smooth path 
$\gamma : [0,1] \to \CC$, 
such that $\gamma(0) = 0$, $\gamma(1) = 1$, 
$0,1 \notin \gamma(]0,1[)$, and the Euclidean scalar product of 
$\gamma'(0)$ with the eigenvalues of $\log(\sigma)$ is $>0$. 

\begin{theorem} \label{thm:pseudo}
$\Psi_\kappa$ satisfies the pseudotwist equation
$$
((\Phi^{\on{KZ}}_\kappa)^{-1})^{1,2,3} \Psi_\kappa^{1,23,4}
\Psi_\kappa^{2,3,4} = \Psi_\kappa^{12,3,4} \Psi_\kappa^{1,2,34}. 
$$
\end{theorem}

\begin{proposition} \label{prop:class}
(Classical limit of $\Psi_\kappa$.)
Let $c$ be a complex number, $\hbar$ a formal variable, and
assume that $\kappa = \hbar c$. Recall that $\wh S^\cdot(\l)_\hbar$
is the $\hbar$-adically complete subalgebra of $U(\l)[[\hbar]]$
generated by $\hbar \l$. Then $\Psi_\kappa$ and 
$(\Psi_\kappa-1)/\hbar$ belong to $U(\g)^{\otimes 2} \bar\otimes 
\wh S^\cdot(\l)_\hbar$. 

Moreover, the reduction modulo $\hbar$ of 
$(\Psi_\kappa - \Psi_\kappa^{2,1}) / \hbar$ is the formal 
function $\rho_{\sigma,c} : \l^* \to \wedge^2(\g)$, such that 
$$
\rho_{\sigma,c}(\lambda) = 
c \Big( \big( f(c\on{ad}(\lambda^\vee))\otimes\id \big)(t_\l)
+ i\pi \big( {{ e^{2\pi i c \on{ad}(\lambda^\vee)} \circ \sigma + \id }
\over  {e^{2\pi i c \on{ad}(\lambda^\vee)} \circ \sigma - \id }}
\otimes \id\big)(t-t_\l)
\Big) ; 
$$
here for $\lambda\in\l^*$, $\lambda^\vee = (\lambda\otimes 
\id)(t_\l) \in \l$, and  
$f(x) = -{1\over x} + \pi \on{cotan}(\pi x)$. 
\end{proposition}

\begin{remark}
This solution of the modified CDYBE was discovered in \cite{AM2}, 
generalizing \cite{ES2}, where $\sigma$ is assumed of finite order. In the case of Example \ref{ex5}, this solution was 
discovered in \cite{S} and quantized using quantum groups 
in \cite{ESS}. Our quantization is different; it should be related to 
the quantization of \cite{ESS} by a gauge transformation 
given by a twisted version of the Kazhdan-Lusztig equivalence
 between the representation categories of an affine algebra and a 
quantum group.  
\end{remark}

\begin{remark} When $\sigma = \id$, $\g$ is semisimple and $\k$
is a Cartan subalgebra, (\ref{basic:eqn}) is the trigonometric KZ equation, 
see \cite{EFK,EV3}. \hfill \qed \medskip 
\end{remark}

Assume now that $\l$ is a quadratic polarized Lie algebra, nondegenerate as a
polarized Lie algebra. So 
$\l = \k\oplus \m_+ \oplus \m_-$, and $t_\l = t_\k + s + s^{2,1}$, with 
$t_\k \in S^2(\k)^\k$ and $s\in \m_+ \otimes \m_-$. We set $\gamma = 
- {1\over 2} \mu(s)$. 
Let $\Psi_{\k,\l,\g}$
be the renormalized holonomy from $0$ to $1$ of the differential 
equation 
\begin{equation} \label{eq:G} \label{31}
z{{\on{d}G}\over{\on{d}z}} = \kappa \Big( \big( X(z) - s\big)^{1,2}
+ \big( t_\k^{2,3} + {1\over 2} m(t_\k)^{(2)} - \gamma^{(2)} 
\big) \Big) G(z). 
\end{equation}
Then $\Psi_{\k,\l,\g} \in U(\g)^{\otimes 2} \otimes U(\l)[[\kappa]]$. 
The map $\eta$ defined in Section \ref{sect:34} restricts to 
$$
\eta : U(\g)^{\otimes 2} \otimes U(\l) \to U(\g)^{\otimes 2} \otimes 
U(\k).
$$ 

\begin{proposition} \label{prop:eta} \label{prop:59}
1) Set $\Psi_\l^\g:= \Psi_\kappa$, then 
$\eta(\Psi_\l^\g) = \Psi_{\k,\l,\g}$. 

2) Set $\Psi_\k^\g := \Psi_{\k,\l,\g} J_\k^\l$. Then $\Psi_\k^\g
\in U(\g)^{\otimes 2} \wh\otimes \wh U_\k[[\hbar]]$ is a dynamical 
pseudotwist, i.e., it satisfies  
$$
((\Phi^{\on{KZ}}_\kappa)^{-1})^{1,2,3} (\Psi_\k^\g)^{1,23,4} 
(\Psi_\k^\g)^{2,3,4}
= (\Psi_\k^\g)^{12,3,4} (\Psi_\k^\g)^{1,2,34} . 
$$
Moreover, if $\kappa = \hbar c$, 
$\Psi_\k^\g$ belongs to $U(\g)^{\otimes 2} \bar\otimes 
\wh S^\cdot(\k)[1/D_\k]_\hbar$, and its classical limit is 
$\tilde \rho_{c,\k,\g} := r_\k^\l(\lambda) 
+ (\rho_{\sigma,c})_{|\k^*}(\lambda)$.  It 
satisfies the modified CDYBE $\on{Alt}(\on{d} \tilde\rho_{c,\k,\g}) + 
\on{CYB}(\tilde\rho_{c,\k,\g}) = - \pi^2 c^2 Z$.  

Here $r_\k^\l(\lambda)$ is as in Section \ref{sect:rat}
and $D_\k \in S^{\on{dim}(\m_+)}(\k)$ is the determinant 
corresponding 
to the polarized Lie algebra $\l = \k\oplus \m_+ \oplus \m_-$.
\end{proposition}

The proofs of Theorem \ref{thm:pseudo}, Proposition 
\ref{prop:class} and Proposition \ref{prop:eta} 
occupy the rest of this section. 

\begin{remark} \label{rem:J}
If we assume that $\g^\sigma = 0$, then $J := \Psi_\l^\g$ is a solution of 
the twist equation $J^{12,3} J^{1,2} = \Phi_\kappa(t^{1,2},t^{2,3})
J^{1,23} J^{2,3}$, so it gives rise to a quasitriangular Hopf algebra
$$
(U(\g)[[\kappa]],m_0,\Delta := \Ad(J^{-1}) \circ \Delta_0, 
R:= (J^{-1})^{2,1} e^{\kappa t/2} J).
$$
Its classical limit is the quasitriangular Lie bialgebra 
structure on $\g$ induced by the $r$-matrix 
$r:= ({{\id + \sigma}\over{\id - \sigma}} \otimes \id)(t)$
($r$ is antisymmetric, and is a solution of the modified CYBE). 
\end{remark}

\subsection{Proof of Theorem \ref{thm:pseudo}}

Consider the system of equations 
\begin{equation} \label{syst:1}
z {{\partial G}\over{\partial z}} = \kappa 
\Big( X(z)^{1,2} + X(z/u)^{3,2} + \big( t_\l^{2,4} 
+ {1\over 2} m(t_\l)^{(2)} \big) \Big) G, 
\end{equation}
\begin{equation} \label{syst:2}
u {{\partial G}\over{\partial u}} = \kappa 
\Big( X(u)^{1,3} +  X(u/z)^{2,3} + \big( t_\l^{3,4} 
+ {1\over 2} m(t_\l)^{(3)}\big) \Big) G, 
\end{equation}
where the unknown function $G(z,u)$ lies in 
$U(\g)^{\otimes 2} \otimes U(\l) \otimes \cO [[\kappa]]$, 
$\cO$ is the ring of analytic functions on $\{(z,u) | 0<z<u<1\}$, 
and $G$ has the form $1 + O(\kappa)$. 

One checks that the system (\ref{syst:1},\ref{syst:2}) is compatible; 
more generally, 
the following is true. Let $\g = \l \oplus \u$ be a Lie algebra with a 
splitting, equipped with $t\in S^2(\g)^\g$, such that $t = t_\l + t_\u$, 
$t_\l\in S^2(\l)$, $t_\u \in S^2(\u)$. Assume that $\ell\in \End(\u)$
commutes with the adjoint action of $\l$ on $\u$, and that 

(a) $(\ell\otimes \id + \id \otimes \ell)(t_\u) = t_\u$, 

(b) $[t^{1,2},t^{2,3}] = [t_\l^{1,2},t_\l^{2,3}] + \tau_1 + \tau_2$, 
where $(\bar\ell^{(1)} + \bar\ell^{(2)} + \bar\ell^{(3)})(\tau_a) = \tau_a$, 
where $\bar\ell\in\End(\g)$ is defined by $\bar\ell_{|\l} = 0$, 
$\bar\ell_{|\u} = \ell$. 

Then the system (\ref{syst:1},\ref{syst:2}), where $X(z) = 
(z^{\ell} \otimes \id)(t_\u) / (z-1) + t_\l \cdot z / (z-1)$, 
is compatible. 
 
The system (\ref{syst:1},\ref{syst:2}) has therefore a solution, 
unique up to right multiplication 
by an element of $U(\g)^{\otimes 2} \otimes U(\l)[[\kappa]]$ of the form 
$1 + O(\kappa)$. 

Following \cite{Dr2}, we consider five asymptotic zones, corresponding 
to the parenthesis orders $P_1 = ((0z)u)1$, $P_2 = (0(zu))1$,
$P_3 = 0((zu)1)$, $P_4 = (0z)(u1)$,  $P_5 = 0(z(u1))$. 

Assume for simplicity that $\sigma$ has no strictly positive eigenvalue.
This guarantees that the function $z^\lambda$, $\lambda$ in the spectrum of 
$\log(\sigma)$, tends to zero as $z\to 0^+$. 

There exist five solutions of the system 
(\ref{syst:1},\ref{syst:2}) $G_1,\ldots,G_5$ corresponding to these
zones. They are uniquely determined by the requirements
$$
G_1(z,u) = z^{\kappa (t_\l^{2,4} + {1\over 2} m(t_\l)^{(2)})}
u^{ \kappa (t_\l^{2,3} + t_\l^{3,4} + {1\over 2} m(t_\l)^{(3)})} 
(1 + g_1(z,u)),
$$
$$
G_2(z,u) = ({{u-z}\over{u}})^{\kappa t^{2,3}} 
z^{\kappa \big( t_\l^{2,4} + t_\l^{3,4}
+ t_\l^{2,3} + {1\over 2} m(t_\l)^{(2)} + {1\over 2} m(t_\l)^{(3)} \big)}
(1 + g_2(z,u)) , 
$$   
$$
G_3(z,u) = (u-z)^{\kappa t^{2,3}} (1-z)^{\kappa(t^{1,2} + t^{1,3})}
(1 + g_3(z,u)), 
$$
$$ 
G_4(z,u) = z^{\kappa(t_\l^{2,4} + {1\over 2} m(t_\l)^{(2)})}
(1-u)^{\kappa t^{1,3}} (1 + g_4(z,u)), 
$$
$$
G_5(z,u) = (1-u)^{\kappa t^{1,3}} (1-z)^{\kappa(t^{1,2} + t^{2,3})}
(1 + g_5(z,u)), 
$$
where $g_1(z,u)$ (resp., $g_2,g_3,g_4,g_5$)
tends to zero as $(u,z/u)$ (resp., $(u,1 - {z\over u})$, 
$(1-z,{{u-z}\over{1-z}})$, $(z,1-u)$, $(1-z,{{1-u}\over{1-z}})$) 
tends to $(0,0)$ in $]0,1[^2$. Here "tends to zero" means that 
these are zeries in $\kappa$ of tensor products of elements
of $U(\g)^{\otimes 2} \otimes U(\l)$ with analytic functions 
in $(z,u)$ tending to zero in the relevant zone. 
The expansions of $G_2$ and $G_3$ are based on the identity 
$$
z{{\partial G}\over{\partial z}} + u {{\partial G}\over{\partial u}}
= \kappa \big( X(z)^{1,2} + X(u)^{1,3} + 
t_\l^{1,2} + t_\l^{1,3} + t_\l^{2,3} + 
{1\over 2} m(t_\l)^{(2)} + {1\over 2} m(t_\l)^{(3)}\big) G.  
$$

We have the relations 
$$
G_1 = G_2 (\Psi_\l^\g)^{3,2,4}, \; 
G_2 = G_3 (\Psi_\l^\g)^{1,23,4}, \; 
G_3 = G_5 \Phi_\kappa^{\on{KZ}}(t^{2,3},t^{1,3}), 
$$
$$
G_1 = G_4 (\Psi_\l^\g)^{1,3,24}, \; 
G_4 = G_5 (\Psi_\l^\g)^{13,2,4}. 
$$
Therefore $G_1 = G_5 \Phi_\kappa^{\on{KZ}}(t^{2,3},t^{1,3})
(\Psi_\l^\g)^{1,23,4} (\Psi_\l^\g)^{2,3,4}
= G_5 (\Psi_\l^\g)^{13,2,4}(\Psi_\l^\g)^{1,3,24}$. 
Simplifying by $G_5$, exchanging factors $2$ and $3$ 
and using the antisymmetry relation $(\Phi_\kappa^{\on{KZ}})^{3,2,1} = 
(\Phi_\kappa^{\on{KZ}})^{-1}$, we get the result. 
 
If $\sigma$ has strictly positive eigenvalues, one applies the same argument
after replacing the segment $[0,1]$ by a smooth path $\gamma : [0,1] \to \CC$, 
as in the definition of $\Psi_\kappa$. 
Then if $\alpha$ is any eigenvalue of $\log(\sigma)$, 
$z^\alpha \to 0$ as $z\to 0$ along $\gamma$.  
\hfill \qed \medskip

\subsection{Proof of Proposition \ref{prop:eta}}

Let us prove 1). Similarly to Proposition \ref{prop:4:3}, one shows that 
if $P$ is any noncommutative polynomial in $\dim(\g) + 1$ 
variables, and $\ell := \log(\sigma)/2\pi i$, then 

\begin{align*}
& \eta \big( P((\ell^k \otimes \id)(t)^{1,2}, k= 0,\ldots,\dim(\g) -1 | 
t_\l^{2,3} + {1\over 2}
m(t_\l)^{(2)})\big)
\\ & =  
P\big( (\ell^k \otimes \id)(t)^{1,2}, k= 0,\ldots,\dim(\g) -1 | 
t_\k^{2,3} - s^{1,2} - m(\bar s)^{(2)} + m(t_\l)^{(2)} \big) 
\\ & = 
P((\ell^k \otimes \id)(t)^{1,2}, k = 0,\ldots,\dim(\g) - 1 
| t_\k^{2,3} - s^{1,2}+ {1\over 2}
m(t_\k)^{(2)} - \gamma^{(2)}). 
\end{align*} 
$\Psi_\l^\g$ may be expressed as $P((\ell^k \otimes \id)(t)^{1,2}, 
k= 0,\ldots,\dim(\g) -1 | t_\l^{2,3} + {1\over 2} m(t_\l)^{(2)})$; 
this implies 1). 

Let us prove 2). Let us denote by $\Psi'$ and $\Psi''$ the renormalized
holonomies from $0$ to $1$ of the differential equations
\begin{equation} \label{eq:G'}
z {{\on{d}G'}\over{\on{d}z}} = \kappa \big( X(z)^{2,3} - s^{1,3} - s^{2,3} 
+ t_\k^{3,4} + {1\over 2} m(t_\k)^{(3)} - \gamma^{(3)} \big) G', 
\end{equation} 
\begin{equation} \label{eq:G''}
z {{\on{d}G''}\over{\on{d}z}} = \kappa \big( X(z)^{1,2} - s^{1,2} - s^{3,2} 
+ t_\k^{2,4} + {1\over 2} m(t_\k)^{(2)} - \gamma^{(2)} \big)G''.  
\end{equation} 
We will prove the identities
\begin{equation} \label{id:a}
\Phi_\kappa^{\on{KZ}}(t^{1,2},t^{2,3})^{-1}
\Psi_{\k,\l,\g}^{1,23,4} \Psi' = 
\Psi_{\k,\l,\g}^{12,3,4} \Psi''
\end{equation} 
and
\begin{equation} \label{id:b}
(J_\k^\l)^{1,23,4} \Psi' =  \Psi_{\k,\l,\g}^{2,3,4}
(J_\k^\l)^{1,23,4}, \; 
(J_\k^\l)^{12,3,4} \Psi'' =  \Psi_{\k,\l,\g}^{1,2,34}
(J_\k^\l)^{12,3,4}. 
\end{equation}

Then combining (\ref{id:a}), (\ref{id:b}) and the twist equation 
$(J_\k^\l)^{1,23,4} (J_\k^\l)^{2,3,4} = 
(J_\k^\l)^{12,3,4} (J_\k^\l)^{1,2,34}$, we get 2).  

Let us prove (\ref{id:a}). Proposition \ref{semi:abrr} implies that if $G(z)$ 
is a solution of (\ref{eq:G}) of the form $1 + O(\kappa)$, then 
$G'(z) := \big( (J_\k^\l)^{1,23,4}\big)^{-1} G(z)^{2,3,4} (J_\k^\l)^{1,23,4}$
is a solution of (\ref{eq:G'}) of the form $1 + O(\kappa)$, and 
$G''(z) := \big( (J_\k^\l)^{12,3,4}\big)^{-1} G(z)^{1,2,34} 
(J_\k^\l)^{12,3,4}$ is a solution of (\ref{eq:G''}) of the form $1 + O(\kappa)$.
This implies (\ref{id:a}).  

Let us prove (\ref{id:b}). We consider the system of equations 
\begin{equation} \label{syst:3}
z{{\partial G}\over{\partial z}} = \kappa \Big(
X(z)^{1,2} + X(z/u)^{3,2} - s^{1,2} - s^{3,2} + t_\k^{2,4} 
+ {1\over 2} m(t_\k)^{(2)} - \gamma^{(2)}  \Big) G, 
\end{equation}
\begin{equation} \label{syst:4}
u{{\partial G}\over{\partial u}} = \kappa \Big( 
X(u)^{1,3} + X(u/z)^{2,3} - s^{1,3} - s^{2,3} + t_\k^{3,4} + {1\over 2}
m(t_\k)^{(3)} - \gamma^{(3)} 
\Big) G, 
\end{equation}
where the unknown function $G(z,u)$ lies in $U(\g)^{\otimes 2}
\otimes U(\k) \otimes \cO[[\kappa]]$ and has the form $1 + O(\kappa)$. 

As before, the system (\ref{syst:3},\ref{syst:4}), supplemented with the
condition $G = 1 + O(\kappa)$, has a solution, unique up to right
multiplication by an element of $U(\g)^{\otimes 2} \otimes U(\k)[[\kappa]]$
of the form $1 + O(\kappa)$. 

The system (\ref{syst:3},\ref{syst:4})  has unique solutions
$G_1,\ldots,G_5$ corresponding to the asymptotic zones $P_1,\ldots,P_5$, 
satisfying 
$$
G_1(z,u) = z^{\kappa 
\big( - s^{1,2} - s^{3,2} + t_\k^{2,4} + {1\over 2} m(t_\k)^{(2)} - \gamma^{(2)}
\big)} u^{\kappa \big(
- s^{1,3} - s^{2,3} + t_\k^{2,3} + t_\k^{3,4} + {1\over 2} m(t_\k)^{(3)} -
\gamma^{(3)}  \big) } (1 + g_1(z,u)),  
$$ 
$$
G_2(z,u) = \big( {{u-z}\over{u}}\big)^{\kappa t^{2,3}} 
z^{\kappa \big( 
- s^{1,2} - s^{1,3} + t_\k^{2,3} + t_\k^{2,4} + t_\k^{3,4}
+ {1\over 2} m(t_\k)^{(2)} + {1\over 2} m(t_\k)^{(3)} - \gamma^{(2)} -
\gamma^{(3)}  \big)  } (1 + g_2(z,u)), 
$$
$$
G_3(z,u) = (u-z)^{\kappa t^{2,3}} (1-z)^{\kappa(t^{1,2} + t^{1,3})}
(1 + g_3(z,u)), 
$$
$$
G_4(z,u) = z^{\kappa \big( 
- s^{1,2} - s^{3,2} + t_\k^{2,4} + {1\over 2} m(t_\k)^{(2)} - \gamma^{(2)} 
\big) } (1-u)^{\kappa t^{1,3}} (1 + g_4(z,u)), 
$$
$$
G_5(z,u) = (1-u)^{\kappa t^{1,3}} (1-z)^{\kappa(t^{1,2} + t^{2,3})}
(1 + g_5(z,u)), 
$$
where $g_i(z,u) \to 0$ in the zone $P_i$. Then we have 
$$
G_1 = G_2 (\Psi')^{1,3,2,4}, \; 
G_2 = G_3 (\Psi_{\k,\l,\g})^{1,23,4}, \; 
G_3 = G_5 \Phi_\kappa^{\on{KZ}}(t^{2,3},t^{1,3}), 
$$
$$
G_1 = G_4 (\Psi'')^{1,3,2,4}, \; G_4 = G_5 (\Psi_{\k,\l,\g})^{13,2,4}. 
$$
As before, this implies (\ref{id:b}), and therefore 2). 
\hfill \qed \medskip 

\subsection{Classical limits}

Let us prove Proposition \ref{prop:class}. Set $H(z) := 
G(z) z^{- \kappa (t_\l^{2,3} + {1\over 2} m(t_\l)^{(2)})}$, then 
$\Psi_\kappa = \on{lim}_{z\to 1^-} \big( (1-z)^{-\kappa t^{1,2}} H(z) \big)$. 

Let us prove that $H(z)$ has the following $\kappa$-adic property: 
it belongs to $U(\g)^{\otimes 2} \wh\otimes 
\wh S^\cdot(\l) \otimes \cO_{]0,1[}[[\kappa]]$.
Set $\bar t_\l := \kappa t_\l$, then $\bar t_\l \in U(\g)^{\otimes 2}
\wh\otimes \wh S^\cdot(\l)$, and $H(z)$ satisfies the equations 
$$
H(0) = 1, \quad  z{{dH}\over{dz}} = 
\kappa X(z)^{1,2} H(z) + [\bar t_\l^{2,3} + {\kappa \over 2}
m(t_\l)^{(2)},H(z)]. 
$$ 
The formal expansion of $H(z)$ therefore belongs to 
$U(\g)^{\otimes 2} \wh\otimes \wh S^\cdot(\l)[\log(z), 
z^a,a\in {\bf D}][[z]][[\kappa]]$, and has the form $1 + O(\kappa)$
(${\bf D}$ is defined in Section \ref{sect:52}).  
Set $H(z) = 1 + \kappa h(z) + O(\kappa^2)$, then 
$$
h(0) = 0, \quad 
z {{dh} \over{dz}} = X(z)^{1,2} + [\bar t_\l,h(z)]. 
$$
View $h(z)$ as a formal function $\l^* \to U(\g)^{\otimes 2}$, then 
$$
h(0,\lambda) = 0, \quad 
z {{\partial h(z,\lambda)}\over{\partial z}} = 
X(z)^{1,2} + [t^\vee_\l(\lambda)^{(2)},h(z,\lambda)]. 
$$
Here $t_\l^\vee(\lambda) := (\lambda\otimes \id)(t_\l)$. 

Therefore 
$$
h(z,\lambda) = \int_0^z \on{Ad}\big( (z/u)^{t^\vee_\l(\lambda)^{(2)}}\big)
(X(u)^{1,2}) {{du}\over u}.  
$$
It follows from the form of $H(z)$ that $\Psi_\kappa\in 
U(\g)^{\otimes 2} \wh\otimes \wh S^\cdot(\l)[[\kappa]]$, and it has the form $1
+ \kappa\psi + O(\kappa^2)$. We have 
$$
\psi(\lambda) = \on{lim}_{z\to 1^-} \Big( 
\on{Ad} \big( (z/u)^{t^\vee_\l(\lambda)^{(2)}}\big)
\big( {{ (u^{\log(\sigma)/2\pi i} \otimes \id)(t-t_\l)}\over
{u-1}} + {{u}\over{u-1}} t_\l\big)  {{du}\over u}
- t^{1,2} \log(1-z)
\Big) . 
$$
Therefore $\psi(\lambda) \in \g^{\otimes 2} \otimes \wh S^\cdot(\l)$. 

We now use the fact that for ${\rm Re}(x)>0$, one has
$$
\lim_{z\to 1}
\biggl(\int_0^z u^{x-1}\frac{du}{1-u}
-\log (1-z)\biggr)=\frac{1}{x}+\sum_{n\ge 1}(\frac{1}{x+n}-\frac{1}{n})=
-\Gamma'(x)/\Gamma(x).
$$

Using the $\l$-invariance of $(u^{\log(\sigma)/2\pi i} \otimes \id)(t-t_\l)$
and of $t_\l$, and the fact that 
$\log(1-z)
( z^{ -\on{ad} t_\l^\vee(\lambda) } - \id)(t) \to 0$ as $z\to 1^-$, we get 
$$
\psi(\lambda) = \Big( {\Gamma'\over\Gamma} \big( {{\log(\sigma)}\over{2\pi i}} + 
\on{ad} t_\l^\vee(\lambda) \big) \otimes \id\Big)(t-t_\l) 
+ \Big( {\Gamma'\over\Gamma} \big( 1 + \ad t_\l^\vee(\lambda)\big) 
\otimes \id\Big)(t_\l). 
$$

Then using $(\log(\sigma) \otimes \id + \id\otimes \log(\sigma))(t-t_\l) 
= 2\pi i (t-t_\l)$, and the $\l$-invariance of $t-t_\l$ and $t_\l$, 
we get 
\begin{align*}
\psi(\lambda) - \psi(\lambda)^{2,1} 
& = (\Big( 
{\Gamma'\over\Gamma} \big({{\log(\sigma)}\over{2\pi i}}
+ \ad t_\l^\vee(\lambda)\big)
- {\Gamma'\over\Gamma} \big( 1 - {{\log(\sigma)}\over{2\pi i}}
- \ad t_\l^\vee(\lambda)\big)
\Big) \otimes\id)(t-t_\l)
\\ & 
+ (\Big( 
{\Gamma'\over\Gamma} \big( 1 + \ad t_\l^\vee(\lambda)\big)
- {\Gamma'\over\Gamma} \big( 1 - \ad t_\l^\vee(\lambda)\big)
\Big) \otimes\id)(t_\l). 
\end{align*}

Using the identities ${\Gamma'\over\Gamma}(1-x) - {\Gamma' \over \Gamma}
(x) = - \pi \on{cotan}(\pi x)$, and ${\Gamma'\over\Gamma}(x+1)
- {\Gamma'\over\Gamma}(x) = {1\over x}$, we obtain Proposition \ref{prop:class}.

The classical part of Proposition \ref{prop:eta} now follows from the fact 
that the classical counterpart of $\eta$ is the restriction to 
$\k^* \subset \l^*$. \hfill \qed \medskip

\subsection{Twists by an element of $Z(\l)$} \label{sect:Z(l)}
\label{sect:shifts}

One checks that the results of Section \ref{sect:52} can be 
generalized as follows. 

Let $(\g,\l,\sigma)$ be as in Section \ref{sect:51}.  
Let us denote by $Z(\l)$ the center of $\l$ and let $\gamma'\in Z(\l)$. 
Denote by $\Psi_{\kappa,\gamma'}$ the renormalized holonomy 
from $0$ to $1$ of the equation 
\begin{equation}
z {{\on{d}G}\over{\on{d}z}} = 
\kappa \big( X(z) + t_\l^{2,3} + {1\over 2} m(l_\l)^{(2)} 
- \gamma^{\prime (2)}\big) G(z), 
\end{equation}
Then $\Psi_{\kappa,\gamma'}$ satisfies the pseudotwist equation 
$$
((\Phi_{\on{KZ}})^{-1})^{1,2,3} \Psi_{\kappa,\gamma'}^{1,23,4}
\Psi_{\kappa,\gamma'}^{2,3,4}
= \Psi_{\kappa,\gamma'}^{12,3,4} \Psi_{\kappa,\gamma'}^{1,2,34}, 
$$
and its classical limit if $\rho_{\sigma,c}(\lambda)$. 

To prove the first statement, one modifies the system of equations
(\ref{syst:1},\ref{syst:2}) by adding $-\kappa \gamma^{\prime (2)}G$ 
in the r.h.s. 
of (\ref{syst:1}), and $-\kappa \gamma^{\prime (3)}G$ in the r.h.s. of 
(\ref{syst:2}). This is again a compatible system, because of the identity 
$X(z) + X(z^{-1})^{2,1} = t_\l$.  

Assume in addition that $\l$ is quadratic polarized as in the sequel of
Section \ref{sect:52}, let $\Psi_{\kappa,\l,\g}(\gamma')$ be the 
renormalized homolomy from $0$ to $1$ of (\ref{31}), modified by the 
addition of $-\kappa \gamma^{\prime (2)}G$ in its r.h.s. Then 
$\eta(\Psi_{\kappa,\gamma'}) = \Psi_{\k,\l,\g}(\gamma')$, 
and $\Psi_{\k,\l,\g}(\gamma')J_\l^\l$ is a dynamical pseudotwist, 
quantizing $\tilde \rho_{c,\k,\l}$. To prove this, one modifies the system 
(\ref{syst:3},\ref{syst:4}) as above.

\subsection{Relation with twisted loop algebras} \label{connection:ABRR}
\label{sect:5:6}

Here, we interpret results of Section \ref{sect:52} in terms of 
the ABRR equations and the dynamical twist for a twisted loop algebra. 
More precisely, we show that the compatibility of the systems
(\ref{syst:1},\ref{syst:2}) and (\ref{syst:3},\ref{syst:4}) are consequences
of the compatibility of multicomponent ABRR equations (Proposition
\ref{ABRR:compat}), and relate $G(z)$ with a dynamical twist. 

Throughout the section, we assume that $\g,t,\sigma,\l,\k,\m_\pm$ 
are as in Sections \ref{sect:51}, \ref{sect:52}. We also assume that 
$t \in S^2(\g)^\g$ is nondegenerate. 

If $s\in\CC^\times$ is an eigenvalue of $\sigma$, let $\g_s\subset\g$ 
be the generalized eigenspace (we set $\g_s = 0$ for other $s\in\CC^\times$). 
Then $t$ decomposes as a sum $\sum_{s\in\CC^\times} t_s$, where $t_s\in
\g_s\otimes \g_{s^{-1}}$.  

\subsubsection{Twisted loop algebras} \label{4:6:1}

Let us say that a function of one variable $x$ 
is a generalized trigonometric polynomial if 
it is a linear combination of functions of the form 
$x^ne^{ax}$, $n\in {\Bbb Z}_+$, $a\in \Bbb C$
(the sum may involve different $a$). 
Let $L_\sigma\g$ be the Lie algebra of 
 $\g$-valued generalized trigonometric polynomials of 
$x$ satisfying the condition
$$
f(x+1)=\sigma(f(x)).
$$
For notational convenience 
we will express such functions as multivalued functions of 
$z=e^{2\pi ix}$. We will denote by $\CC[\log(z),z^a,a\in\CC]$ the 
ring of generalized trigonometric polynomials.

Set $e(x) = e^{2\pi i x}$, $\Gamma := \{\alpha\in\CC | e(\alpha)$ is an 
eigenvalue of $\sigma\}$, $\CC_+ =\{u + i v | u>0$ or $(u=0$ and $v\geq 0)\}$, 
$\CC_- = - \CC_+$, $\Gamma_\pm = \Gamma \cap \CC_\pm$. If $\lambda,\mu\in\CC$,
we write $\lambda \leq \mu$ (resp., $\lambda < \mu$) iff 
$\mu - \lambda\in\CC_+$ (resp., $\CC_+ - \{0\}$). 

The operator $z{{\on{d}}\over{\on{d} z}}$ acts on $L_\sigma\g$, and 
its eigenvalues belong to $\Gamma$. If $\alpha\in \Gamma$, 
we denote by $(L_\sigma\g)_\alpha$ the corresponding generalized eigenspace. 
Then $L_\sigma\g = \oplus_{\alpha\in\Gamma}(L_\sigma\g)_\alpha$. 

$(L_\sigma\g)_\alpha$ is the subspace of $L_\sigma\g$ of all elements $u$, 
which can be expressed as $u = \sum_i a_i \otimes f_i$, where 
$a_i \in\g_{e(\alpha)}$ and $f_i \in z^\alpha \CC[\on{log}(z)]$. 
Then if $u  = \sum_i a_i \otimes f_i\in (L_\sigma\g)_\alpha$ and 
$v = \sum_j b_j \otimes g_j \in (L_\sigma\g)_{-\alpha}$,  the function 
$\sum_{i,j} (a_i,b_j) f_i g_j(z)$ is in $\CC[\on{log}(z)]$ and is 
invariant under $x\mapsto x+1$, and is therefore constant
(we denote by $(-,-)$ the pairing on $\g$ inverse to $t$). This defines a
nondegenerate pairing $(L_\sigma\g)_\alpha \times (L_\sigma\g)_{-\alpha}
\to\CC$. We denote by $(-,-)$ the direct sum of these pairings, which 
is a nondegenerate pairing $(L_\sigma\g)^2 \to\CC$. 

Let $C_\g$ be the endomorphism of $\g$ equal to $\on{ad}(m(t))$
(we denote by $\on{ad} : U(\g) \to \on{End}(\g)$ the algebra morphism 
extending the Lie algebra morphism $\g \to \on{End}(\g)$ induced by 
the adjoint action of $\g$). 
Then $C_\g \otimes z {{\on{d}}\over{\on{d}z}}$ is a derivation of 
$\g \otimes \CC[\log(z),z^a,a\in\CC]$, which restricts to a derivation of
$L_\sigma\g$.  

Set $\omega(x,y) = {1\over 2}((C_\g \otimes z{{\on{d}}\over{\on{d}z}})(x),y)$
for any $x,y\in L_\sigma\g$. Then $\omega$ is a cocycle on $L_\sigma\g$,
independent on a rescaling of $t$, and is a generalization of the critical level
cocycle (which corresponds to $\g$ simple, $\sigma = \id$). 
  
Define the affine Lie algebra 
$$
\bar\g = \overline{L_\sigma\g} = L_\sigma\g \oplus \CC k \oplus \CC {\bf 1}
\oplus \CC d \oplus \CC \delta,  
$$
where 
$$ 
[a(z),b(z)]=[a,b](z)+
(z a'(z), b(z)) k + \omega (a(z),b(z)) {\bf 1}, 
$$ 
$$
[d,a(z)]=za'(z), \; [\delta,a(z)] = {1\over 2}C_\g(za'(z)), 
\; [d,\delta] = 0,
$$
$k$ and ${\bf 1}$ are central. If $\g$ is semisimple, $\bar\g$ is closely
related to a (possibly twisted) affine Kac-Moody algebra. 

An invariant, nondegenerate symmetric bilinear form is defined on $\bar\g$ 
by the following requirements: its extends the bilinear form on 
$L_\sigma\g$, $(k,d) = ({\bf 1},\delta) = 1$, the other
pairings of $k,{\bf 1},d,\delta$ are zero, and $k,{\bf 1},d,\delta$ 
are orthogonal to $L_\sigma\g$.  

Define a Lie subalgebra $\wh\g = L_\sigma\g\oplus \CC k \oplus \CC{\bf 1}
\subset \bar\g$. We also let $\bar\l := \l\oplus \on{Span}(k,{\bf 1},d,\delta)$.

We have $(L_\sigma\g)_0 = \g^\sigma$. 
Set $\bar\g_0 = \bar\l = (L_\sigma\g)_0 \oplus \on{Span}(k,{\bf 1},d,\delta)$, 
$\bar\g_\nu = (L_\sigma\g)_\nu$ if $\nu\in\Gamma - \{0\}$. 
Then $\bar\g = \oplus_{\nu\in\Gamma}\bar\g_\nu$.

\subsection{Critical cocycle}

If $\alpha\in\Gamma$, let $s_\alpha\in (L_\sigma\g)_\alpha \otimes 
(L_\sigma\g)_{-\alpha}$ be the element dual to the pairing $(-,-)$. 
Then $s_0 = t_\l$. Set $T = {1\over 2} m(s_0) + \sum_{\alpha >0}
m(\bar s_\alpha)$. Then $T$ belongs to the normal-order completion
$\wh{U(L_\sigma\g)}$ of $U(L_\sigma\g)$. 

We set 
$$
\gamma_\sigma := {1\over 2} \mu((\ell \otimes \id)(t_\u)) 
$$
(where $\mu$ denotes the Lie bracket). 

\begin{proposition} \label{prop:critical}
1) Denote by $Z(\l)$ the center of the Lie algebra $\l$, then 
$\gamma_\sigma$ belongs to $Z(\l)$. 

2) The derivation $u\mapsto [T,u]$ preserves $L_\sigma\g$, and we have 
for $u\in L_\sigma\g$
\begin{equation} \label{id:gamma:sigma}
[T,u] = -{1\over 2} (C_\g \otimes z {{\on{d}}\over{\on{d}z}})(u)
+ [\gamma_\sigma \otimes 1,u]. 
\end{equation}
\end{proposition}

\begin{remark} The critical level cocycle for $L_\sigma\g$ is defined as 
$(u,v) \mapsto -([T,u],v)$; so this cocycle is cohomologous to $\omega$. 
\end{remark}

\begin{remark} Let $\g$ be a simple, simply laced Lie algebra, 
$h$ be its Coxeter 
number, $\Delta_+ \subset \h^*$ a system of positive roots, and $\rho = {1\over
2}\sum_{\alpha\in\Delta_+} \alpha$. Equip $\h$ with its scalar product 
$(-,-)$ such that all roots have length $2$. 
Let $h_\rho \in \h$ be the element corresponding to $\rho$ 
(so $[h_\rho,e_\alpha] = (\rho,\alpha)e_\alpha$ for any root $\alpha$). 
Set $\sigma = \exp({{2\pi i}\over h} \on{ad}(h_\rho))$. Then 
$L_\sigma\g$ is isomorphic to $L\g$ with the principal gradation. 
In that case $\gamma_\sigma = 0$. Indeed, 
$$
\gamma_\sigma = {1\over {2}} \sum_{\alpha\in\Delta_+}
{{(\rho,\alpha)}\over h} [e_\alpha,f_\alpha]
+ (1 - {{(\rho,\alpha)}\over h}) [f_\alpha,e_\alpha]
= {1\over{h}} (\sum_{\alpha\in\Delta_+} (\rho,\alpha) h_\alpha) 
- h_\rho . 
$$
Now if $\beta\in\h$, $(\gamma_\sigma,\beta) 
= {1\over h} \sum_{\alpha\in\Delta_+} (\rho,\alpha)(\alpha,\beta) 
- (\rho,\beta)$ vanishes because of the identity $\sum_{\alpha\in\Delta_+}
\alpha\otimes \alpha = h (-,-)$; this identity holds up to scaling by
$W$-invariance, and the contraction of both sides yields $2\on{card}(\Delta_+)
= h \times \on{rank}(\g)$, which is Kostant's identity. 
\end{remark}

{\em Proof of Proposition \ref{prop:critical}.} 
Let us prove 1). $\gamma_\sigma$ clearly belongs to $\g_1 = \l$. 
On the other hand, $\ell$ commutes with the adjoint action of $\l$, 
and $t_\u$ is $\l$-invariant, so $(\ell\otimes \id)(t_\u)$ is $\l$-invariant.
Therefore $\gamma_\sigma$ commutes with $\l$. 

Let us prove 2). Let us set $s_\alpha = \sum_i e_{\alpha,i}(z) 
\otimes e_{-\alpha,i}(z)$. We first prove: 

\begin{lemma} \label{comp:T}
Let $\lambda\in\Gamma$ and $u(z)\in (L_\sigma\g)_\lambda$. If $\lambda\geq 0$, 
then 
$$
[T,u(z)] = -{1\over 2}\sum_{\alpha | 0 \leq \alpha < \lambda} 
[e_{\alpha,i}(z),[e_{-\alpha,i}(z),u(z)]],  
$$
and if $\lambda <0$, then 
$$
[T,u(z)] = {1\over 2}\sum_{\alpha | -\lambda \leq \alpha < 0} 
[e_{-\alpha,i}(z),[e_{\alpha,i}(z),u(z)]].   
$$
\end{lemma}

{\em Proof of Lemma \ref{comp:T}.} Assume that $\lambda \geq 0$, then 
$$
[T,u(z)] = 
{1\over 2} \sum_i [e_{0,i},u(z)] e_{0,i} + e_{0,i}[e_{0,i},u(z)]
+ \sum_{\alpha >0} \sum_i 
[e_{-\alpha,i}(z),u(z)]e_{\alpha,i}(z) + 
e_{-\alpha,i}(z) [e_{\alpha,i}(z),u(z)]. 
$$
Now if $\alpha\in\Gamma$, then 
\begin{equation} \label{useful:id}
[e_{-\alpha,i}(z),u(z)]\otimes 
e_{\alpha,i}(z) = - e_{\lambda - \alpha,i}(z) \otimes
[e_{\alpha-\lambda,i}(z),u(z)]. 
\end{equation}
So an infinity of cancellations take place, and we get 
\begin{align*} 
& 
[T,u(z)] = {1\over 2} 
\big( [e_{0,i},u(z)] e_{0,i} + e_{0,i}[e_{0,i},u(z)]\big) 
+ \sum_{0<\alpha\leq \lambda}
[e_{-\alpha,i}(z),u(z)] e_{\alpha,i}(z)
\\ & 
= {1\over 2} \sum_i [e_{0,i},u(z)] e_{0,i}
+ \sum_{0<\alpha<\lambda} [e_{-\alpha,i}(z),u(z)]e_{\alpha,i}(z)
+ {1\over 2} \sum_i [e_{-\lambda,i}(z),u(z)] e_{\lambda,i}(z). 
\\ & 
= {1\over 2} \sum_i [[e_{0,i},u(z)], e_{0,i}]
+ \sum_{0<\alpha<\lambda} [e_{-\alpha,i}(z),u(z)]e_{\alpha,i}(z), 
\end{align*}
where we have used (\ref{useful:id}) for $\alpha = \lambda$. 
Using again (\ref{useful:id}) for $0<\alpha<\lambda$, 
and using the change of variables $\alpha \mapsto \lambda - \alpha$, 
we get the first identity of Lemma \ref{comp:T}. The second identity is 
proved in the same way. 
\hfill \qed \medskip 

\begin{lemma} \label{for:critical}
Let $\alpha\in\Gamma$ be such that $0\leq \alpha <1$, 
and $u\in \g_{e(\alpha)}$, then 
$$
(\sum_{\lambda | 0\leq \lambda < \alpha} \on{ad}(m(t_\lambda))(u) 
- \ell \circ \on{ad}(m(t))(u) = -2 [\gamma_\sigma,u]. 
$$
Here $t_\lambda \in \g_{e(\lambda)} \otimes \g_{e(-\lambda)}$
is such that $\sum_{\lambda | 0\leq \lambda < 1} t_\lambda = t$. 
\end{lemma}

{\em Proof of Lemma \ref{for:critical}.}
For $\lambda\in\Gamma$, set $t_\lambda = \sum_i t_{\lambda,i}
\otimes t_{-\lambda,i}$, where $t_{\pm\lambda,i} \in \g_{e(\pm\lambda)}$. 
We have $t_\u = \sum_{\lambda | 0 < \lambda < 1} t_\lambda$, so 
$$
\gamma_\sigma = {1\over 2} \sum_{\lambda | 0<\lambda<1} \sum_i
[\ell(t_{\lambda,i}),t_{-\lambda,i}].  
$$
Then 
$$
 -2[\gamma_\sigma,u]
= - \sum_{\lambda | 0<\lambda<1} \sum_i 
[\ell(t_{\lambda,i}),[t_{-\lambda,i},u]]
+ [[\ell(t_{\lambda,i}),u],t_{-\lambda,i}]. 
$$
Now using $(\ell \otimes \id + \id \otimes \ell)(t_\lambda) = t_\lambda$, 
and the change of variables $\lambda \mapsto 1-\lambda$, we get 
$$
-2[\gamma_\sigma,u]
= \sum_{\lambda | 0<\lambda<1} \sum_i 
[(1-2\ell)(t_{\lambda,i}),[t_{-\lambda,i},u]]. 
$$
We split this sum as 
\begin{equation} \label{sum:of:3}
\sum_{\lambda | 0<\lambda<\alpha } \sum_i 
[(1-2\ell)(t_{\lambda,i}),[t_{-\lambda,i},u]] 
+ 
\sum_{\lambda | \alpha<\lambda<1 } \sum_i 
[(1-2\ell)(t_{\lambda,i}),[t_{-\lambda,i},u]] 
+ 
\sum_i 
[(1-2\ell)(t_{\alpha,i}),[t_{-\alpha,i},u]]. 
\end{equation}
Using the change of variables $\lambda \mapsto \alpha - \lambda$, 
we rewrite the first summand of (\ref{sum:of:3}) as 
$$
S_1 = {1\over 2} 
\sum_{\lambda | 0<\lambda<\alpha} 
\big( \sum_i 
[(1-2\ell)(t_{\lambda,i}),[t_{-\lambda,i},u]] 
+ \sum_i 
[(1-2\ell)(t_{\alpha-\lambda,i}),[t_{\lambda-\alpha,i},u]] 
\big) . 
$$
Now we have $\sum_i t_{\alpha - \lambda,i} \otimes [t_{\lambda-\alpha,i},u]
= \sum_i [u,t_{-\lambda,i}] \otimes t_{\lambda,i}$. Therefore 
$$
S_1 = {1\over 2} 
\sum_{\lambda | 0<\lambda<\alpha } \big( \sum_i 
[(1-2\ell)(t_{\lambda,i}),[t_{-\lambda,i},u]] 
+ \big[ (1-2\ell)([u,t_{-\lambda,i}]),t_{\lambda,i}\big] 
\big) , 
$$
which we rewrite as 
$$
S_1 = {1\over 2} 
\sum_{\lambda | 0<\lambda<\alpha } \big( \sum_i 
[(1-2\ell)(t_{\lambda,i}),[t_{-\lambda,i},u]] 
+ \big[ t_{\lambda,i},(1-2\ell)([t_{-\lambda,i},u])\big] 
\big) .  
$$
Now we have: if $t'\in \g_{e(\lambda)}$, $v\in \g_{e(\alpha-\lambda)}$, 
then $[(1-2\ell)(t'),v]+[t',(1-2\ell)(v)] = 2(1-\ell)([t',v])$. 
Therefore 
$$
S_1 = (1-\ell) \big( 
\sum_{0<\lambda<\alpha} \sum_i [t_{\lambda,i},[t_{-\lambda,i},u]] \big) . 
$$
In the same way, we use the change of variables 
$\lambda \mapsto 1+\alpha-\lambda$ and the identity: if $\alpha<\la<1$, 
$t'\in \g_{e(\lambda)}$, $v\in \g_{e(\alpha-\lambda)}$, 
then $[(1-2\ell)(t'),v]+[t',(1-2\ell)(v)] = -2\ell([t',v])$
to prove that the second summand of (\ref{sum:of:3})
is 
$$
S_2 = -\ell
\big( \sum_{\alpha<\lambda<1} 
\sum_i [t_{\lambda,i},[t_{-\lambda,i},u]] \big) . 
$$ 
Finally (\ref{sum:of:3}) is equal to 
\begin{equation} \label{sum:of:3:result}
\sum_{\lambda | 0<\lambda<\alpha}
\sum_i [t_{\lambda,i},[t_{-\lambda,i},u]]  
- \ell \big( 
\sum_{\lambda | 0<\lambda<1,\lambda\neq\alpha}
\sum_i [t_{\lambda,i},[t_{-\lambda,i},u]]  
\big) 
+ \sum_i 
[(1-2\ell)(t_{\alpha,i}),[t_{-\alpha,i},u]] 
\big) . 
\end{equation}
Now the last sum is (using the invariance of $t$, then the fact that $\ell$ 
is an $\l$-module endomorphism, then the invariance of $t$ again) 
\begin{align*}
& \sum_i [t_{\alpha,i},[t_{-\alpha,i},u]] + [- 2\ell([u,t_{0,i}]),t_{0,i}] 
=  
\sum_i [t_{\alpha,i},[t_{-\alpha,i},u]] -2\ell [ t_{0,i}, [t_{0,i},u]].  
\\ & = 
(1-\ell) \big( \sum_i [t_{0,i},[t_{0,i},u]]\big) - \ell \big(
\sum_i [t_{\alpha,i},[t_{-\alpha,i},u]] \big).   
\end{align*}
Finally (\ref{sum:of:3}) is equal to  
$$
\sum_{\lambda | 0\leq\lambda<1}
\sum_i [t_{\lambda,i},[t_{-\lambda,i},u]]  
- \ell \big( 
\sum_{\lambda | 0\leq\lambda<\alpha}
\sum_i [e_{\lambda,i},[e_{-\lambda,i},u]]  
\big) , 
$$
i.e., to the l.h.s. of the identity of Lemma \ref{for:critical}. 
\hfill \qed \medskip 

{\em End of proof of Proposition \ref{prop:critical}.}
Denote by $u\mapsto D_1(u)$, $u\mapsto D_2(u)$ both sides of 
(\ref{id:gamma:sigma}). Then $D_1,D_2$ are derivations of $L_\sigma\g$, 
such that 
$$
\forall u(z)\in L_\sigma\g, 
\quad D_i(z u(z)) = z D_i(u(z)) - {1\over 2} zC_\g(u(z)). 
$$
Here we view $L_\sigma\g$ as a module over $\CC[z,z^{-1}]$. 
Since $\oplus_{\alpha\in \Gamma | 0\leq \alpha <1} (L_\sigma\g)_\alpha$
generates $L_\sigma\g$ as a $\CC[z,z^{-1}]$-module, it suffices to prove
(\ref{id:gamma:sigma}) for $u(z)\in (L_\sigma\g)_\alpha$, $0\leq \alpha <1$. 

We have then $z v'(z) = \ell(v(z))$ for any $v(z) \in (L_\sigma\g)_\alpha$, 
In particular, $C_\g(u(z)) \in (L_\sigma\g)_\alpha$, therefore 
$$
(C_\g \otimes z{{\on{d}}\over{\on{d}z}})(u(z)) = 
\big( (\ell \circ C_\g) \otimes \id \big)(u(z))
$$
(in fact, one can show that $C_\g$ commutes with $\ell$). 

On the other hand, Lemma \ref{comp:T} implies that 
$$
[T,u(z)] = - {1\over 2} \sum_{\lambda | 0 \leq \lambda < \alpha} 
\on{ad}(m(t_\lambda))(u(z)). 
$$
Finally, we get 
\begin{align*}
& 
[T,u(z)] + {1\over 2} \big( C_\g \otimes z {{\on{d}}\over{\on{d}z}}\big)
(u(z))
\\ & = 
- {1\over 2} \big( \sum_{0\leq \lambda<\alpha} \on{ad}(m(t_\lambda))
\otimes \id\big) (u(z))
+ {1\over 2} \big(  \big( \ell \circ \sum_{0\leq \lambda <1} 
\on{ad}(m(t_\lambda))\big) \otimes \id\big) (u(z)) 
\\ & = [\gamma_\sigma\otimes 1,u(z)] \on{\ (by\ Lemma\ \ref{for:critical})}, 
\end{align*}
which proves (\ref{id:gamma:sigma}). 
\hfill \qed \medskip

\subsubsection{Infinite dimensional polarized Lie algebras}
\label{sect:53}

One checks that the theory of dynamical pseudotwists extends as follows. 

Let $\Gamma$ be a subset of $\CC$. We assume that 
$\ZZ \subset \Gamma$, $\Gamma$ is stable under the translations 
by elements of $\ZZ$, and $\Gamma/\ZZ$ is finite.  
 We set as before $\Gamma_\pm = \Gamma \cap \CC_\pm$. 

Let $\bar\g$ be a $\Gamma$-graded Lie algebra, $\bar\g = 
\oplus_{\nu\in\Gamma} \bar\g_\nu$. Here $\Gamma$-graded means that 
$[\bar\g_\nu,\bar\g_{\nu'}] \subset \bar\g_{\nu + \nu'}$ if $\nu + \nu'
\in\Gamma$, and equals 
$0$ otherwise.  

Set $\bar\l := \bar\g_0$. Assume that $\bar\l$ is a 
polarized Lie algebra 
$\bar\l = \bar\k \oplus \m_+ \oplus \m_-$. The Lie brackets 
induce linear maps $\m_+ \otimes \m_- \to \bar\k$ and 
$\bar\g_\nu \otimes \bar\g_{-\nu} \to \bar\k$ for 
$\nu\in \Gamma_+ - \{0\}$. 

We assume that (a) $\on{dim}(\m_+) = \on{dim}(\m_-)$ and 
$\on{dim}(\bar\g_\nu) = \on{dim}(\bar\g_{-\nu})$ for any 
$\nu\in \Gamma_+ - \{0\}$; 
we set $d'_\m = \on{dim}(\m_\pm)$ and $d'_\nu = \on{dim}(\bar\g_{\pm \nu})$, 
and (b) the corresponding determinants $\underline D_\m\in S^{d'_\m}(\bar\k)$ 
and $\underline D_\nu \in S^{d'_\nu}(\bar\k)$ are all nonzero. 

These linear maps can therefore be "inverted" and yield
$\rho_\m \in \m_+\otimes \m_- \otimes S^\cdot(\bar\k)[1/\underline D_\m]$, 
$\rho_\nu\in \bar\g_\nu \otimes \bar\g_{-\nu} \otimes S^\cdot(\bar\k)
[1/
\underline D_\nu]$, such that $\rho_\m + \sum_{\nu\in\Gamma_+ - \{0\}} 
\rho_\nu$ is a formal solution 
of the CDYBE. 

Let $D_\m,D_\nu$ be lifts in $U(\bar\k)$ of $\underline D_\m,\underline
D_\nu$. 
We denote by $\wh U$ the microlocalization of $U(\bar\k)$ w.r.t. all 
$D_\m,D_\nu$. 

Let us set $\bar\u_\pm = \m_\pm \oplus (\oplus_{\nu\in\Gamma_\pm - \{0\}} 
\bar\g_\nu)$, then 
$U(\bar\u_\pm)$ are $\Gamma_\pm$-graded algebras with
finite dimensional homogeneous parts. 
As in Section \ref{sect:K}, we can construct 
$$
K_\nu = \sum_i a_i^{(\nu)} \otimes b_i^{(\nu)} 
\otimes \ell_i^{(\nu)}
\in \big( U(\bar\u_+)_\nu \otimes U(\bar\u_-)_{-\nu} \big) \wh\otimes \wh U, 
$$ 
such that for any $x\in U(\bar\u_-)_{-\nu}$, $y\in U(\bar\u_+)_\nu$, we have 
$$
\sum_i H(xa_i^{(\nu)}) \ell_i^{(\nu)} H(b_i^{(\nu)}y) = H(xy). 
$$
Here $H$ is the Harish-Chandra map $U(\bar\g) \to U(\bar\k)$ corresponding to 
$\bar\g = \bar\k \oplus \bar\u_+ \oplus \bar\u_-$. 

Let us set $\bar\p_\pm = \bar\k \oplus \bar\u_\pm$. We also set 
$J_\nu = \sum_i a_i^{(\nu)} \otimes S(b_i^{(\nu)}) 
S(\ell_i^{(\nu)(2)})
\otimes S(\ell_i^{(\nu)(1)})$, then $J_\nu \in 
(U(\bar\u_+)_\nu \otimes U(\bar\p_-)_{-\nu}) \wh\otimes\wh U$. 
Moreover, $\wt J := \sum_{\nu\in\Gamma_+} J_\nu$ belongs to 
$\bar\oplus_{\nu\in\Gamma_+} (U(\bar\u_+)_\nu 
\otimes U(\bar\p_-)_{-\nu}) \wh\otimes \wh U$, and
satisfies the identity 
$$
\wt J^{12,3,4} \wt J^{2,3,4} = \wt J^{2,3,4} \wt J^{1,2,34}
$$
in $\bar\oplus_{\nu,\nu'\in\Gamma_+}
(U(\bar\u_+)_\nu \otimes U(\bar\g)_{\nu'-\nu}
\otimes U(\bar\p_-)_{-\nu'})\wh\otimes \wh U$. 
Here $\bar\oplus$ means the direct product.

\subsubsection{The ABRR equation in the infinite dimensional case}
\label{sect:inf:ABRR}

Assume that $\bar\g$ is equipped with a nondegenerate invariant 
pairing of degree $0$, $\langle-,-\rangle$, such that $\langle 
\bar\k,\m_\pm\rangle = \langle \m_\pm,\m_\pm\rangle = 0$. 
Let $s\in \m_+\otimes\m_-$, $t_\nu\in \bar\g_\nu\otimes\bar\g_{-\nu}$
and $t_{\bar\k}\in S^2(\bar\k)$ be dual to this pairing. 

Then the analogue of the normally-ordered ABRR equation (\ref{no:modif}) is
$$
(s^{1,2} + \sum_{\nu>0} t_\nu^{1,2}) \wt J 
= 
[ - m(\bar s)^{(2)} - \sum_{\nu>0} m(\bar t_\nu)^{(2)}
+ t_{\bar\k}^{2,3}, \wt J].  
$$
This is an identity in $U(\bar\g)^{\bar\otimes 2}\wh\otimes\wh U$. 

Moreover, the component $\wt J_0$ of $\wt J$
coincides with the twist $J_{\bar\k}^{\bar\g_0}$ corresponding to the 
nondegenerate polarized Lie algebra $\bar\g_0 = \bar\k\oplus \m_+ 
\oplus \m_-$. 

\subsubsection{} 

Let us return to the setup of Section \ref{4:6:1}. 
Assume that $\g^\sigma$ is polarized and nondegenerate, $\g^\sigma 
= \k\oplus \m_+ \oplus \m_-$. Then we set $\bar\k = \k\oplus 
\on{Span}(k,{\bf 1},d,\delta)$. Then 
$\bar\g\ = \oplus_{\nu\in\Gamma} \bar\g_\nu$, 
$\bar\g_0 = \bar\k \oplus \m_+ \oplus \m_-$ is an example of the situation of
Sections \ref{sect:53}, \ref{sect:inf:ABRR}. 

\subsubsection{The algebra $\wh U$}

We have $S^\cdot(\bar\k) = S^\cdot(\k) \otimes S^\cdot(\CC k\oplus \CC {\bf 1})
\otimes S^\cdot(\CC d \oplus \CC\delta)$. 

Let $\underline D^0_\m\in S^{d'_\m}(\k)$ be the determinant
corresponding to the nondegenerate Lie algebra $\l = \k\oplus \m_+
\oplus \m_-$, then $\underline D_\m = \underline D_\m^0 \otimes 1 \otimes 1$. 

Let us now describe $\underline D_\nu$ when $\nu \in \Gamma_+ - \{0\}$. 
Let us set
$d_\g := \on{dim}(\g)$. Then $\underline D_n$ is an element of 
$S^{d_\g}(\k\oplus \CC k) \subset S^{d_\g}(\bar\k)$, of the form 
$(nk)^{d_\g} + $ polynomial of partial degree $<d_\g$ in $k$. 
So $\underline D_\nu$ is nonzero. 

Let $D_\m,D_\nu$ be lifts of $\underline D_\m,\underline D_\nu$ in 
$U(\bar\k)$, let $\wh U$ be the corresponding 
microlocalization. 

Let $\wh U_\k$ (resp., $\wh U_{\k\oplus\CC d\oplus \CC\delta}$) 
be the microlocalization of $U(\k)$ (resp., 
$U(\k\oplus\CC d\oplus \CC\delta)$) w.r.t. $D_\u$. 
Then the form taken by the $\underline D_\nu$
allows us to embed $\wh U$ in $\wh U_{\k\oplus\CC d\oplus \CC\delta}((1/k))$. 
We will still denote by $\bar K$ the image of $\bar K$ in 
$U(\bar\g)^{\bar\otimes 2} \wh\otimes \wh U_{\k\oplus\CC d\oplus \CC\delta}
((1/k))$; 
actually, $\bar K$ belongs to 
$\big( U(\bar\u_+) \bar\otimes U(\wh\p_-) [1/k]\big) 
\wh\otimes \wh U_\k$, where $\wh\p_- = (\k\oplus \CC k)
\oplus\bar\u_-$.

\subsubsection{The affine ABRR equations}

We have a morphism $\pi_1 : U(\bar\u_+) \to U(\g)$, 
taking $a(x) \in (L_\sigma\g)_\alpha$ to $a(0)$ (recall that $\alpha > 0$),
$a\in \m_+$ to $a$, $k$ to $0$, ${\bf 1}$ to $1$. 

If $\gamma'\in Z(\l)$, we also have a morphism 
$\pi^{\gamma'}_z : U(\wh\p_-) \to U(\g) \otimes \CC[z^a,a\in\CC,
\on{log}(z)]$, 
taking $a f(x) \in (L_\sigma\g)_\alpha$ to $a \otimes f(x)$ 
(recall that $\alpha < 0$),
$a\in \m_+$ to $a \otimes 1$, $a\in \l$ to $\big( a - (\gamma',a) \big) 
\otimes 1$, $k$ to $0$, ${\bf 1}$ to $1$. 

Set $J(z) := (\pi_1 \otimes \pi^{-\gamma_\sigma}_z \otimes \id)(\wt J)$. Then 
$J(z)$ belongs to $U(\g)^{\otimes 2} \wh\otimes \wh U_\k((1/k))
[[z^a,a\in\CC_+]][\log(z)]$. Here $A[[z^a,a\in\CC_+]] = 
A[z^a,a\in {\bf D}][[z]]$ (${\bf D}$ is defined in Section \ref{sect:52}). 

The ABRR equation is an equality in $U(\wh\g)^{\bar\otimes 2}
\wh\otimes \wh U_\k((1/k))$. We have $(\pi_z^{-\gamma_\sigma} \otimes \pi_1)
(\sum_{\nu>0} t_\nu) = - X(z)$. Moreover, 
$$
(\pi^{-\gamma_\sigma}_z\otimes \pi_1\otimes \id)
([\sum_{\nu>0} m(\bar t_\nu)^{(2)},\wt J]) 
= [-{1\over 2}m(t_\l)^{(2)} ,J(z)], 
$$
so the image of ABRR by 
$\pi_z^{-\gamma_\sigma} \otimes \pi_1 \otimes \id$ is 
\begin{equation} \label{***}
- k^{(3)} z {{dJ}\over{dz}}
= \big(X(z) - s\big)^{1,2} J(z) + [t_\k^{2,3} + {1\over 2}
m(t_\k)^{(2)} - \gamma^{(2)},J(z)],    
\end{equation}
where $\gamma = \mu(s)$. 
Moreover, the constant coefficient of $J(z)$ is 
$J_0 = J_\k^\l$. These
conditions determine the series $J(z)$ uniquely. 

Set $\kappa =  - 1/k^{(3)}$, then $G(z)$ is a solution of (\ref{basic:eqn})
iff $G(z) z^{- \kappa (t_\k^{2,3} + {1\over 2} m(t_\k)^{(2)} - 
\gamma^{(2)})}$
is a solution of (\ref{***}). 

Let us write $J_\k^\g(z) = J(z)$. 
We then have $\Psi_\k^\g = \lim_{z\to 1^-} (1-z)^{-\kappa t^{1,2}} 
J_\k^\g(z)$. The composition formula implies that 
$J_\k^\g(z) = \eta (J_\l^\g(z)) J_\k^\l$, therefore 
$\Psi_\k^\g = \eta(\Psi_\l^\g) J_\k^\l$. 

We can now reprove Proposition \ref{prop:59}, 1) as follows: 
we have 
$$
\eta(\Psi_\l^\g) = \Psi_\k^\g (J_\k^\l)^{-1} = 
\lim_{z\to 1^-} (1-z)^{-\kappa t^{1,2}} J_\k^\g(z) (J_\k^\l)^{-1},
$$ 
now $g(z) := J_\k^\g(z) (J_\k^\l)^{-1}$ is such that $g(0) = 1$, 
and (thanks to the ABRR equation for $J_\k^\l$) 
$$
z {{\on{d}g}\over {\on{d}z}} = 
\kappa (X(z) - s)^{1,2} g(z) + \kappa 
[t_\k^{2,3} + {1\over 2} m(t_\k)^{(2)} 
- \gamma^{(2)},g(z)]. 
$$
It follows that $\eta(\Psi_\l^\g)$ is the renormalized holonomy from $0$ to $1$
of (\ref{31}), i.e., $\Psi_{\k,\l,\g}$. 

The theory of infinite dimensional ABRR equations also 
underlies the systems (\ref{syst:1},\ref{syst:2}) and 
(\ref{syst:3},\ref{syst:4}). Indeed, set 
$J(z,u) := (\pi_{1} \otimes \pi^{-\gamma_\sigma}_{z^{-1}} \otimes 
\pi^{-\gamma_\sigma}_{u^{-1}} \otimes \id)
(\wt J^{[2]})$, then the multicomponent ABRR implies that 
$J(z,u)$ satisfies the equations
\begin{align*}
z{{\partial J}\over{\partial z}} = & \kappa \Big(
X(z)^{1,2} + X(z/u)^{3,2} - s^{1,2} - s^{3,2} + t_\k^{2,4} 
+ {1\over 2} m(t_\k)^{(2)} - \gamma^{(2)}  \Big) J(z,u)
\\ & - \kappa J(z,u) (t_\k^{2,3} + t_\k^{2,4} + {1\over 2}m(t_\k)^{(2)}
- \gamma^{(2)}), 
\end{align*}
\begin{align*} 
u{{\partial J}\over{\partial u}} = & \kappa \Big( 
X(u)^{1,3} + X(u/z)^{2,3} - s^{1,3} - s^{2,3} + t_\k^{3,4} + {1\over 2}
m(t_\k)^{(3)} - \gamma^{(3)} 
\Big) J(z,u)
\\ & - \kappa J(z,u) (t_\k^{3,4} + {1\over 2}m(t_\k)^{(3)} - \gamma^{(3)}), 
\end{align*}
so that $J(z,u)$ satisfies this system iff it has the form 
$G(z,u) z^{\kappa(t_\k^{2,3} + t_\k^{2,4} + {1\over 2}m(t_\k)^{(2)} 
- \gamma^{(2)})}
u^{\kappa(t_\k^{3,4} + {1\over 2}m(t_\k)^{(3)} - \gamma^{(3)})}$, 
where $G(z,u)$ is a solution of (\ref{syst:3},\ref{syst:4}). 

The compatibility of the systems (\ref{syst:1},\ref{syst:2}) and 
(\ref{syst:3},\ref{syst:4}) is the image by 
$\pi_{1} \otimes \pi^{-\gamma_\sigma}_{z^{-1}} \otimes 
\pi^{-\gamma_\sigma}_{u^{-1}} \otimes \id$ 
of the compatibility relations for the multicomponent 
ABRR equations (Proposition \ref{ABRR:compat}).  

\begin{remark}
If $\gamma'$ is an element of $Z(\l)$, and 
$J_{\gamma'}(z)$ is the analogue of $J(z)$, 
where $-\gamma_\sigma$ is replaced by $\gamma'-\gamma_\sigma$, 
then $\on{lim}_{z\to 1^-} (1-z)^{-\kappa t^{1,2}} J_{\gamma'}(z)$
coincides with $\Psi_{\kappa,\gamma'}$, as defined in Section 
\ref{sect:shifts}. The ABRR arguments of this section can be modified to 
provide other proofs of the statements of Section \ref{sect:shifts}. 
\end{remark}

\section{Cayley $r$-matrices}

Proposition \ref{lemma:0:3} follows from Proposition \ref{prop:0:4}. Let us 
prove this proposition. It will be enough to treat the case $c = 1$. 
We set $\rho_C := \rho_{C,1}$. 

Set $\rho_\l(\lambda) := (f(\on{ad}(\lambda^\vee))\otimes \id)(t_\l)$, 
$$
\rho_\u(\lambda) := i \pi \Big( 
{{ (C + \id) +  e^{-2\pi i \on{ad}(\lambda^\vee)} (C-\id)} \over
{(C + \id) - e^{-2\pi i \on{ad}(\lambda^\vee)} (C - \id) }} 
\otimes \id \Big) (t_\u). 
$$
Then $\rho_C = \rho_\l + \rho_\u$. We want to prove that 
$\on{CYB}(\rho_\l + \rho_\u) - \on{Alt}(\on{d}(\rho_\l + \rho_\u))
= -\pi^2 Z$. We have $Z = Z_\l + Z_{\l,\u} + Z_\u$, where 
$Z_\l = [t_\l^{1,2},t_\l^{2,3}]$, $Z_{\l,\u} = \on{Alt}
([t_\l^{1,2},t_\i^{2,3}])$, $Z_\u = (\id\otimes p_\u\otimes
\id)([t_\u^{1,2},t_\u^{2,3}])$, where $p_\u : \g\to\u$ is the projection on 
$\u$ parallel to $\l$. 

Applying \cite{AM1} to the quadratic algebra $(\l,t_\l)$, we 
already have $\on{CYB}(\rho_\l) - \on{Alt}(\on{d}\rho_\l) = -\pi^2 Z_\l$. 
It remains to prove that 
$$
\on{CYB}(\rho_\l,\rho_\u) + \on{CYB}(\rho_\u) - \on{Alt}(\on{d}\rho_\u)
= -\pi^2 (Z_{\l,\u} + Z_\u). 
$$
Both sides of this equality belong to $\wedge^3(\g) = 
\oplus_{\alpha = 0}^3 \wedge^\alpha(\l) \otimes \wedge^{3-\alpha}(\u)$. 
Since the equality already
holds when projecting it on the components $\alpha = 3$ and $\alpha = 2$, 
it remains to prove its projection on the components $\alpha = 0$
and $\alpha = 1$. 

The projection on the component $\alpha = 0$ is the equality
\begin{equation} \label{alpha=0}
\on{Alt} \circ p_\u^{(2)} \big( 
[\rho_\u^{1,2}(\lambda),\rho_\u^{2,3}(\lambda)] 
+ \pi^2[t_\u^{1,2},t_\u^{2,3}]\big) = 0, 
\end{equation}
and the projection on the component $\alpha = 1$ is 
\begin{equation} \label{alpha=1}
p_\l^{(1)}([\rho_\u^{1,2}(\lambda),\rho_\u^{1,3}(\lambda)])
+ [\rho_\l^{1,23}(\lambda),\rho_\u^{2,3}(\lambda)]
- (\on{d}\rho_\u(\lambda))^{2,3,1} + \pi^2 [t_\l^{1,2},t_\u^{2,3}]
= 0. 
\end{equation}
To prove (\ref{alpha=0}), we apply to it $\big( 
(C + \id) - e^{- 2\pi i \on{ad}(\lambda^\vee)} (C- \id)\big)^{\otimes 3}$. 
We get 
\begin{align*} 
& 
\Big( 
\big( (C + \id) - e^{-2\pi i \on{ad}(\lambda^\vee)}
(C - \id)\big) \otimes 
\big(
(C + \id) + e^{-2\pi i \on{ad}(\lambda^\vee)}
(C - \id)
\big)^{\otimes 2}
+ \text{cyclic\ permutation}
\\ &  
+ \big( (C + \id) - e^{-2\pi i \on{ad}(\lambda^\vee)}
(C - \id)\big)^{\otimes 3} \Big) (Z_\u) = 0, 
\end{align*}
i.e., 
$4 \big( (C + \id)^{\otimes 3} - (e^{-2\pi i \on{ad}(\lambda^\vee)}
(C - \id))^{\otimes 3} \big)(Z_\u) = 0$, which follows from 
$(C + \id)^{\otimes 3}(Z_\u) = (C - \id)^{\otimes 3}(Z_\u)$, which follows 
from the assumptions on $C$. 

Let us now prove (\ref{alpha=1}). Let us apply to it 
$\id\otimes \big( (C + \id) - e^{-2\pi i \on{ad}(\lambda^\vee)}
(C - \id) \big)^{\otimes 2}$, we get 
\begin{align} \label{47}
& 
\pi^2 
\Big( \id\otimes \big( (C + \id) + e^{-2\pi i\on{ad}(\lambda^\vee)} 
(C - \id)\big)^{\otimes 2} \Big) ([t_\l^{1,2},t_\u^{2,3}])
\\ & \nonumber
+ \pi^2 
\Big( \id\otimes \big( (C + \id) - e^{-2\pi i\on{ad}(\lambda^\vee)} 
(C - \id)\big)^{\otimes 2} \Big) ([t_\l^{1,2},t_\u^{2,3}])
\\ & \nonumber
- i \pi (e^{2\pi i \on{ad}(\lambda^\vee)})^{(2)}
\big( (C + \id) - e^{-2\pi i \on{ad}(\lambda^\vee)}(C - \id)\big)^{(2)} 
\on{d} \big( 
\big( (C + \id) + e^{-2\pi i \on{ad}(\lambda^\vee)}(C - \id)\big)^{(2)}
(t_\u^{2,3}) \big) 
\\ & \nonumber
+ i\pi (e^{2\pi i \on{ad}(\lambda^\vee)})^{(2)}
\big( (C + \id) + e^{-2\pi i \on{ad}(\lambda^\vee)}(C - \id)\big)^{(2)} 
\on{d} \big( 
\big( (C + \id) - e^{-2\pi i \on{ad}(\lambda^\vee)}(C - \id)\big)^{(2)}
(t_\u^{2,3}) \big)   
\\ & \nonumber
+ i \pi \Big( \id \otimes 
\big( (C + \id) - e^{-2\pi i \on{ad}(\lambda^\vee)}(C - \id)\big)^{\otimes 2}
\Big) ([\rho_\l^{1,23}(\lambda),
\Big( 
{ {C + \id + e^{-2\pi i \on{ad}(\lambda^\vee)}(C - \id)} \over
{C + \id - e^{-2\pi i\on{ad}(\lambda^\vee)}(C - \id)}}
\Big)^{(2)}(t_\u^{2,3})]) = 0,  
\end{align}
i.e., 
\begin{align*}
& 2 \pi^2 (\id - C^2)^{(2)}([t_\l^{1,2},t_\u^{2,3}])
- i \pi (e^{2\pi i\on{ad}(\lambda^\vee)})^{(2)}
(C^2 - \id)^{(2)} \on{d} \Big( (e^{-2\pi i\on{ad}(\lambda^\vee)})^{(2)}
(t_\u^{2,3})\Big) 
\\ & + \text{last\ line\ of\ }(\ref{47}) = 0. 
\end{align*}
Now 
\begin{align*}
& \text{last\ line\ of\ (\ref{47})}
\\ & 
= i\pi 
\Big( \id \otimes 
\big( C + \id + e^{-2\pi i \on{ad}(\lambda^\vee)}(C - \id)\big) 
\otimes 
\big( C + \id - e^{-2\pi i \on{ad}(\lambda^\vee)}(C - \id)\big) 
\Big) ([\rho_\l^{1,3}(\lambda),t_\u^{2,3}])
\\ & 
- i\pi 
\Big( \id \otimes 
\big( C + \id - e^{-2\pi i \on{ad}(\lambda^\vee)}(C - \id)\big) 
\otimes 
\big( C + \id + e^{-2\pi i \on{ad}(\lambda^\vee)}(C - \id)\big) 
\Big) ([\rho_\l^{1,2}(\lambda),t_\u^{2,3}])
\\ & 
= - 2i\pi 
\Big( \id \otimes (C + \id)^{\otimes 2}
- \id \otimes 
\big( e^{-2\pi i\on{ad}(\lambda^\vee)}(C - \id) \big)^{\otimes 2}\Big) 
([\rho_\l^{1,2}(\lambda),t_\u^{2,3}])
\\ & 
= -2 i \pi (\id - C^2)^{(2)} \big( 1 - \id \otimes (e^{- 2 \pi i
\on{ad}(\lambda^\vee)})^{\otimes 2}\big)
([\rho_\l^{1,2}(\lambda),t_\u^{2,3}]) = 0, 
\end{align*}
which follows from the CDYBE identity for the Alekseev-Meinrenken 
$r$-matrix for $(\g,t)$, restricted to $\lambda\in\l^*$ and projected on 
$\l\otimes \wedge^2(\u)$. 
\hfill \qed \medskip

\section{Quantization of homogeneous spaces}
\label{sect:homog}

In this section, we show that the (pseudo)twists constructed in 
Sections \ref{sect:K}, \ref{sect:pseudo} and \ref{sect:twists} enable us 
to quantize (quasi)Poisson structures on homogeneous spaces.  

\subsection{Quantization of coadjoint orbits} \label{prev:sect}

Let $\g = \l\oplus \u$ be a Lie algebra with a splitting; we assume that $\g$
is nondegenerate.
Let $G$ be the formal group with Lie algebra 
$\g$, and $L\subset G$ the subgroup corresponding to 
$\l$ (with suitable restrictions, the following constructions
may be extended to other categories, like algebraic or complex Lie 
groups). 

Let $D_0 : \l^* \to \CC$ be the determinant corresponding to 
$\g = \l\oplus \u$. The dynamical $r$-matrix, $r_\l^\g(\lambda)$, 
enables one to define a Poisson structure on 
$\big( \l^* - D_0^{-1}(0)\big) \times G$. The group $L$
acts on this space by Poisson automorphisms, and the moment map 
is the projection on the first factor. According to P.\ Xu, a dynamical twist
quantizing $r_\l^\g(\lambda)$ yields a quantization of this Poisson space. 

Let $\chi\in\l^*$ be a character of $\l$. Then $L\subset \on{Stab}(\chi)$, 
where $\on{Stab}(\chi) := \{g\in G | \on{Ad}^*(g)(\chi) = \chi\}$. Let 
$O_\chi \subset \g^*$ be the coadjoint orbit of $\chi$. Then we have a natural
$G$-map $G/ L \to G / \on{Stab}(\chi) = O_\chi$, taking the class of $g$ 
to $\Ad^*(g)(\chi)$. 

Let us further assume that $D_0(\chi) \neq 0$. Then we have $L =
\on{Stab}(\chi)$, so $G/L \to O_\chi$ is an isomorphism. 
The assumptions on $\chi$ imply that 
$r_\l^\g(\chi)\in\wedge^2(\g)$ is well-defined and $\l$-invariant. 
This implies that the bivector ${\bold R}(r_\l^\g(\chi))$ on $G$ (${\bold R}$
denotes the translation from the right) descends to $G/L$, and that it equips
$G/L$ with a Poisson structure. Moreover, the map $G/L \to O_\chi$
is Poisson (this is an observation of J.-H.\ Lu). 

According to Remark \ref{rem:char}, we may extend $\hbar^{-1}\chi$
to a character of $\wh{U}_\l$, since $D_0(\chi)\neq 0$, therefore 
$(\id\otimes \id \otimes \hbar^{-1}\chi)(J)$
is well-defined; it coincides with $J(\chi)$ as defined in 
Section \ref{sect:quant}, 
and belongs to $U(\g)^{\otimes 2}[[\hbar]]$. It satisfies 
$J(\chi)^{12,3} J^{1,2}(\chi + l^{(3)}) = J(\chi)^{1,23} J(\chi)^{2,3}$. 
Here $J^{1,2}(\chi + l^{(3)})$ has the usual meaning, and its
action on $f\otimes g \otimes h$ is the same as that of $J(\chi)^{1,2}$
if $h$ is $\l$-invariant. This relation implies that one can define a 
$G$-equivariant star-product on $G/L$ by the formula 
$f*g=m({\bold R}(J(\chi))(f\otimes g))$ (where 
$\bold R$ stands for right translations), quantizing the Poisson structure on
$G/L$ described above. 

By virtue of the results of Section \ref{sect:K}, these considerations 
allow us to get equivariant star-products in all nondegenerate polarized cases. 
In the polarized quadratic case, $J(\chi)$ satisfies the equation 
(derived from ABRR) 
$$
s^{1,2} J(\chi) = 
[ \big( {1\over 2} m(t_\l) + \hbar^{-1} t_\l^\vee(\chi) - \gamma\big)^{(2)}, 
J(\chi)].  
$$

In the reductive case, this quantization (which yields an explicit 
equivariant star-product for all semisimple coadjoint orbits) has been 
obtained by Alekseev-Lachowska (\cite{AL}) and Donin-Mudrov (\cite{DM}).

\subsection{Quantization of Poisson homogeneous spaces}
\label{sect:6:2} 

Let $\g = \l\oplus \u$ be a Lie algebra with a splitting, such that 
$\g$ is nondegenerate. We assume that $\g$ is quadratic, i.e., 
we have $t\in S^2(\g)^\g$. Let $c$ be a complex number. 

The dynamical $r$-matrix $\rho_c(\lambda)$ (see Corollary \ref{cor:rat:trigo})
can be used to equip $U \times G$ with the structure of a quasi-Poisson 
homogeneous space under the pair $(G,t)$ (\cite{AKM}). 
Here $U = \{\lambda\in \l^* | \on{ad}(\lambda^\vee)$ 
has no eigenvalue of the form $n/c$, $n$ a nonzero integer, and 
$D_0(\lambda)\neq 0\}$. 
Let $\chi\in \l^*$ be a character of $\l$, and assume that $\chi\in U$. 
Then $G/L$ is equipped with a quasi-Poisson homogeneous 
space structure under $(G,t)$, given by the bivector $\Pi_\chi=
{\bold R}(\rho_c(\chi))$.
These quasi-Poisson structures may be viewed as trigonometric versions of
the Poisson structures of Section \ref{prev:sect}. 

A quantization $J$ of $\rho_c(\lambda)$ gives rise to a quantization of 
$U_{\on{formal}} \times G$ (in the sense of \cite{EE}, Section
4.5), where $U_{\on{formal}}$ is the intersection of 
$U$ with the formal neighborhood of $0 \in \l^*$. 

Moreover, if $J(\lambda)$ is regular at $\chi$, 
then it can be used to construct a quantization of this quasi-Poisson 
space, according to the formula 
$$
f*g=m({\bold R}(J(\chi))(f\otimes g)).
$$
(Recall that we do not know a quantization of $\rho_c$ in the 
nonpolarized case, even if $c=0$.)

Assume now that $\g$ is polarized, i.e., $\u = \u_+ \oplus \u_-$
and $t = t_\l + s + s^{2,1}$, with $t_\l \in S^2(\l)$, 
$s\in \u_+ \otimes \u_-$.
Then $J(\lambda)$ has been constructed in Section \ref{sect:pseudo}, 
and if we take $\Phi = \Phi^{\on{KZ}}$, $J(\lambda)$ is regular on 
an explicit neighborhood of $0$ (see \cite{EE}). This yields a 
quantization of $U_{\on{formal}} \times G/L$ and of 
$(G/L,\Pi_\chi)$ for characters $\chi$ in this 
neighborhood. 

Let $\varrho\in \g^{\otimes 2}$ be a quasitriangular structure on $\g$, 
i.e., $\varrho + \varrho^{2,1} = t$ and  $\on{CYB}(\varrho) = 0$. 
Let $(G,({\bold R} - {\bold L})(\varrho))$ be the 
corresponding Poisson-Lie group. We have 
$\on{CYB}(\pi i c(\varrho - \varrho^{2,1})) = - \pi^2 c^2 Z$
(equality in $\wedge^3(\g)$) and 
$\on{CYB}(\rho_c(\chi)) = -\pi^2 c^2 Z$ (equality in 
$\wedge^3(\g/\l)$), where $Z = [t^{1,2},t^{2,3}]$. 
Therefore $G/L$, equipped with the 
Poisson bivector $\Pi_{\varrho,\chi} := -{\bold L}(\pi i c (\varrho -
\varrho^{2,1})) + {\bold R}(\rho_c(\chi))$, is a Poisson homogeneous 
space under $(G,({\bold R} - {\bold L})(\varrho))$. 
(Here ${\bold L}$ stands for left translations.)

A quantization of the Poisson homogeneous space $(G/L,\Pi_{\varrho,\chi})$
may be obtained as follows. According to \cite{EK} (in the reductive 
case, \cite{Dr2,ESS}), there exists a pseudotwist $J_{\on{EK}} 
\in U(\g)^{\otimes
2}[[\hbar]]$ quantizing $\varrho$, i.e., $J_{\on{EK}}^{12,3} 
J_{\on{EK}}^{1,2} =  \Phi_\kappa^{\on{KZ}}(t^{1,2},t^{2,3})^{-1} 
J_{\on{EK}}^{1,23} J_{\on{EK}}^{2,3}$. 
Then the star-product on $G/L$ is defined 
by the formula 
$$
f*g=m({\bold R}(J(\chi)){\bold L}(J_{\on{EK}}^{-1})(f\otimes g)). 
$$
This quantization is equivariant with respect to the quantum group 
$U(\g)^{J_{\on{EK}} }$ ($U(\g)$ twisted by $J_{\on{EK}}$).  

In the case when $G,L$ are reductive, the homogeneous spaces we 
considered include generic dressing orbits of $G$, and we get 
their quantization equivariant under the quantum group $U_q(\g)$.  
A different way of quantizing such Poisson 
(and quasi-Poisson) homogeneous spaces 
was proposed in \cite{DGS}.


\subsection{Quantization of Poisson homogeneous spaces
corresponding to an automorphism} \label{sect:autom}

Let us assume that $(\g,t\in S^2(\g)^\g)$ is a quadratic Lie 
algebra, equipped with $\sigma\in \on{Aut}(\g,t)$. We set 
$\l := \g^\sigma$ and assume that $\sigma - \id$ is invertible on 
$\g/\g^\sigma$. 

As above, the dynamical $r$-matrix $\rho_{\sigma,c}(\lambda)$ can be used 
to equip $G/L$ with a structure of a 
quasi-Poisson homogeneous space 
of the group $(G,-\pi^2c^2 Z)$ (where $Z = [t^{1,2},t^{2,3}]$). 
Namely, the quasi-Poisson bivector on 
$G/L$ is given by the formula $\Pi=
{\bold R}(\rho_{\sigma,c}(0))$.  
We will set $c = 1/(2\pi i)$, therefore 
$$
\Pi = {\bold R}\big(
{1\over 2}({{\sigma + \id}\over{\sigma - \id}}\otimes \id)
(t_\u) \big)
$$
The dynamical pseudotwist 
$\Psi_\kappa$ provides 
a quantization of this quasi-Poisson structure. 
Namely, set  $\Psi_\kappa(0) := (\id\otimes \id\otimes \eps)(\Psi_\kappa)$. 
The non-associative star-product on $G/L$
(which is associative in the representation category of Drinfeld's
quasi-Hopf algebra) is given by the formula
$$
f*g=m({\bold R}(\Psi_\kappa(0))(f\otimes g)).
$$

Let $\varrho\in \g^{\otimes 2}$ be a quasitriangular structure on $\g$, 
i.e., $\varrho + \varrho^{2,1} = t$ and  $\on{CYB}(\varrho)=0$. Let 
$(G,({\bold R} - {\bold L})(\varrho))$ be the corresponding 
Poisson-Lie group.
Since $\on{CYB}({{\varrho - \varrho^{2,1}}\over 2}) = Z/4$ 
(in $\wedge^3(\g)$) and $\on{CYB}(\rho_{\sigma,c}(0)) = Z/4$ 
in $\wedge^3(\g/\l)$, we have a 
Poisson homogeneous space $G/L$ under 
$(G,({\bold R} - {\bold L})(\varrho))$, 
with Poisson bivector 
\begin{equation} \label{bivect}
\Pi=-{\bold L}({{\varrho - \varrho^{2,1}}\over 2})
+ {\bold R}({1\over 2} ({{\sigma+\id}\over{\sigma-\id}}
\otimes \id)(t_\u)).
\end{equation}

The above construction yields a star-product quantization 
of this Poisson homogeneous structure.  
Namely, the star-product on $G$ is defined 
by the formula 
$$
f*g=m({\bold R}(\Psi_\kappa(0)){\bold L}(J^{-1})(f\otimes g)),
$$
where $J$ is a pseudotwist quantizing $\varrho$ (e.g., $J = J_{\on{EK}}$).
As in Section \ref{sect:6:2}, this quantization is equivariant under 
the quantum group $U(\g)^{J}$. 

\subsection{Relation to the De Concini homogeneous spaces}

Recall that according to Drinfeld
\cite{Dr1}, if $G$ is a Poisson-Lie group and $L$ is a subgroup, 
then Poisson homogeneous space structures on $G/L$ correspond 
to Lagrangian Lie subalgebras $\h\subset D(\g)$ of the double of $\g$ 
such that $\g\cap \h=\l$.

C. De Concini explained to us the following construction
of Poisson homogeneous spaces. 
Let $\g$ be a factorizable quasitriangular Lie bialgebra. 
This means that $\g$ is a Lie algebra, $\varrho\in\g^{\otimes 2}$
is such that $\on{CYB}(\varrho) = 0$, and $t:= \varrho + \varrho^{2,1}
\in S^2(\g)^\g$ is  nondegenerate. Assume also that 
$\sigma\in \on{Aut}(\g,t)$. 
Then $D(\g)$ is isomorphic to $\g\oplus\g$, with bilinear form  
given by $\langle (x_1,x_2), (y_1,y_2)\rangle = 
\langle x_1,y_1 \rangle - \langle x_2,y_2\rangle$. 
The graph $\h$ of $\sigma$ is a Lagrangian subalgebra of $\g\oplus \g$, 
which induces a Poisson homogeneous space structure on $G/L 
= G/G^{\sigma}$. 

\begin{theorem} \label{thm:DC} 
The construction of Section \ref{sect:autom} yields quantizations of 
all the De Concini homogeneous spaces, such that $\sigma$ is 
invertible on $\g/\g^\sigma$. 
\end{theorem}

{\em Proof.} The Drinfeld subalgebra $\h\subset D(\g)$ corresponding to 
a Poisson homogeneous space $(G/L,\Pi)$ is defined as 
$$
\h = \{(x,\xi)\in \g\oplus \g^* | \xi\in\l^* \text{\ and\ }
x = (\xi\otimes \id)(\Pi(0)) \text{\ modulo\ }\l\}, 
$$
where $\Pi(0)\in\wedge^2(\g/\l)$ is the value at origin of $\Pi$. 
In the case of the Poisson structure (\ref{bivect}),  
$\Pi(0)$ is equal to the class of $P$ in $\wedge^2(\g/\l)$, where 
$$
P = {{\varrho - \varrho^{2,1}}\over 2} + {1\over 2} (\id \otimes 
{{\sigma + \id}\over {\sigma - \id}})(t_\u), 
$$
therefore $\h$ is the image 
of the linear map $\l\oplus \u^* \to \g\oplus \g^*$, $(x,0)\mapsto 
(x,x)$, $(0,\xi) \mapsto (x(\xi),\xi)$. Here $x(\xi) = (\xi\otimes
\id)(P)$. 

Let us set  $L(\alpha) := (\id\otimes \alpha)(\varrho)$ and $R(\alpha) := 
(\alpha\otimes \id)(\varrho)$, for $\alpha\in\g^*$. 
If $\xi\in\u^*$, then $(\xi\otimes \id)(t_\u) = (\xi\otimes \id)(t)$, 
therefore $(L+R)(\xi)\in\u$ for any $\xi\in\u^*$. If follows that 
$x(\xi) = {1\over 2}(R-L)(\xi) + {1\over 2} {{\sigma + \id}\over
{\sigma - \id}} \circ (L+R)(\xi)$.  

According to \cite{RS}, the isomorphism $D(\g) \to \g\oplus \g$
is given by $(a,0) \mapsto (a,a)$ and $(0,\alpha) \mapsto (-R(\alpha),
L(\alpha))$. 

Let us view $\h$ as a subalgebra of $\g\oplus \g$ using this isomorphism. 
Then $\h$ is the image of the linear map $\l \oplus \u^* \to \g\oplus \g$, 
$(x,0) \mapsto (x,x)$ and 
$$
(\xi,0) \mapsto \big( x(\xi) - R(\xi), x(\xi) + L(\xi)\big) 
= \big( 
{{\id}\over{\sigma - \id}} \circ (L+R)(\xi), 
{{\sigma }\over{\sigma - \id}}  \circ (L+R)(\xi)
\big).  
$$
Since the image of $\u^*\to \u$, $\xi\mapsto (L+R)(\xi)$ is exactly $\u$, 
and $\sigma - \id$ is invertible on $\u$, $\h$ is the image of 
$\l\oplus \u\to \g\oplus \g$, $(x,0)\mapsto (x,x)$ and $(0,y)\mapsto 
(y,\sigma(y))$, i.e., $\{(z,\sigma(z)) | z\in\g\}$, i.e., $\h$
is the graph of $\sigma$. 
\hfill \qed \medskip 

\begin{remark} It is useful to describe $\Psi_\kappa(0)$ directly, 
since this is the only information about $\Psi_\kappa$ needed in 
Theorem \ref{thm:DC}: $\Psi_\kappa(0)$ is the renormalized holonomy 
from $0$ to $1$ of the differential equation 
${{dG}\over{dz}} = \kappa (X(z)^{1,2} + {1\over 2} m(t_\l)^{(2)}) 
G(z)$.   
\end{remark}

\subsection{Poisson homogeneous spaces corresponding to a Cayley 
endomorphism}

Assume that $\g = \l\oplus \u$ is a Lie algebra with a splitting and a 
factorizable structure, such that $t = t_\l\oplus t_\u$, $t_\x\in S^2(\x)$
for $\x = \l,\u$, and $C\in\End_\l(\u)$ is a Cayley endomorphism,
such that $(C\otimes \id + \id\otimes C)(t_\u) = 0$.  
Then $\h\subset \g\oplus \g$ defined by 
$$
\h := \l^{\on{diag}} \oplus \{(x,y) \in \u\times \u| (C+\id)(x)
= (C-\id)(y)\}
$$ is a Lagrangian subalgebra, generalizing De Concini's 
subalgebras (here $\l^{\on{diag}} = \{(x,x) | x\in \l\}$). 
It gives rise to a Poisson homogeneous structure 
on $G/L$. A quantization of the dynamical $r$-matrices of Proposition 
\ref{prop:0:4} should lead to a quantization of these homogeneous spaces. 

An example of this situation is: $\g$ is semisimple, with 
Cartan decomposition $\g = \t\oplus \n_+ \oplus \n_-$ and the corresponding 
standard $r$-matrix $\varrho$; $w$ is a Weyl group element; $\l = \t$, 
$\u = \n_+ \oplus \n_-$; $C$ has eigenvalues $\pm 1$ on $w(\n_\mp)$
($\wt w(\n_\pm)$ is independent on the choice of a Tits lift of $\wt w$, 
so we denote it $w(\n_\pm)$). Then $\h = \{(h+x_+,h+x_-) | 
h\in \t, x_\pm \in w(\n_\pm)\}$. The corresponding Poisson 
homogeneous structure on $G/T$ is given by $-{\bold L}(\varrho) 
+ {\bold R}(w^{\otimes 2}(\varrho))$
(if $V$ is a $\g$-module and $v\in V$ is a vector of weight $0$, 
$\wt w(v)$ is independent on the choice of $\wt w$ and is denote by 
$w(v)$). In this case a quantization 
can be obtained using the formula
$f * g := m({\bold L}(J^{-1}) {\bold R}(w^{\otimes 2}(J))(f\otimes g))$, 
where $J$ is a pseudotwist quantizing $\varrho$. 
Another quantization is given by 
$f*g = m_\hbar ({\bold R}(J_w)(f\otimes g))$, where $f,g\in U_\hbar(\g)^*$
and $m_\hbar$ is the product in $U_\hbar(\g)^*$, and 
$J_w\in U_\hbar(\g)^{\otimes 2}$
is such that $T_w^{\otimes 2} \circ \Delta \circ T_w^{-1} = 
\on{Ad}(J_w) \circ \Delta$, where $\Delta$ is the coproduct of 
the Drinfeld-Jimbo quantum group $U_\hbar(\g)$
and $T_w$ is a Lusztig-Soibelman automorphism corresponding to $w$.

The following fact implies the equivalence of both quantizations. 

\begin{proposition}
$w^{\otimes 2}(J)$ and $J_w J$ are gauge-equivalent
pseudotwists quantizing $w^{\otimes 2}(\varrho)$. 
\end{proposition}

{\em Proof.} Recall the construction of $J$: let $\Phi$ be 
a Drinfeld associator, $\Phi_\g$ its specialization to $(\g,t)$. 
Then $J$ is a series $J_\Phi(\varrho)$, such that 
$\wt d(J) := (J^{2,3}J^{1,23})^{-1} J^{1,2} J^{12,3} = \Phi_\g$. 
Since $J_w$ satisfies the twist equation, we have 
$\wt d(J_w J) = \Phi_\g$. On the other hand, 
$\wt d(w^{\otimes 2}(J)) = w^{\otimes 3}(\Phi_\g) = 
\Phi_\g$. So  $w^{\otimes 2}(J)$ and $J_w J$ are 
pseudotwists quantizing $w^{\otimes 2}(\varrho)$. 

Let us now show that they are gauge-equivalent. 
Let us still denote by $T_w$ the automorphism of 
$U(\g)[[\hbar]]$ obtained by transporting $T_w$
by the isomorphism $U(\g)[[\hbar]] \simeq U_\hbar(\g)$. 
The reduction mod $\hbar$ of $T_w$ is a Tits lift of $w$, 
which we denote by $\wt w$. 
Then $T_w \circ \wt w^{-1}$ is an automorphism of $U(\g)[[\hbar]]$, 
whose reduction modulo $\hbar$ is the identity, hence is inner
(as $\g$ is semisimple). Let $u\in U(\g)[[\hbar]]$, 
$u = 1 + O(\hbar)$, be such that $T_w \circ \wt w^{-1} = \on{Ad}(u)$. 

Let $\Delta_0$ be the undeformed coproduct of $U(\g)[[\hbar]]$. 
Set $J_1 = J_w J$, $J_2 = w^{\otimes 2}(J)$. 
We have $\on{Ad}(J_1) \circ \Delta_0 = T_w^{\otimes 2} \circ \Delta_0 \circ
T_w^{-1}$,  
$\on{Ad}(J_2) \circ \Delta_0 = \wt w^{\otimes 2} \circ \Delta_0 \circ
\wt w^{-1}$, so 
$$
\on{Ad}(J_1) \circ \Delta_0 = \on{Ad}(u)^{\otimes 2}
\circ \on{Ad}(J_2) \circ \Delta_0 \circ \on{Ad}(u)^{-1}, 
$$
therefore $J_1 = u^{\otimes 2} J_2 \Delta_0(u)^{-1} \xi$, 
where $\xi\in U(\g)^{\otimes 2}[[\hbar]]$ has the form $1 + O(\hbar)$
and is $\g$-invariant. 

We now prove inductively that $\xi = \on{exp}(\on{d}(\eta))$, 
where $\eta \in \hbar (U(\g)^{\otimes 2})^\g[[\hbar]]$
and $\on{d}(\eta) = \Delta_0(\eta) - \eta \otimes 1 - 1 \otimes \eta$. 
Assume that we have found $\eta_1,\ldots,\eta_{n-1}$ in $U(\g)^\g$, 
such that $\xi = \on{exp}(\on{d}(\hbar \eta_1 + \cdots + \hbar^{n-1}
\eta_{n-1}))(1 + O(\hbar^n))$. 
Then set $u' := u \on{exp}(-(\hbar \eta_1 + \cdots + \hbar^{n-1}\eta_{n-1}))$. 
We get $J_1 = (u')^{\otimes 2} J_2 \Delta_0(u')^{-1} \xi'$, 
with $\xi' \in (U(\g)^{\otimes 2})^\g[[\hbar]]$ has the form $1 + O(\hbar^n)$. 
Let $\xi_n\in (U(\g)^{\otimes 2})^\g$ be such that $\xi = 1 + \hbar^n \xi_n +
O(\hbar^{n+1})$, then since $\wt d(J_1) = \wt d(u^{\prime \otimes 2}
J_2 \Delta_0(u')^{-1})$, we get $\on{d}(\xi_n) = \xi_n^{12,3}
- \xi_n^{1,23} - \xi_n^{2,3} + \xi_n^{1,2}= 0$. Since 
the cohomology group involved is $\wedge^2(\g)$
and since $\wedge^2(\g)^\g = 0$, we get $\xi_n = \on{d}(\eta_n)$, 
with $\eta_n \in U(\g)^\g$. This proves the induction step. 

Now let $\eta := \sum_{i\geq 1} \hbar^i \eta_i$. We get 
$J_1 = (u e^{-\eta})^{\otimes 2} J_2 \Delta_0(u e^{-\eta})^{-1}$, 
therefore $J_1$ and $J_2$ are gauge-equivalent. 
\hfill \qed 


\medskip 


\begin{thebibliography}{EnGH}

\bibitem [AF]{AF}
A. Alekseev, L. Faddeev, {\it $(T\sp *G)\sb t$: a toy model for 
conformal field theory,} Commun. Math. Phys. \textbf{141} (1991), no. 2, 
413--422.

\bibitem [AKM]{AKM}
A. Alekseev, Y. Kosmann-Schwarzbach, E. Meinrenken, {\it Quasi-Poisson 
manifolds,} Canad. J. Math. \textbf{54} (2002), no. 1, 3--29.

\bibitem[AL]{AL}
A. Alekseev, A. Lachowska,  {\it Invariant $*$-products on coadjoint orbits 
and the Shapovalov pairing,} preprint math.QA/0308100 (2003). 

\bibitem [AM1]{AM1}
A. Alekseev, E. Meinrenken,
{\it The non-commutative Weil algebra,} Invent. Math. \textbf{139} (2000), 
135--172.

\bibitem [AM2]{AM2}
A. Alekseev, E. Meinrenken,
{\it Clifford algebras and the classical dynamical Yang-Baxter equation,} 
preprint math.RT/0209347, Math. Res. Lett. \textbf{10} (2003), no. 2-3, 
253--268. 

\bibitem [ABRR] {ABRR} D. Arnaudon, E. Buffenoir, E. Ragoucy, and P. Roche,
{\it Universal solutions of quantum dynamical Yang-Baxter equations,}
q-alg/9712037, Lett. Math. Phys. \textbf{44} (1998), no. 3, 201--214.

\bibitem[BS]{BS}
A. Borel, J. De Siebenthal, {\it Les sous-groupes ferm\'es de rang 
maximum des groupes de Lie clos,} Comment. Math. Helv. 
\textbf{23} (1949), 200--221.

\bibitem[DGS]{DGS}
J. Donin, D. Gurevich, S. Shnider, {\it Double quantization on some orbits 
in the coadjoint representations of simple Lie groups,} Commun. Math. Phys. 
\textbf{204} (1999), no. 1, 39--60.

\bibitem[DM]{DM} J. Donin, A. Mudrov,
{\it  Dynamical Yang-Baxter equation and quantum vector bundles,}
preprint math.QA/0306028 (2003).

\bibitem[Dr1]{Dr1}
V. Drinfeld, {\it On Poisson homogeneous 
spaces of Poisson-Lie groups,} Teoret. Mat. Fiz. \textbf{95} (1993), 
no. 2, 226--227; translation in Theoret. and
Math. Phys. \textbf{95} (1993), no. 2, 524--525. 

\bibitem [Dr2]{Dr2}
V. Drinfeld, {\it Quasi-Hopf algebras,} Leningrad Math. J. \textbf{1}
(1990), 1419--1457.

\bibitem[EE]{EE} 
B. Enriquez, P. Etingof, {\it Quantization of Alekseev-Meinrenken 
dynamical r-matrices,} preprint math.QA/0302067 (2003).  
 
\bibitem[EFK]{EFK} P. Etingof, I. Frenkel, A. Kirillov,
{\it Lectures on Representation Theory and Knizhnik-Zamolodchikov 
Equations,} AMS, 1998.

\bibitem [EK]{EK}
P. Etingof, D. Kazhdan, 
{\it Quantization of Lie bialgebras, I, II,}  
Selecta Math. (N.S.) \textbf{2} (1996), no. 1, 1--41;
\textbf{4} (1998), no. 2, 213--231. 

\bibitem[ESS]{ESS}
P. Etingof, T. Schedler, O. Schiffmann, {\it Explicit 
quantization of dynamical r-matrices,} 
J. Amer. Math. Soc. \textbf{13} (2000), 595--609.

\bibitem[ES1]{ES1} P. Etingof, O. Schiffmann,
{\it Lectures on the Dynamical Yang-Baxter Equations,} 
math.QA/9908064, in "Quantum groups and Lie theory (Durham, 1999)", 
89--129,  London Math. Soc. Lecture Note Series, \textbf{290} (2001), 
Cambridge Univ. Press, Cambridge. 

\bibitem [ES2]{ES2}
P. Etingof, O. Schiffmann, 
{\it  On the moduli space of classical dynamical 
$r$-matrices,} Math. Res. Lett. \textbf{8} (2001), 
no. 1-2, 157--170.

\bibitem [EV1]{EV1}
P. Etingof, A. Varchenko, {\it Geometry and classification of
solutions of the classical dynamical Yang-Baxter equation,}
Commun. Math. Phys. \textbf{192} (1998), 77--120.

\bibitem [EV2]{EV2}
P. Etingof, A. Varchenko, {\it Exchange dynamical quantum groups,}  
Commun. Math. Phys. \textbf{205} (1999), no. 1, 19--52.

\bibitem [EV3]{EV3}
P. Etingof, A. Varchenko, {\it Dynamical Weyl groups and applications,}
Adv. Math., \textbf{167} (2002), 74--127. 

\bibitem [Fad]{Fad}
L. Faddeev,  {\it On the exchange matrix 
of the WZNW model}, Commun. Math. Phys. \textbf{132} (1990), 
131--138.

\bibitem [FGP]{FGP}
L. Feh\'er, A. G\'abor and P. Pusztai, {\it On dynamical $r$-matrices 
obtained from Dirac reduction and their generalizations to affine
Lie algebras,} math-ph/0105047, J. Phys. A \textbf{34}
(2001), no. 36, 7335--7348. 
  
\bibitem[RS]{RS}
N. Reshetikhin, M. Semenov-Tian-Shansky, 
{\it Quantum $R$-matrices and factorization problems,} 
J. Geom. Phys. \textbf{5} (1988), no. 4,
533--550.
  
\bibitem[S]{S}
O. Schiffmann, {\it On classification of dynamical r-matrices,}  
Math. Res. Lett. \textbf{5} (1998), no. 1-2, 13--30.

\bibitem[Spr]{Spr}
T. Springer, {\it Microlocalisation alg\'ebrique,} 
S\'em. alg\`ebre Dubreil-Malliavin, Lecture Notes in Mathematics 
\textbf{1146} (1983), 299-316, Springer-Verlag.  

\bibitem [Xu]{Xu}
P. Xu, {\it Quantum dynamical Yang-Baxter equation over a nonabelian 
base,} Commun. Math. Phys. \textbf{226} (2002), no. 3, 475--495.

\end{thebibliography}
\end{document}